\input amstex
\def\fbox#1#2{\vbox to 0pt{\hrule\hbox{\hsize=#1\vrule\kern2pt
  \vbox{\kern2pt#2\kern2pt}\kern2pt\vrule}\hrule}}
\mag=\magstep1
\documentstyle{amsppt}

\topmatter
\title 3-fold log flips according to 
V. V. Shokurov
\endtitle
\author Hiromichi Takagi$^*$ 
\endauthor
\rightheadtext {Review of $3$-fold log flips}
\leftheadtext {Hiromichi Takagi}
\address 
Department of Mathematical Science,
University of Tokyo,
Komaba,
Meguro-ku,
Tokyo 153,
Japan\endaddress
\email
htakagi\@ms.u-tokyo.ac.jp
\endemail
\thanks 
\newline
$*$ Research Fellow of the Japan Society for the Promotion of Science
\endthanks
\abstract  
We review $\S 8$ of V. V. Shokurov's paper "$3$-fold log flips".
\endabstract
\pagewidth{14.5cm}
\hcorrection{-1cm}
\vcorrection{-0.5cm}
\toc
\head 0. Introduction \endhead
\head 1. Preliminaries \endhead
\head 2. Auxiliary flips and back tracking method \endhead
\head 3. Set up and dividing into cases \endhead
\head 4. Treatment of case (a) and (b1) \endhead
\head 5. Preliminaries for case (b2) and (c) \endhead
\head 6. Treatment of case (b2, $\delta$) \endhead
\head 7. Treatment of case (c) \endhead
\endtoc

\endtopmatter
\def \pa#1#2#3{#1^*(K_{#2}+#3)}
\def \pb#1#2#3#4{#1^*(K_{#2}+#3+#4)}
\def \pc#1#2#3#4#5{#1^*(K_{#2}+#3+#4+#5)}
\def \lda#1#2{K_{#1}+#2}
\def \ldb#1#2#3{K_{#1}+#2+#3}
\def \ldc#1#2#3#4{K_{#1}+#2+#3+#4}
\def \-#1{{#1}^{-1}}
\def \supp{\text{Supp} \ }
\def \llc{\text{LLC} \ }
\def \clc#1{\text{CLC}_{#1} \ }
\def \diff#1{\text{Diff}_{#1} \ }
\def \t#1{\tilde{#1}}
\def \+n#1{{#1}^{+\nu}}
\def \o#1{\overline{#1}}

\document

\head 0. Introduction 
\endhead
In his paper '$3$-fold log flips', V. V. Shokurov proved the 
existence of the flips for flipping contractions from
log $3$-folds with only Kawamata log terminal singularities. 
First we review his strategy after $3$ definition needed for explanation.
 
\definition{Definition (0.1) (Special flip) 
(See [S,Definition 6.1] or [FA, Definition 18.8])}
Let $X$ be a $\Bbb Q$-factorial $n$-fold and $S$ a integral Weil divisor.
Let $f:X\to Y$ be a small contraction with $\rho(X/Y)=1$.
We say that $f$ is a special flipping contraction if
\roster
\item $\lda{X}{S}$ is LT; 
\item $\lda{X}{S}$ and any component of $S$ are $f$-negative. 
\endroster
We call the flip of such an $f$ the special flip.
\enddefinition

\definition{Definition (0.2) (Special termination)}
We say that the special termination holds if
for any sequence of flips
$$(X_1, B_1) \dashrightarrow\dots
\dashrightarrow(X_k, B_k)\dots ,$$
where $(X_i, B_i)$ is a LC pair,
there is $k_0 \in \Bbb N$ (depending on a sequence)
such that flipping curves on $X_{k_0}$ are disjoint from
$\llcorner B_{k_0} \lrcorner$.
\enddefinition

\definition{Definition (0.3) ($\o{\Cal S}_n ^0 (\text{local})$)}
We denote by $\o{\Cal S}_n ^0 (\text{local})$ the set 
of sequences $(b_1, b_2, \dots b_m)$ such that there is a 
$\Bbb Q$-factorial variety $X$ of dimension at most $n$, a subset $Z\subset X$
and a divisor $B_0 + \sum b_i B_i$ ($B_0 \not = \phi$ is reduced but possibly
reducible) such that every $B_i$ intersects $Z$, $Z\subset B_0$,
$K_X + B_0$ is PLT and $K_X + B_0 + \sum b_i B_i$ is maximally log canonical
near $Z$.
\enddefinition

Shokurov proved the following Reduction Theorem:

\proclaim{Theorem (0.4) (Reduction Theorem)}
Let $X$ be a $n$-dimensional normal variety and $B$ a boundary
such that $K_X + B$ is KLT.
Let $f: X\to Y$ be a small contraction such that $-(K_X + B)$ is $f$-ample.
Then the flip of $f$ exists if the following three hold:

\roster
\item (Existence of Special flips) 

there exists the flip for any special flipping contraction
of dimension $n$;

\item (Special termination)

The special termination holds in dimension $n$;

\item (A.C.C for $\o{\Cal S}_n ^0 (\text{local})$)

The ascending chain condition holds for the set
$\o{\Cal S}_n ^0 (\text{local})$.
\endroster
\endproclaim

\demo{Proof}
See [S,Theorem 6.4] or [FA, \S 18].
\qed
\enddemo
 
In dimension $3$, the special termination holds by
[S, Theorem 4.1] or [FA, Theorem 7.1].
Furthermore the A.C.C holds for $\o{\Cal S}_3 ^0 (\text{local})$
by [S, Chicago Lemma 4.9]
or [FA,Proposition 18.19 and Corollary 18.25.1].
Hence in dimension $3$, it suffices to prove the existence
of the special flip.
So we consider only the special flipping contraction in the below and
we use the notation in Definition 0.1 for explanation.

First we can easily prove the existence in case $S$ is reducible
(See [FA, Proposition 21.2].)
Hence we may assume that $S$ is irreducible in the below.
To flip in this case, 
Shokurov introduced the concept of the complement and
divided special flipping contractions into cases using it:

\definition{Definition (0.5) (Complement) (See [S, \S 5] or [FA, \S 19])}
Let $X$ be a normal variety and $B$ a subboundary on $X$.
Let $S$ be the smallest effective Weil divisor on $X$ such that
$\llcorner B-S \lrcorner \leq 0$, and let $B_0 := B-S$.
An $n$-complement of $K_X + B$ is a divisor
$\o{B} \in |-nK_X -nS - \llcorner (n+1)B_0 \lrcorner|$ such that
$K_X + B^+$ is LC, where
$B^+ := S + \frac 1n (\llcorner (n+1)B_0 \lrcorner + \o{B})$.
We say $K_X + B$ is $n$-complemented if an $n$-complement exists.
\enddefinition

The definition is arranged so that the vanishing theorem can be applied.
But in the explanation below, we consider only the case
$B$ is a reduced boundary. So we restate the definition for this case.

\definition{Definition (0.5') (Complement)}
Let $X$ be a normal variety and $S$ a reduced boundary.
An $n$-complement of $K_X + S$ is a divisor
$\o{B} \in |-nK_X -nS|$ such that
$K_X + S+B$ is LC, where
$B :=\frac 1n \o{B}$.
We say $K_X + B$ is $n$-complemented if an $n$-complement exists.
\enddefinition

Then the following theorem holds:

\proclaim{Theorem (0.6)}
Let $f: X\to Y$ be a $3$-dimensional special flipping contraction
such that $S$ is irreducible.
Then $K_X +S$ is $n$-complemented, where $n \in \{ 1,2,3,4,6 \}$.
\endproclaim

\demo{Proof}
See [S, Theorem 5.12] or [FA, Theorem 19.6 and Theorem 19.8].
\qed
\enddemo

We call the minimum of $n$ as in Theorem 0.6 the index of $f$.

In case the index is $1$, it is easy to show the existence of the flip.
(See [S, Proposition 6.8] or [FA, Proposition 21.4].)

In case the index is $2$, 
if $K_X + S+ B$ is LT, it is easy to prove the existence of the
flip.
(See below Proposition (2.2) which we quot from [S] and [FA].)

In case $n\geq 2$, Shokurov used the back tracking method
(see [FA, \S 6 and \S 21] or the section 2 of this paper
for the back tracking method) and proved the following:

\proclaim{Theorem (0.7)}
If index $2$ special flips exist, all special flips exist.
\endproclaim

\demo{Proof}
See [S, Reduction 7.6] or [FA, Theorem 21.10].
\qed
\enddemo

Hence the existence of index $2$ special flips are left
which we explain in this paper.

By some minor reductions, it suffices to consider the flipping contraction
as in Set up (3.0).
After some preparations in the sections 1 and 2, we start the proof
of the existence of the special flip. 
The case (a) or (b) in Set up (3.0) is relatively easy to treat.
(But the author had trouble understanding the proof of [S, Proposition 8.6]
or [FA, Lemma 22.6] which is Theorem (4.0) in this paper.)
The case (c) is very hard to treat and
we have to divide it into finer cases.
The proof in the section 7 is almost taken from [S, Proposition-Reduction 8.8]
but we separate Shokurov's argument into Lemmas
(Lemma 7.6 $\sim$ 7.9 are principal ones).
Futhermore the author cannot understand the inductive argument in
[S, (8.8.4) and (8.8.5)] and hence find the alternative argument
(see the proof of Theorem 7.1).

\definition{Acknowledgement}
I express my hearty thanks to Professor Yujiro Kawamata
to recommend me to write down this paper.
This is the extended version of my talk in the seminar of the doctor course
(1997) (I explained the contents of the section 0$\sim$6).
I am grateful to the member of the seminar,
Professor Yujiro Kawamata,
Professor Keiji Oguiso, Doctor Nobuyoshi Takahashi, Mr. Tatsuhiro Minagawa
and Mr. Kimikazu Kato for listening patiently to my talk and 
giving me useful comments.
\enddefinition
 
\definition{Notation and Conventions}
\roster
\item In this paper, we will work over $\Bbb C$, the complex number field;
\item The $\Bbb P^1$-bundle 
$\Bbb P(\Cal O_{\Bbb P^1}\oplus \Cal O_{\Bbb P^1}(-a))$ over $\Bbb P^1$, 
a Hirzebruch surface of degree $a$, is denoted by $\Bbb F_a$. 
Its unique negative section is denoted by $C_0$ 
and a ruling is denoted by $f$. 
The projective cone obtained from $\Bbb F_a$ by the contraction of $C_0$ 
is denoted by $\Bbb F_{a,0}$. 
A generating line on $\Bbb F_{a,0}$ is denoted by $l$;
\item we use the following abbreviations:

ODP for ordinary double point,
LT for log terminal, PLT for purely log terminal, 
KLT for kawamata log terminal, DLT for divisorial log terminal,
LC for log canonical, NC for normal crossing and 
MRS for minimal resolution. 
\endroster
\enddefinition  

\head 1. Preliminaries 
\endhead 
In this section, we gather results used in the proof of the main theorem.
The reader had better skip this section and refer
later when the results in this section are quoted. 

We start by easy lemmas about log surfaces.

\proclaim{Definition and Proposition (1.0)}
Let $(S, \lda{S}{B}, s)$ be a germ of LC surface singularity.
Then there is a unique extraction $\mu:T\to S$ 
satisfying the following conditions:
\roster
\item
$\lda{T}{B_T}=\pa{\mu}{S}{B}$ and $\lda{T}{B_T}$ is LT, where
all exceptional curves appear in $B_T$ with coefficients $1$;
\item if there is a $\mu$-exceptional $(-1)$-curve, then 
$\lda{T}{B_T}$ become non LT after contracting this. 
\endroster
We call this $\mu:T\to S$ the minimal LT model (abbreviated MLT) 
of $\lda{S}{B}$.
\endproclaim
\comment
\demo{Proof}
It is well-known except the uniqueness of $\mu$.
If $B$ has two reduced boundaries, $\mu$ is MRS
of $s$ so the uniqueness is clear.
Assume that $B$ has at most one reduced boundary.
Let $\mu':T'\to S$ be another extraction satisfying $(1)$ and $(2)$ and 
$U$ a common partial resolution over $S$ such that there is no curve
exceptional for both $U\to T$ and $U\to T'$.
Let $E$ be a curve exceptional for $U\to T$.
Then by assumption, the image of $E$ on $T'$ is a $\mu'$-exceptional
curve so $E$ is log crepant for $\lda{T}{B_T}$.
If $B$ has no reduced component, 
$\lda{T}{B_T}$ is PLT, so there is no such an $E$.
By interchanging $T$ with $T'$, the same argument shows that there is no
curve exceptional for $U\to T'$. Hence we are done.
Assume that $B$ has one reduced boundary $B'$.
Let $t$ be the (unique) intersection of 
$\-{\mu}B'$ and a $\mu$-exceptional divisor. 
Outside $t$, $\lda{T}{B_T}$ is PLT so
$E$ is contracted to $t$.
\qed
\enddemo
\endcomment

\proclaim{Proposition (1.1) (Description of LC surface 1)}
Let $(S, \lda{S}{B}, s)$ be a germ of LC surface singularity.
Then $\llcorner B \lrcorner$ consists of at most $2$ components.
Furthermore if $\llcorner B \lrcorner$ consists of $2$ components,
then $B=\llcorner B \lrcorner$.
\endproclaim

\demo{Proof}
See [Ka1, Theorem 9.6].
\qed
\enddemo
We will repeatedly use this Proposition for $3$-dimensional situations 
as follows:

\proclaim{Corollary (1.2)}
Let $(X, \lda{X}{D}, x)$ be a $3$-dimensional LC germ.
Then every irreducible curve on $X$ is contained in 
at most two reduced components of $D$.
Furthermore if an irreducible curve is contained in two reduced components of 
$D$, the curve is contained in no other component of $D$. 
\endproclaim

\demo{Proof}
Take a general hyperplane section of $X$ and apply (1.1).
\qed
\enddemo

\proclaim{Proposition (1.3) (Description of LC surface 2)}
Let $(S, \lda{S}{B}, s)$ be a germ of LC surface singularity.
Then for an exceptional curve $F$ of MRS of $s$,
we have $d(F, B) \leq 1$.
\endproclaim

\demo{Proof}
Easy.
\qed
\enddemo 

\proclaim{Proposition (1.4) (Description of LC surface 3)}
Let $(S, \lda{S}{B}, s)$ be a germ of LC surface singularity
such that $2(\lda{S}{B}) \sim 0$.
Then its MLT is described as follow:
FIGURE (I)

Assume that $(S,  \lda{S}{B}, s)$ is not LT.
Then all exceptional curves of MRS except the end $(-2)$-curves of the dual
graph are log crepant for $\lda{S}{B}$.
\endproclaim
\demo{Proof}
Easy.
\qed
\enddemo

\definition{Definition (1.5) (LLC, CLC and PLC) (See [Ka4, Definition 1.3]}
Let $X$ be a normal variety and $D$ be a $\Bbb Q$-divisor such that
$\lda{X}{D}$.

Then we define 

$\text{CLC} \ (\lda{X}{D})$ to be the set of subvarieties $W$ on $X$
such that 
there is a birational morphism from a normal variety
$\mu: Y \to X$ and a prime divisor $E \subset Y$  
with the discrepancy coefficient $e \leq -1$
such that $\mu(E)=W$.

We denote the set of $m$-dimensional elements in $\text{CLC} \ (\lda{X}{D})$
by \break
$\clc{m}(\lda{X}{D})$. 

We define $\llc(\lda{X}{D})$ to be 
the union of subvarieties in $\text{CLC} \ (\lda{X}{D})$.  

We define $\text{PLC} \ (\lda{X}{D})$ to be 
the set of discrete valuations $E$ in $\Bbb C(X)$
such that there is a birational morphism from a normal variety 
$\mu: Y \to X$  
such that $E$ is realized as a prime divisor on $Y$ 
with the discrepancy coefficient $e \leq -1$.

\enddefinition

The next two results descript $\llc$.

\proclaim{Theorem (1.6) (Connectedness Lemma (a))}
Let $X$ and $Y$ be normal varieties and $f: X\to Y$
a proper surjective morphism with only connected fibers.
Let $D$ be a $\Bbb Q$-divisor on $X$ 
such that $\lda{X}{D}$ is $\Bbb Q$-Cartier.
Write $D:=\sum \limits ^{n}_{i=1} d_iD_i$, where $D_i$ is an
irreducible component of $D$.
  
Assume the following:
\roster
\item if $d_i<0$, $\text{codim} \ f(D_i) \geq 2$;
\item $-(\lda{X}{D})$ is $f$-nef and $f$-big;  
\endroster
 
Then for any point $y \in Y$,
$\llc (\lda{X}{D}) \cap \-{f}(y)$ is connected. 
\endproclaim

\demo{Proof}
See [S, 5.7] or [FA, 17.4 Theorem] and also [Ka, Theorem 1.4].
\qed
\enddemo
 
\proclaim{Theorem (1.7) (Connectedness Lemma (b))}
Let $(S, \Theta)$ be a proper LC surface.
Assume that $\lda{S}{\Theta} \equiv 0$.
Then 
one of the following holds:
\roster
\item $\llcorner \Theta \lrcorner$ is irreducible and there is no
point contained in $\clc{0}(\lda{S}{\Theta})$ on $\llcorner \Theta \lrcorner$;
\item $\llcorner \Theta \lrcorner$ is connected and there is a
point contained in $\clc{0}(\lda{S}{\Theta})$ on every component
of $\llcorner \Theta \lrcorner$. Furthermore 
$\llc (\lda{S}{\Theta}) = \llcorner \Theta \lrcorner$;
\item $S$ is a generically $\Bbb P^1$-bundle over a smooth rational or 
elliptic curve and 
$\llcorner \Theta \lrcorner=C_1\cup C_2$, where $C_1$ and $C_2$ are
two disjoint sections of above generically $\Bbb P^1$-bundle.
There is no point contained in $\clc{0}(\lda{S}{\Theta})$ on $C_1$ or $C_2$.
Furthermore $\llc (\lda{S}{\Theta}) = \llcorner \Theta \lrcorner$.
\endroster
\endproclaim 

\demo{Proof}
See [S, Theorem 6.9] or [FA, 12.3.1 Proposition].
\qed
\enddemo
  
\proclaim{Lemma (1.8) (Ampleness)}
Let $f:S\dashrightarrow T$ be a birational map of normal projective surfaces
and $D$ an effective ample divisor on $S$ such that $D_T$ is irreducible.
($D_T$ is the log birational transform, i.e., all the blown up curve
are contained in it with multiplicity 1.) Then $D_T$ is numerically ample.
\endproclaim
\demo{Proof}
See [S, 8.10 Lemma].
\qed
\enddemo

The following theorem is slightly weaker than [S, Lemma 8.9].
But in our application, it is sufficient.

\proclaim{Theorem (1.9) ($\clc{0}$ for log surface)}
Let $(S, B)$ be a projective LC surface
and $C_1$ and $C_2$ connected contractible curves.
(Caution: we allow that $C_i$ is empty. The contractibility
doesn't mean $C_1$ and $C_2$ are contractible simultaneously.)
We assume the following:
\roster 
\item"(i)" $2(K_S+B)\sim 0$; 
\item"(ii)" $\llcorner B \lrcorner=B_1+B_2$, where $B_i$ is an
irreducible curve and $B_1^2>0$ and $B_2$ is ample;
\item"(iii)" $B_1\cap C_1=\phi$;
\item "(iv)" $P:=B_1\cap B_2$ is one point and on $B_1$,
$P$ is the unique point contained in $\clc{0}(\lda{S}{B})$;
\item "(v)" all the component of $C_2$ pass through $P$.
\endroster
(See FIGURE (II).)
Then $C_1\cap B_2\not \in \clc{0}(K_S+B)$.
\endproclaim

\demo{Proof}
First note that if $C_1\cap B_2$ (resp. $C_2\cap (B_1\cup B_2)$)
is not empty,
$C_1 \cap B_2$ (resp. $C_2 \cap (B_1 \cup B_2)$)
is one point by (1.6).

Let $\mu:T\to S$ be the MRS of $S$.

\proclaim{Claim 1}
Let $C'_2$ be any component of $C_2$.
Then 
\roster
\item $\-{\mu}C'_2$ is a $(-1)$-curve;
\item $P$ is singular point of $S$.
Outside $P$, there are only Du Val singularities on $C'_2$;
\item $C'_2\cap \supp\{ B \}= \phi$.

\endroster
\endproclaim
\demo{Proof}
By the assumption of (v), $C'_2$ passes through $P$ and
there are two reduced boundaries of $\lda{S}{B}$ through $P$.
Hence $C'_2\not \subset \supp B$.
So $K_S.C'_2=-B.C\leq -(B_1+B_2).C'_1<0$, which shows (1).

If $P$ is nonsingular on $S$, we obtain $K_S.C'_2\leq -(B_1+B_2).C'_2\leq-2$,
a contradiction to the contractibility of $C'_2$. 
Hence $P$ is singular.
 
Define $B_T$ by $K_T+B_T=\pa{\mu}{S}{B}$.
Then $B_T.\-{\mu}C'_2\geq 1$ since $\-{\mu}C'_2$ intersects log
crepant curves over $P$.
But by (1), this must be equality.
From this we can deduce the latter half of (2) and (3).
\qed
\enddemo

Let $\nu:U\to S$ be MLT for $\lda{S}{B}$.

Below we assume that $C_1\cap B_2 \in \clc{0}(K_S+B)$
to get a contradiction.

\proclaim{Claim 2}
$C_1\cap C_2=\phi$
(See FIGURE (III).)
\endproclaim

\demo{Proof}
Assume there is a component $C'_2$ of $C_2$ such that $C_1\cap
C'_2\not = \phi$. Let $C'_1$ be a component of $C_1$ such that
$C'_1\cap C'_2\not = \phi$.
Then by Claim 1, $C'_1\not \subset \supp \{ B \}$.
Hence by (ii) and (iii), $K_S.C'_1\leq -B_2.C'_1<0$,
which shows that $\-{\mu}C'_1$ is a $(-1)$-curve.
On the other hand, $\-{\nu}C'_1\cap \-{\nu}C'_2\not = \phi$ 
since by Claim 1, $C'_1\cap C'_2 (\not = \phi)$ is at worst Du Val
singularity and not contained in $\supp B$.

Furthermore we know that $C'_1\cap B_2$ is singular as follow:

If $C'_1\cap B_2$ is smooth, we have
$K_S.C'_1\leq -B_2.C'_1\leq -1$ by $C'_1\not \subset \supp B$.
So this must be equality since $\-{\mu}C'_1$ is a $(-1)$-curve.
Hence $C'_1\cap \supp \{ B \} = \phi$.
But this shows that $C'_1\cap B_2 \not \in \clc{0}(K_S+B)$, a contradiction.
 
Next we claim that $(\-{\nu}B_1)^2<0$ as follows:

By the assumption (iii), $\-{\nu}B_1 \cap  \-{\nu}C'_1=\phi$.
By Claim 1 (2), $\-{\nu}B_1\cap \-{\nu}C'_2=\phi$.
Clearly $\nu$-exceptional curves  over $C'_1\cap B_2$ and $\-{\nu}B_1$
doesn't intersect.
Now contract $\nu$-exceptional curves over $C'_1\cap C'_2$ and
$\-{\nu}C'_1$ then $(\-{\nu}C'_2)^2$ becomes $\geq 0$ (cf. Claim 1 (2)
and the fact $\-{\mu}C$ is a $(-1)$-curve).
Next contract $\nu$-exceptional curves over $C'_1\cap B_2$ 
(Note that there exist such exceptional curve by Claim 1 (2)
and the fact $P\in \clc{0}(\lda{S}{B})$).
then $(\-{\nu}C'_2)^2$ becomes $>0$.
Even after such contractions, $\-{\nu}B_1$ and $\-{\nu}C'_2$ doesn't
intersect since such contractions are isomorphic near $\-{\nu}B_1$. 
Hence by Hodge index theorem, $(\-{\nu}B_1)^2<0$.

By the way we know that $B_1\simeq \Bbb P^1$ 
and $K_S+B|_{B_1}=K_{B_1}+\frac12 P_1+\frac12 P_2+P$
by the assumption (iv).

Here we claim that one of $P_1$ and $P_2$, say $P_1$ is smooth
and $P_2$ is an ODP as follow:
If both $P_1$ and $P_2$ are singular, by the above subadjunction
they are ODP's.
Hence we can write 
$$\mu^*{B_1}=\-{\mu}B_1+\sum e_j E_j+\frac12 E^1+\frac12 E^2,$$
where $E^i$ is the exceptional curve resolving $P_i$
and $E_j$'s are exceptional curves over $P$ and $e_j\leq 1$.
(We assume that $E_1$ intersects $\-{\mu}B_1$.)
By $(B_1)^2>0$ and $(\-{\mu}B_1)^2<0$,
$0<(\-{\mu}B_1)^2+e_1+1$ and hence $\-{\mu}B_1$ is a $(-1)$-curve.
On the other hand
since $\-{\nu}B_1$ is also (possibly analytically) contractible,
$E^1$, $\-{\mu}B_1$ and $E^2$ is a contractible chain.
But this is absurd.
By similar calculation, $P_2$ must be an ODP.

By the above subadjunction, we know that there is a component $B_3$
of $\{ B \}$ which passes through $P_1$.

Next we claim that $(\-{\nu}B_3)^2<0$ below.
First by Claim 1, $\-{\nu}B_3 \cap \-{\nu}C'_2=\phi$.
Furthermore we obtain also $\-{\nu}B_3\cap \-{\nu}C'_1=\phi$ as
follows:

Looking the degree of $K_S+B|_{B_2}$,
Only $P$ and $B_2\cap C'_1$ can contribute to its boundary. 
On the other hand by the ampleness of $B_2$, $B_2\cap B_3 \not =\phi$
and hence $B_2\cap B_3 = B_2 \cap C'_1$.
Consider the description of MLT over $B_2\cap B_3$.
If $\-{\nu}B_3\cap \-{\nu}C'_1 \not =\phi$, 
we have $K_T.\-{\mu}C'_1\leq -\frac32$, a contradiction to
that $\-{\mu}C'_1$ is a $(-1)$-curve.
Hence $\-{\nu}B_3\cap \-{\nu}C'_1=\phi$.
Now contract first $\-{\nu}C'_1$.
Then $(\-{\nu}C'_2)^2$ become $\geq 0$.
Next contract $\nu$-exceptional curve over $P$. Then 
$(\-{\nu}C'_2)^2$ become $>0$.
These contractions are isomorphic near $\-{\nu}B_3$ 
since $\-{\nu}B_3 \cap \-{\nu}C'_1=\phi$
and
hence even after the contractions, $\-{\nu}B_3\cap \-{\nu}C'_2$ remain
empty.
Hence by Hodge index theorem, $(\-{\nu}B_3)^2<0$.

But this gives the final contradiction.
For $\-{\nu}B_3$ intersects at least one log crepant curve over
$B_2\cap B_3$ and also $\-{\nu}B_1$ at another point.
Hence applying (1.6) for the contraction of $B_3$,
we get a contradiction.
We obtain Claim 2.
\qed
\enddemo

By this Claim 2, even after contracting $C_2$,
every assumption of (1.9) and assumption that 
$C_1\cap B_2\in \clc{0}(K_S+B)$ hold
by setting $C_2=\phi$.
Hence it is sufficient to get a contradiction assuming also $C_2=\phi$.

Consider extremal rays for $K_S+(1-\epsilon)(B-B_1)$,
where $\epsilon$ is a sufficiently small positive rational number.
Note that $K_S+(1-\epsilon)(B-B_1)$ is KLT by (1.7).
\proclaim{Claim 3}
There is an extremal ray $R$ for $K_S+(1-\epsilon)(B-B_1)$
such that $B_1.R>0$.
\endproclaim

\demo{Proof}
Assume that any extremal ray $R'$ for $K_S+(1-\epsilon)(B-B_1)$
satisfies $B_1.R'\leq 0$.
For an irreducible curve $l$ such that 
$(K_S+(1-\epsilon)(B-B_1)).l\geq 0$,
we have $((1-\epsilon)B_1+\epsilon B).l\leq 0$.
Since $\epsilon$ is sufficiently small, we have
$B_1.l\leq 0$. But by the cone theorem, for every curve $l$,
we must have $B_1.l\leq 0$, a contradiction to $(B_1)^2>0$.
\qed
\enddemo

Fix a extremal ray $R$ as in Claim 3.

If the contraction of $R$ is birational,
by the choice of $R$ and ampleness of $B_2$, 
$B_1\cap R \not = \phi$ and $B_2\cap R \not = \phi$.
By (1.6), $B_1\cap R=B_2\cap R=\{ P \}$.
Now set $C_2=\supp R$.
Then again all the assumptions of this theorem and the assumption
$C_1\cap B_2\in \clc{0}(K_S+B)$ hold.
Hence by considering inductively,
we have the maximal composition of birational extremal contractions 
$S\to \o{S}$
such that all extremal ray intersect $B_1$ positively
and for $\o{S}$, all the assumptions of this theorem and the assumption
$\o{C_1}\cap \o{B_2}\in \clc{0}(K_{\o{S}}+\o{B})$ hold.
(We denote all the image on $\o{S}$ with $\o{ \ }$.)
For this $\o{S}$, we can take an extremal ray $\o{R}$ as in Claim 3. 
By the choice of $\o{S}$, the contraction of $\o{R}$ is not birational.
By (iii), it is not a contraction to a point.
Hence it is a contraction to a curve.
Let $l$ be its fiber such that $\o{B_2}\cap \o{C_1} \in l$. 
By $\o{B_1}.l>0$, $\o{B_1}\cap l \not = \phi$.
Let $\o{\nu}:\o{T}\to \o{S}$ be MLT for $K_{\o{S}}+\o{B}$.
Then
$(\-{\o{\nu}}l)^2<0$ 
since $\o{B_2}\cap \o{C_1}$ is contained in $\clc{0}(K_{\o{S}}+\o{B})$ 
and singular.
Note that $\o{B_1}\cap l \not = \o{B_2}\cap \o{C_1}$ and 
$l \not \subset \llc{\o{S}}{\o{B}}$.
Furthermore $\-{\o{\nu}}l$ intersects 
$\text{LLC} \ {\o{\nu}}^*(K_{\o{S}}+\o{B})$
at two point (and is not contained in it), 
a contradiction to (1.6).
We finish the proof of this theorem.
\qed
\enddemo 

The following proposition is used in the treatment of case (b) in Set up (3.0).

\proclaim{Proposition (1.10)}
Let $(S, K_S + B, s)$ be a pointed LC surface.
Assume that 
$C := \llcorner B \lrcorner$ is a smooth projective curve.
Let $a_i:S_i \to S$ for $i = 1, 2$ be a birational morphism
satisfying the following:

\roster
\item $a_i$ ($i = 1, 2$) is log crepant for $K_S + B$ and
the exceptional curve $F_i$ for $a_i$ is irreducible;
\item $\clc{0}(a_1^*(K_S + B)) \cap F_1 = \-{a_1}C \cap F_1$;
\item $a_2^*(K_S + B)$ is LT at $\-{a_2}C \cap F_2$ (and hence $S_2$ is 
smooth at $\-{a_2}C \cap F_2$). 

\endroster

Then  $(\-{a_2}C)^2 > 0$ (resp. $(\-{a_2}C)^2 \geq 0$) implies 
that $(\-{a_1}C)^2 > 0$ (resp. $(\-{a_1}C)^2 \geq 0$). 
\endproclaim

\demo{Proof}
Let $m:\t{S_1} \to S_1$ be the MLT for $a_1^*(K_S + B)$
at $\-{m}C \cap F_1$. 
By this and assumption (2), the assumption that $C$ is irreducible and 
the uniqueness of MLT (1.0),
$a_1 \circ m$ must be the MLT for $K_S + B$.
Let $n:\t{S_2} \to S_2$ be the MLT for $a_2^*(K_S + B)$ near $F_2$.
Possibly after contracting some $a_2 \circ n$-exceptional curves,
we obtain the MLT for $K_S + B$.
Assume that $(\-{a_2}C)^2 > 0$ (resp. $(\-{a_2}C)^2 \geq 0$). 
Then by assumption (3), there is no
$n$-exceptional curve over $\-{a_2}C \cap F_2$ and hence 
$(\-{(a_2 \circ m)}C)^2 > 0$ (resp. $(\-{(a_2 \circ m) }C)^2 \geq 0$). 
Since the MLT for $K_S + B$ is dominated by  
$\t{S_2}$, we have also $(\-{(a_1 \circ n)}C)^2 > 0$ 
(resp. $(\-{(a_1 \circ n)}C)^2 \geq 0$) and hence
$(\-{a_1}C)^2 > 0$ (resp. $(\-{a_1}C)^2 \geq 0$).
\qed
\enddemo
 
\proclaim{Proposition (1.11) (Smoothness)}
Let $X$ be a normal variety and $\ldb{X}{S}{B}$ is LT,
where $S$ is a $\Bbb Q$-Cartier integral Weil divisor which contains no
codimension $2$ singular locus.
Assume that $S$ is nonsingular.
Then $X$ is nonsingular.
\endproclaim
\demo{Proof}
See [S, 3.7, Corollary].
\qed
\enddemo 

\proclaim{Proposition (1.12) (Classification of LT singularity)}
Let $(X, \lda{X}{B}, x)$ be a germ of a normal LT $3$-fold such that
$\lda{X}{\llcorner B \lrcorner}$ is $\Bbb Q$-Cartier (and hence LT).
Then
\roster
\item if $B$ has at least $3$ reduced components which are $\Bbb Q$-Cartier,
then locally analytically 
$$(x \in B \subset X) \simeq (o \in (xyz =0) \subset \Bbb C^3);$$ 
\item if $B$ has $2$ components which are $\Bbb Q$-Cartier, then
locally analytically
$$(x \in \llcorner B \lrcorner \subset X) 
\simeq (o \in (xy=0) \subset \Bbb C^3/ \Bbb Z_n(q_1, q_2, 1),$$
where $q_1$ and $q_2$ are natural numbers such that $(q_1, q_2, 1)=1$.
\endroster
\endproclaim

\demo{Proof}
In [FA,16.10 Corollary] and [FA,16.15 Theorem], this theorem is proved
under the assumption that $\lda{X}{B}$ is DLT. But as we will see below,
we can prove 
under the above assumptions that 
$(X, x)$ is analytically $\Bbb Q$-factorial and 
hence $\lda{X}{\llcorner B \lrcorner}$ is DLT.

If $(X, x)$ is not analytically $\Bbb Q$-factorial, then
we can take 
a small analytic \newline
$\Bbb Q$-factorialization $\pi \t{X} \to X$ 
by applying (1.13) for $\lda{X}{B_1}$.
In fact everything is clear except that $\pi$ is small.
Since $\pi$ is log crepant for $\lda{X}{B_1}$ and $\lda{X}{B_1}$ is PLT,
it must be small.
Let $B_1$ and $B_2$ be two of reduced components in $\llcorner B \lrcorner$
which are $\Bbb Q$-Cartier. Since $B_i$'s are $\Bbb Q$-Cartier,
exceptional curves of $\pi$ are contained in $\-{\pi}B_i$'s.
Hence $\-{\pi}B_1 \cap \-{\pi}B_2$  is not irreducible.
Now we apply [FA,16.15 Theorem] for $\ldb{\t{X}}{\-{\pi}B_1}{\-{\pi}B_2}$,
a contradiction. 

Hence we can directly apply [FA,16.10 Corollary] and [FA,16.15 Theorem]
for $\lda{X}{\llcorner B \lrcorner}$. 
Note that in (1), $B=\llcorner B \lrcorner$.
In fact the simple blow up
of $x$ is log crepant for $\lda{X}{\llcorner B \lrcorner}$ and hence
there can be no component other than $\llcorner B \lrcorner$
since $\lda{X}{B}$ is also LT.
\qed
\enddemo

\proclaim{Theorem (1.13) ($\Bbb Q$-factorial LT model)}
Let $X$ be a $3$-fold and $B$ a boundary on $X$.
Assume that $B$ is $\Bbb Q$-Cartier and $\llcorner B \lrcorner$ is an
LSEPD divisor, i.e.,
there is a principal divisor $(h)$ near $\supp  \llcorner B \lrcorner$
such that $\supp (h) = \supp \llcorner B \lrcorner$.
Then in a neighborhood of $\llcorner B \lrcorner$, 
there exists a $\Bbb Q$-factorial LT model for $K_X + B$
, i.e., there is a birational morphism $f: Y \to X$
such that $K_Y + B_Y$ is $f$-nef and over a neighborhood
of $\llcorner B \lrcorner$, $K_Y +B_Y$ is LT,
where $B_Y$ has all $f$-exceptional divisors with coefficients $1$. 
\endproclaim

\demo{Proof}
See [S, Corollary 5.19] or [FA, Corollary 20.9].
\qed
\enddemo

\head 2. Auxiliary flips and Back tracking method
\endhead

\proclaim{Proposition (2.0) (Flip of type 1)}
Let $X$ be a $\Bbb Q$-factorial $3$-fold. Let $D=S+B$ be a $\Bbb Q$-divisor,
where $S$ is reduced and $\llcorner B \lrcorner = 0$.
Let $f:X \to Y$ be a small contraction with $\rho(X/Y)=1$
and $C$ its exceptional curve. 
Assume that
\roster
\item $\lda{X}{D}$ is LC and $(\lda{X}{D}).C \leq 0$;
\item $S$ has at least $2$ irreducible components $S_1$ and $S_2$ such that
$S_1.C<0$ and $S_2.C>0$.
\endroster
Then the flip of $f$ exists.
\endproclaim

\demo{Proof}
See [FA, 20.7 Theorem].
\qed
\enddemo

\proclaim{Proposition (2.1) (Flip of type 2)}
Let $f:X \to Y$ be a small contraction 
and $C$ its exceptional curve. Let $S_1$ and $S_2$ are effective divisors
such that $S_1 \cap S_2=C$. Assume that for some $m_1, m_2 \in \Bbb N$,
$m_1 S_1 \sim m_2 S_2$ near $C$. Then the flip of $f$ exists.  
\endproclaim
\demo{Proof}
See [FA, 20.11 Theorem].
\qed
\enddemo

\proclaim{Proposition (2.2) (Flip of type 3)}
Let $X$ be a $\Bbb Q$-factorial $3$-fold. Let $D=S+B$ be a $\Bbb Q$-divisor,
where $S$ is reduced.
Let $f:X \to Y$ be a small contraction with $\rho(X/Y)=1$
and $C$ its exceptional curve.
Assume that
\roster
\item $f$ is special with respect to $\lda{X}{S}$;
\item $2(\lda{X}{D}) \sim 0$ near $C$;
\item $\lda{X}{D}$ is LT;
\endroster
Then the flip of $f$ exists.
\endproclaim

\demo{Proof}
See [FA, Proposition 21.5].
\qed
\enddemo

\proclaim{Proposition (2.3) (Flip of type 4)}
Let $X$ be a $\Bbb Q$-factorial $3$-fold. Let $D=S+B$ be a $\Bbb Q$-divisor,
where $S$ is irreducible and reduced.
Let $f:X \to Y$ be a small contraction with $\rho(X/Y)=1$
and $C$ its exceptional curve.
Assume that
\roster
\item $f$ is special with respect to $\lda{X}{S}$;
\item $2(\lda{X}{D}) \sim 0$ near $C$;
\item $\lda{X}{D}$ is LC;
\item $\lda{X}{D}|_S$ is exceptional, i.e., 
there is exactly one divisor in $\text{PLC} (\lda{X}{D}|_S)$;
\endroster
Then the flip of $f$ exists.
\endproclaim

\demo{Proof}
See [FA, 21.1.2.1 Corollary].
\qed
\enddemo

\proclaim{Theorem (2.4)}
Let $X$ be an algebraic $3$-fold, 
$S_1$ and $S_2$ irreducible $\Bbb Q$-Cartier
surfaces.
Let $f: X\to Y$ be a small contraction.
We assume the following:
\roster
\item"(a)"
 $\ldb{X}{S_1}{S_2}$ is LC;
\item"(b)" $C:=\text{excep} \ f$ is an irreducible curve and 
$C \subset S_1\cap S_2$;
\item"(c)" $(\lda{X}{S_1}).C<0$, $S_1.C<0$ and $(\ldb{X}{S_1}{S_2}).C=0$;
\item"(d)" $2(\ldb{X}{S_1}{S_2})\sim 0$ near $C$;
\item"(e)" Along $C$, $S_1$ and $S_2$ are generically normal crossing 
and $S_1$ is normal;
\item"(f)" $K_X+S_1+S_2|_C=K_C+\frac12 P_1+\frac12 P_2+P$, where $P_1$,
$P_2$ and $P$ are distinct points.
\endroster

Then the flip of $f$ exists. 
Let $f^+:X^+\to Y$ be the flip and $C^+$ the exceptional curve for
$f^+$.
We denote all the strict transforms of divisors on $X$ with
superscript $+$.
Then we obtain the following descriptions;
\roster
\item One of $P_1$ and $P_2$, say $P_1$ is an ODP on $S_1$.
$P_2$ is a smooth point on $S_1$.
$P$ is smooth or a singularity with the following description on the 
MRS of $S_1$: FIGURE (IV)
of type $(3,2,2,2,2,...,2)$;
\item Near $P$, $S_2|_{S_1}=C+C'$ as a Weil divisor on $S_1$, 
where $C'$ is a smooth curve.
In particular along $C'$, 
$S_1$ and $S_2$ are also generically normal crossing and
$\ldb{X}{S_1}{S_2}$ is LT outside $P$;
\item The natural birational map $S_1\dashrightarrow (S_1)^+$
is a morphism whose exceptional locus is $C$.
The image point of $C$ (we call this $Q$) is smooth on $(S_1)^+$;
\item $C^+$ is irreducible and near $Q$, $C^+$ coincides with the
nonnormal locus of $S_2^+$.
Let $\nu: S_2^{+\nu}\to S_2^+$ be the normalization.
Then $\-{\nu}Q$ is one point and a nonsingular point of 
$S_2^{+\nu}$.
\endroster
(See FIGURE (VI).)
\endproclaim

\demo{Proof}
Since $f$ is a flipping contraction of type 1, 
the flip exists by (2.0).
We prove the descriptions $(1)\sim (4)$ below.

First we prove $(1)$ and $(2)$.
\proclaim{Claim 1}
$S_1\cap S_2$ is not irreducible.
\endproclaim
\demo{Proof}
If $S_1\cap S_2$ is irreducible, $C^2_{S_1}=S_2.C>0$
by the assumption (c) and (e).
But this contradicts the contractibility of $C$ in $S_1$.
\qed
\enddemo
By LC property of $\ldb{X}{S_1}{S_2}|_{S_1}$ and Claim 1,
there are two components of $S_1\cap S_2$ at $P$.
Let $C'$ be the component other than $C$.
By [FA, 16.6.3] and the assumption (e), 
we can write $S_2|_{S_1}=C+\frac1mC'$, where $m$ is a natural number.
We prove the description (1).
Let $\mu: (S_1)^{\mu}\to S_1$ be the MRS.
Assume that both $P_1$ and $P_2$ are singular, i.e., ODP's by (f). 
Let 

$$
\mu^*(C+\frac1m C')=
\-{\mu}C+\frac1m\-{\mu}C'+\frac 12 E^1+\frac 12 E^2 + 
\sum \limits_{i=1}^{n} e_iE_i, \tag $\bigstar$ $$
where $E^i$ is the $(-2)$-curve
resolving $P_i$ and $E_j$'s are exceptional curves resolving $P$
and $E_1$ is the nearest component to $\-{\mu}C$.
We note that $e_i\leq 1$ by the minimality of $\mu$ and LC property
of $\ldb{S_1}{C}{\frac1m C'}$.
Intersecting $\bigstar$ with $\-{\mu}C$, we obtain

$0<S_2.C=(\-{\mu}C)^2+1+$ 
$$
\cases
\frac 1m & {\text{if} \  n=0}\\
e_1 & {\text{if} \ n\geq 1.}
\endcases $$

Hence by this and the contractibility of $\-{\mu}C$,
$(\-{\mu}C)^2 = -1$.
But we obtain a non-contractible chain consisting of $E^1$, 
$\-{\mu}C$ and $E^2$, a contradiction.

Assume that both $P_1$ and $P_2$ are smooth.
Then by similar calculations, we obtain $(\-{\mu}C)^2\geq 0$,
a contradiction.

Hence we proved the description about $P_1$ and $P_2$.

Furthermore by the similar calculations,
we obtain $(\-{\mu}C)^2 = -1$ and 
$$\cases
\frac 1m >\frac 12 & {\text{if} \ n=1} \\ 
e_1>\frac 12 & {\text{if} \ n\geq 1.} \\
\endcases$$

If $n=1$, we immediately obtain $(1)$ and $(2)$. 

If $n\geq 1$, we see that $E_1^2 \leq -3$
since otherwise $E^1$, $\-{\mu}C$ and $E_1$ construct a
non-contractible chain.
Based on this and $e_1 > \frac12$, we obtain
$P$ is a singularity as in the statement of (1)
and 
$$\mu^*C=
\-{\mu}C+\frac 12 E^1+\sum \limits_{i=1}^{n} \frac {n-i+1}{2n+1} E_i$$ 
and 
$$\mu^*C'=
\-{\mu}C'+\sum \limits_{i=1}^{n} \frac {2i-1}{2n+1} E_i$$
by elementary calculations.
Hence $e_1=\frac {n}{2n+1}+\frac 1m \frac1 {2n+1} >\frac 12$,
i.e., $m=1$.
So in this case we proved $(1)$ and $(2)$.

Next we prove $(3)$.
By $(1)$ and $(2)$, $f(C)$ is a smooth point of $f(S_1)$.
In particular $f(S_1)$ is normal. 
Hence if we prove that no component of $C^+$ is contained in $S_1^+$,
then by the Zariski's Main Theorem, $S_1^+\simeq f(S_1)$ and $(3)$ follows.
\comment
Since  is , 
$\ldb{X}{S_1}{S_2}|_{S_1}=\pb{f'}{Y}{f(S_1)}{{f(S_2)|_{f(S_1)}}},
where $f':= f|_{S_1}$.
On the other hand 
\endcomment
Note that
$\ldb{X}{S_1}{S_2}|_{S_1}=\ldc{S_1}{C}{\frac12 C''}{C'}$ by (f) and (1),
where $C''$ is an irreducible curve 
which crosses normally with $C$ at $P_2$.
Here
$\ldb{Y}{f(S_1)}{f(S_2)|_{f(S_1)}}=\ldb{f(S_1)}{\frac 12 f(C'')}{f(C')}$
is LC since $\ldb{X}{S_1}{S_2}$ is $f$-trivial.
Hence $\lda{Y}{f(S_1)}|_{f(S_1)}=
\lda{f(S_1)}{\frac 12 f(C'')}$ is canonical since $f(C)$ is smooth
on $f(S_1)$.
By this and the $f^{'+}$-ampleness of $\lda{X^+}{S_1^+}|_{S_1^+}$
(we define that $f^{'+} := f^+|_{S_1^+}$),
$S_1^+$ contains no exceptional curve of $f$.
Hence we finish the proof of $(3)$.

Next we prove $(4)$.
Since $S_1^+$ is $\Bbb Q$-Cartier, any analytic branch of $S_2^+$  
intersects $S_1^+$ along curves. 
On the other hand by $(2)$ and $(3)$,
$S_1^+ \cap S_2^+ =C^{\prime +}$ and the intersection is generically NC
so the number of analytic branch
of $S_2^+$ is $1$, i.e., $\-{\nu}Q$ is one point.

\proclaim{Claim 2}
The index of $\lda{X^+}{S_2^+}$ at $Q$ is $1$.
\endproclaim
\demo{Proof}
Assume the contrary.
Let $\pi : \t{X^+} \to X^+$  be the index $1$ cover of  
$\lda{X^+}{S_2^+}$. 
On $S_1^+$, $X^+$ has only generically ODP curve singularity 
(along ${C''}^{+}$) near $Q$.
Hence $\-{\pi} S_1^+ \to S_1^+$ is only ramified at $Q$ near $Q$.
But $S_1^+$ is smooth at $Q$ so $\-{\pi}S_1^+$ has at least $2$
components intersecting only at $\-{\pi} Q$. 
Since $\lda{\t{X^+}}{\-\pi{S_1^+}}$ is PLT by (c) and $(2)$,
we can take a small $\Bbb Q$-factorialization $\rho: \o{X}^+\to \t{X}^+$
by (1.13). 
The exceptional curves of $\rho$
are contained in $\-{\rho}\-{\pi}{S_1}$ and $\-{\rho}\-\pi{S_2}$
(they have at most $3$ components in all.)
But this contradicts LC property of 
$\ldb{\o{X^+}}{\-{\rho}\-{\pi}S_1^+}{\-{\rho}\-{\pi}S_2^+}$
along  the exceptional curves of $\rho$.
Hence  we obtain Claim 2.
\qed
\enddemo

\proclaim{Claim 3}
Near $\-{\nu}Q$, $\ldb{X^+}{S_1^+}{S_2^+}|_{S_2^{+\nu}}=
\ldb{S_2^{+\nu}}{C^*}{C^{\prime +}}$, where $C^*$ is some curve
through $\-{\nu}Q$. (We will prove later that $C^*:=\-{\nu}C^+$.)

\endproclaim
\demo{Proof}
Since $\ldb{X^+}{S_1^+}{S_2^+}|_{S_2^{+\nu}}$ is LC 
at $\-{\nu}Q$ and $\lda{X^+}{S_2^+}|_{S_2^{+\nu}}$ is LT
of index $1$ at $\-{\nu}Q$ by Claim 2, 
the boundary near $\-{\nu}Q$ is reduced 
and the number of it is at most $1$.
Hence we have only to prove
that $\lda{X^+}{S_2^+}|_{S_2^{+\nu}}=
\lda{S_2^{+\nu}}{C^*}$ for some irreducible curve $C^*$.
Assume the contrary. Then $\lda{X^+}{S_2^+}|_{S_2^{+\nu}}=
K_{S_2^{+\nu}}$.
By this, we know that $S_2^+$ is normal at $Q$.
Since $\ldb{X^+}{S_1^+}{S_2^+}|_{S_2^+}=
\lda{S_2^+}{C^{\prime +}}$ and $K_{S_2^+}$ is Cartier at $Q$,
$K_{S_2^+}$ is canonical at $Q$.

On the other hand 
$\ldb{X^+}{S_1^+}{S_2^+}|_{S_2^+}$ is not LT.
In fact if it is LT, it is PLT.
Hence its restriction to $C^{\prime +}$ is KLT
and by the inversion of adjunction [FA, 17.6 Theorem], 
$\ldb{X^+}{S_1^+}{S_2^+}|_{S_1^+}$ is PLT.
But this contradicts the description $(1)$ and $(2)$
(see FIGURE (V)).

So by (1.12), $Q$ is a canonical singularity
of type $D_m$ for some $m$.
Let $C_1^+$ be any component of $C^+$ and 
$\sigma: S_2^{+\sigma} \to S_2$ the MRS.
Since $K_{S_2^+}.C_1^+=(\lda{X^+}{S_2^+}).C_1^+<0$,
$C_1^+$ is $(-1)$-curve on the MRS of $S_2$.
But the union of the exceptional curves of MRS of 
a canonical singularity of type $D_m$
and $C_1^+$ cannot be contracted 
at one time, a contradiction.
\qed
\enddemo
By Claim 2 and Claim 3, we see that $S_2^{+\nu}$ is smooth at $\-{\nu}Q$.

In fact it is derived 
by the fact that$\lda{S_2^{+\nu}}{C^*}$ is Cartier and LT. 

By Claim 3, $\nu(C^*)$ coincides with the nonnormal locus near $Q$.
Since $S_1^+$ is normal outside $C^+$ by (e) and
all the components of $C^+$ pass through $Q$, 
every $1$ dimensional component of the nonnormal locus of $S_1^+$ is
contained in $C^+$ (and passes through $Q$).
Hence $\nu(C^*) \subset C^+$.
If there is a component $C_1^+$ of $C^+$ which is not contained in $\nu(C^*)$,
we have $K_{S_2^{+\nu}}.\-{\nu}C_1^+ =- ({C^*}+{C^{\prime +}}).\-{\nu}C_1^+
\leq -2$, a contradiction to the contractibility of $\-{\nu}C_1^+$ in 
$S_2^{+\nu}$.
Hence $\nu(C^*)=C^+$ and in particular $C^+$ is irreducible.

Now we finish the proof of (2.4).
\qed
\enddemo

\definition{Description of BTM (2.5)}
We review the back tracking method for introducing the notation. 
The material here is contained in [FA, Chapter 21].
Let $X$ be an algebraic $\Bbb Q$-factorial normal $3$-fold and  
$f:X\to Y$ a projective small contraction such that $\rho(X/Y)=1$.
Let $S$ be an irreducible divisor on $X$, $B$ a boundary 
and $C$ the exceptional curve of $f$.
We assume the following assumptions.
\roster
\item "(i)" $\lda{X}{S}$ is PLT and $\ldb{X}{S}{B}$ is LC.
\item "(ii)" $(\lda{X}{S}).C<0$, $S.C<0$ and $(\ldb{X}{S}{B}).C=0$.
\endroster

Then we can make the following procedure:

Choose a primitive log crepant divisorial extraction 
$g_1:Z_1\to X$ for $\ldb{X}{S}{B}$.
Let $E_1$ be the exceptional divisor.

By assumption, $\ldc{Z_1}{S_1}{B_1}{E_1}=\pb{g_1}{X}{S}{B}$,
where $S_1$ and $B_1$ are the strict transforms of
$S$ and $B$.
Hence by choosing a rational number $b \ (0<b<1)$ appropriately,
we have $\ldc{Z_1}{S_1}{bB_1}{e'E_1}=\pb{g_1}{X}{S}{bB}$,
where $e'\in (0,1)$ is a rational number.
Here note that $(\ldb{X}{S}{bB}).C<0$ by the assumption $(2)$.
So by choosing a rational number $e\in (e',1)$,
$-(\ldc{Z_1}{S_1}{bB_1}{eE_1})$ become $f\circ g$-ample.
Furthermore $\ldc{Z_1}{S_1}{bB_1}{eE_1}$ is PLT by the following argument:

By (i), $\ldb{X}{S}{bB}$ is PLT. Hence $\ldc{Z_1}{S_1}{bB_1}{e'E_1}$
is also PLT. The PLT property is an open one so 
$\ldc{Z_1}{S_1}{bB_1}{eE_1}$ is also PLT.
 
Consequently we can run the $\ldc{Z_1}{S_1}{bB_1}{eE_1}$-MMP over $Y$.
Since $\rho(Z_1/Y)=2$, the surgeries are uniquely determined
once we choose $g_1$.
We obtain the following diagram:
$$\matrix 
& & Z_1 & &\dashrightarrow & & Z_2 & & \dashrightarrow \\
& {\swarrow g_1} & & {h_1 \searrow} & & {\swarrow g_2} & & {h_2 \searrow} \\
X & & & & X_1 & & & &X_2.... \\
& {\searrow f} & &  \swarrow & & \searrow  & & \swarrow \\
& & Y & & & & Y &.\\
\endmatrix
$$
We denote the extremal ray associated to $h_i$ by $R_i$.
If $h_i:Z_i\to X_i$ is a flipping contraction,
we denote the exceptional curve of $h_i:Z_i\to X_i$ by $C_i$  
and the exceptional curve of $g_{i+1}:Z_{i+1}\to X_i$ by $C_i^+$.
We also denote the strict transforms of $S_1$, $B_1$ and $E_1$ on $Z_i$
by $S_i$, $B_i$ and $E_i$.
Then we have the following important numerical properties:
\enddefinition
\definition{Numerical properties}
\roster
\item For some positive rational numbers $p$ and $q$,
$(pS_i+qE_i).C_i<0$;
\item For some positive rational numbers $r$ and $s$,
$(S_i+rB_i+sE_i).C_i=0$;
\item For fixed $i$, we have $S_{i+1}.C_i^+>0$ or $S_{i+1}.C_{i+1}>0$.
\endroster
For the proof, see [FA, 21.6.2 Lemma, 21.7.1 Proposition and 21.7.2 Lemma].

If we encounter divisorial contraction $h_n:Z_n\to X_n$,
then $X_n\to Y$ is the flip of $f$ ([FA, 6.5.5]).
Hence if all flips in BTM exist and terminate, then the flip of $f$
exists. So we need to examine such flips in detail.

Roughly speaking there are three types of such flips as follow:
\enddefinition

\definition{ (a) beginning flip}

Let $$j:= 
\text \ {min} \ \{ k | h_k \ \text {is flipping and}\  
 E_k.C_k \leq 0, \ \text {or} \ h_k \ 
\text{is divisorial.} \ \}.$$
We call $g_i$ with $1 \leq i \leq j-1$ a beginning flipping contraction. 

By the numerical property $(1)$,
we have $E_i.C_i >0$ and $S_i.C_i<0$. 
Hence this is a flipping contraction of type 1 and the
flip exists by (2.0) and terminates (for termination, see [FA, \S 7]).
\enddefinition

After beginning flips, we encounter $h_j$.
If it is a divisorial contraction, then we are done. 
Assume that $h_j$ is a flipping contraction.

\definition {(b) middle flip}

$C_j$ satisfies $E_j.C_j \leq 0$ by definition.
If $C_j$ satisfies also $S_j.C_j\leq 0$, we call $h_j$ a middle flipping
contraction.
\enddefinition

We treat this flip by dividing into the following cases:

\definition {(b1) $E_j.C_j < 0$ and $S_j.C_j <0$}

If $C_j \subsetneqq S_j \cap E_j$,
there is no general strategy to prove the existence of the flip of $g_j$.
In our application, we can choose $g$ so as to avoid this situation.
 
If $C_j=S_j \cap E_j$, it is a flip of type 2 and the flip
of $h_j$  exists by (2.1) and $h_{j+1}$ is a divisorial contraction
(see [FA, 21.8.2.1]). Hence in this case we are done.
\enddefinition

\definition {(b2) $E_j.C_j=0$ and $S_j.C_j<0$}

There is no general strategy to prove the existence of the flip of $g_j$.
In our application, we can choose $g$ so as to avoid this situation.
\enddefinition

\definition {(b3) $E_j.C_j<0$ and $S_j.C_j=0$}

If some component of $C_j$ is contained in $S_j$,
there is no general strategy to prove the existence of the flip of $g_j$.
In our application, we can choose $g$ so as to avoid this situation.

The left case is that any component of $C_j$ does not intersect
$S_j$. The heart of the Shokurov's paper is to prove the existence
of the flip for such a flipping contraction by induction
(see $\S 4$ and $\S 7$).
\enddefinition
 
Assume that $h_j$ satisfies $(b3)$ and   
any component of $C_j$ does not intersect $S_j$.
Assume furthermore that the flip of $h_j$ exists and $h_{j+1}$ is also
a flipping contraction. By the numerical property (3),
we have $S_{j+1}.C_{j+1} >0$.

\definition {(c) final flip}
We can prove that if $h_k$ is a flipping contraction such that
$S_k.C_k >0$ for some $k$,
every $C_l$'s for $l \geq k$ satisfies
$S_l.C_l>0$ (see [FA, 21.8.3.1 Lemma]).
We call such a $h_l$ a final flipping contraction.
By the numerical property (1),
we have $E_k.C_k<0$.
Hence $h_k$ is of type 1.
So the flip exists by (2.0) and terminate (for termination, see [FA, \S 7]). 
After final flips, we obtain a divisorial contraction. 
By contracting the divisor we obtain the flip of $f$.
\enddefinition

We remark that if $S_j \cap E_j$ is irreducible,
then $h_j$ satisfies the assumption of $(b1)$ and that $S_j \cap E_j = C_j$ 
or the assumption of $(b2)$ and that 
any component of $C_j$ does not intersect $S_j$.
If the latter case occurs, once we clear the flip of $C_j$ final flipping
contractions or a divisorial contraction appear and we are done.

\head 3. Set up and dividing into cases
\endhead
\definition{Set up (3.0)}
Let $X$ be an algebraic $\Bbb Q$-factorial normal $3$-fold,
$f:X\to Y$ is a projective small contraction onto a normal $3$-fold
$Y$ such that $\rho(X/Y) = 1$ and $C$ the exceptional locus.
Let $S$ be an irreducible divisor on $X$ and $B$ a boundary on $X$ such that
$\llcorner B \lrcorner = 0$. 
We assume all the following assumptions in the below sections:
\roster
\item $\lda{X}{S}$ is PLT. $\ldb{X}{S}{B}$ is LC and LT outside $S$.
But $\ldb{X}{S}{B}$ is not LT;
\item $(\lda{X}{S}).C <0$, $S.C < 0$ and $(\ldb{X}{S}{B}).C =0$;
\item $2(\ldb{X}{S}{B}) \sim 0$ near $C$;
\item there are at least $2$ divisors in $\text{PLC} (\lda{S}{\diff{S}(B)})$;
\item there is a unique component $L$ of 
$\llcorner \diff{S}(B) \lrcorner$ such that
$L \not \subset C$.
For any component $D$ of $\supp B|_S$, 
if $D \cap \llc (\lda{S}{\diff{S}(B)}) \not = \phi$ 
then $D=L$ or $D \subset C$.
\comment
We see that for a connected component $C'$ of $C$, $C' \cap L$ is empty
or consists of one point by applying (1.6) for $f|_S$
\endcomment
\endroster
As we reviewed in the introduction,
we have only to prove the existence of the flip
under the assumptions 

(1)' $K_X +S$ is PLT but $K_X +S+B$ is not LT

and (2) and (3).
But by 
[S, Proposition 6.12] or [FA, Proposition 21.9], we may assume (1),
by[S, Corollary 7.3] or [FA, 21.12.1], we may assume (4)
and by [S, Reduction 8.2] or [FA, 22.1], we may assume (5).
\enddefinition

\proclaim{Main Theorem}
Under the assumptions in (3.0), the flip of $f$ exists.
\endproclaim

\definition{(3.1) Reduction to the case that $C$ is irreducible}
Here we consider $f$ as in (3.0) but only assume (1)$\sim$ (3).
Replace $X$ by an analytic neighborhood of $C$ and take 
a $\Bbb Q$-factorial LT model $\t{X} \to X$ for $K_X + S$ near $S$
which exists by (1.13). Let $\t{S}$ be the strict transform of $S$.
Run the $K_{\t{X}}+\t{S}$-LMMP. Then the intermediate flips
satisfy (1)$\sim$(3) in (3.0).
Hence by taking algebraization,
it suffices to prove the main theorem under
the additional assumption that $C$ is irreducible
(note that we can assume also (4) and (5).)
\enddefinition

\definition{Dividing into cases (3.2)}
If $C\cap L$ consists of more than $1$ points,
by (3.1), we may assume $L$ is not irreducible.
But in this case the flip exists by (5) of (3.0).
Hence we may assume that $C\cap L$ is one point and we denote it by
$P$.

Note that if $C \not \subset \supp B$,
$\ldb{X}{S}{B}$ is PLT at the generic point of $C$
since $\lda{X}{S}$ is PLT. Hence the following $4$ cases occur.
\enddefinition

\definition{(a)}
$C$ is irreducible.
$C \not \subset \supp B$. $\ldb{X}{S}{B}$ is PLT at the generic point of $C$.

By applying (1.6) for $f|_S$ and the assumption (5)
in (3.0),
$\clc{0}(\lda{S}{\diff{S}(B)}) = \{ C \cap L \}$.
\enddefinition

\definition{(b)}
$C$ is irreducible.
$C \subset \supp B$. $\ldb{X}{S}{B}$ is not PLT at the generic point of $C$.
\enddefinition
 
In this case we furthermore divide into $2$ cases.

\definition{(b1)}
There is a point $Q \in C \ \text{and} \ \not \in L$ such that $Q \in 
\clc{0}(\lda{S}{\diff{S}(B)})$.

We see that such a $Q$ is unique by considering the degree of the boundary
of $\ldb{X}{S}{B}|_S|_C$.
\enddefinition
 
\definition{(b2)}
The contrary case to (b1).
\enddefinition

\definition{(c)}
$C \subset \supp B$. $\ldb{X}{S}{B}$ is PLT at the generic point of $C$.
By applying (1.6) for $f|_S$ and the assumption (5)
in (3.0),
$\clc{0}(\lda{S}{\diff{S}(B)}) = \{ C \cap L \}$.

In this case we do not assume that $C$ is irreducible for induction.
But we may assume that all the components intersect $L$.
In fact a component which does not intersect $L$ is analytically of type $3$
(see (2.2)).
Hence by (1.4),
the number of irreducible components is at most $2$.
\enddefinition

\definition{Definition (3.3) (Good extraction)}
Consider the situation as in (3.0).
Let $g:Z\to X$ be a log crepant divisorial (primitive)
extraction for $\ldb{X}{S}{B}$ and
$E$ the exceptional divisor.
We say that $g$ is a good extraction (resp. semi-good extraction)
if and only if the following hold:
\roster
\item $D:= \-{g}S\cap E\simeq \Bbb P^1$ 
(resp. $\-{g}S\cap E$ is 
a chain of $D_1\simeq \Bbb P^1$ and a smooth curve $D_2$);
\item $\ldb{Z}{\-{g}S}{E}$ is LT 
(resp.$\ldb{Z}{\-{g}S}{E}$ is LT outside $P$);
\item $\ldc{Z}{\-{g}S}{E}{B}|_D=\ldc{D}{P}{\frac12 P_1}{\frac12 P_2}$, 
where $P$, $P_1$ and $P_2$ are distinct points.
(resp. $\ldc{Z}{\-{g}S}{E}{B}|_{D_1}=\ldc{D_1}{P}{\frac12 P_1}{\frac12 P_2}$, 
where $P:= D_1\cap D_2$, $P_1$ and $P_2$ are distinct points).
\endroster
\enddefinition

\proclaim{Corollary for the Definition (3.4)}
Consider the situation as in Definition (3.3).
Assume that $g(E)$ is a point.
Then 
if $g$ is good (resp. $g$ is semi-good),
$$\llc(\ldc{Z}{\-{g}S}{\-{g}B}{E}|_E)=D \ \text{or} \ D\cup M,$$
where $M$ is an irreducible curve and intersects $D$ only at $P$ 
$$\text{(resp.} \ \llc(\ldc{Z}{\-{g}S}{\-{g}B}{E}|_E)=
D_1\cup D_2 \ \text{or} \ 
D_1\cup D_2\cup M,$$  
where $M$ is an irreducible curve and intersects $D_2$ only at one point
$\not =P$).

Furthermore on $M$ 
there is at most one $\clc{0}(\ldc{Z}{\-{g}S}{\-{g}B}{E}|_E)$
except $P$ (resp. except $D_2 \cap M$).
\endproclaim
\demo{proof}
We prove the claim only for the case $g$ is good (the proof for another case
is similar). By (1.7) and 
(3) of the definition of good extraction in (3.3),
$\llc(\ldc{Z}{\-{g}S}{\-{g}B}{E}|_E)$ is a union of curves. 
Assume that
$\llc(\ldc{Z}{\-{g}S}{\-{g}B}{E}|_E)$ has 
another irreducible component except $D$. 
Let $M$ be one of it.
Since $D = \-{g}S|_E$ is ample and
$\clc{0}(\ldc{Z}{\-{g}S}{\-{g}B}{E}|_E)\cap D=P$, $M$ must intersect
$D$ only at $P$.
By this, $M$ must be unique.
For the latter half, we have only to note
$\text{deg} \ \ldc{Z}{\-{g}S}{\-{g}B}{E}|_E|_M=0$.
\qed
\enddemo

\head 4. Treatment of case $(a)$ and $(b1)$
\endhead
In this section we consider cases (a) and (b1).
We use the notation as in (3.2).
\proclaim{Theorem (4.0)}
Let $h_t:Z_t\to X$ be 
a $\Bbb Q$-factorial LT model for $\ldb{X}{S}{B}$ near $S$
(which exists by (1.13)).
 
Run the $(\ldc{Z_t}{S_t}{B_t}{E_t-\epsilon B_t})$-MMP over $X$.

Then by interchanging the surgeries of MMP 
if necessarily, there is a nonempty sequence
of divisorial extractions 
$$Z_n\overset {g_n}\to\longrightarrow 
Z_{n-1} \overset {g_{n-1}}\to \longrightarrow \dots
Z_1 \overset {g_1}\to \longrightarrow X$$ satisfying the following properties:
\roster 
\item Let $E_i$ be the exceptional divisor of $g_i$.
Then $Q_{i-1}:=g_i(E_i)$ is a point for $i\geq 2$
and $Q=g(E_1)$. 
Furthermore $Q_{i-1}\in E_{i-1}\ \text{but} \ \not \in \cup_{j\not
=i-1} E_j$;
\item $\ldb{Z_i}{S_i}{\sum_{j=1}^{i} E_j}$ is LT,
where $S_i$ is the strict transform of $S$ on $Z_i$;
\item Let $D_i:=E_i\cap E_{i-1}$ for $i\geq 2$ and $D_1:=E_1\cap S_1$.
Then for $1\leq i\leq n-1$ (resp. $i=n$),
$\clc{0}(\ldc{Z_i}{S_i}{B_i}{E_i}|_{E_i})=\{Q_i, P_i \}$ (resp. $=\{P_i\}$),
and 
$\llc(\ldc{Z_i}{S_i}{B_i}{E_i}|_{E_i})=D_i\cup L_i$, where
$L_i:=\supp B_i|_{E_i}$ and $P_i:=L_i\cap D_i$.
Furthermore $L_i$ is irreducible and $Q_i\in L_i$.
In particular $g_1$ is a good extraction.
\endroster
(See FIGURE (VII).)
\endproclaim

\demo{Proof}
As long as a component of $E_t$ is left, $-\epsilon B_t$ is not nef by
[FA, 2.19 Lemma]
so the MMP continue till every component of $E_t$ is contracted.
After contracting every component of $E_t$, the morphism going back to
$X$ is at most a small contraction but since $X$ is $\Bbb Q$-factorial,
it must be isomorphism.
Hence the last morphism is a divisorial contraction to $X$.
Let it $g_1:Z_1\to X$ and $E_1$ the exceptional divisor.

\proclaim{Claim 1}
We can assume that 
$g_1(E_1)=Q$ by interchanging the surgeries of MMP if necessarily. 
\endproclaim

\demo{Proof}

\definition{case (a)}

If $g_1(E_1)\not = Q$, we must have $g_1(E_1)=L$.
But we prove that it is impossible by the nonexceptional assumption.
Assume that $g_1(E_1)=L$.
Since $\ldb{Z}{S_1}{E_1}$ is LT by the choice of MMP, 
$S_1\cap E_1$ is irreducible.
Hence $g'_1:=g_1|_{S_1}: S_1 \to S$ is isomorphism.
By subadjunction, 
$$\lda{S_1}{\diff{S_1}(E_1+B_1)}
=\pa{{g'}_1}{S}{\diff{S}(B)}.
\tag 4.1.1$$
Let $L'$ be the divisor corresponding to $L$ by the isomorphism
$S_1\simeq S$.
Then $L'=E_1\cap S_1$.
Hence $L'\not \subset \supp B_1$ by LC property of
$\ldc{Z}{S_1}{B_1}{E_1}$ along $L'$.
On the other hand in case (a), 
$L=\supp B|_S$ so $\supp B_1 \cap S_1$ is at most a set of finite points.
But $Z_1$ is also $\Bbb Q$-factorial, it must be empty.
Hence $\ldc{Z_1}{S_1}{B_1}{E_1}$ 
(and also $\ldc{Z_1}{S_1}{B_1}{E_1}|_{S_1}$)
is LT.
Consequently by $(4.1.1)$,
$\lda{S}{\diff{S}(B)}$ is also LT.
But since $C\not \subset \llc(\lda{S}{\diff{S}{B}})$
it must be PLT, a contradiction to the nonexceptional assumption.
\enddefinition

\definition{Case (b1)}

If $g_1(E_1)\not = Q$, we must have $g_1(E_1)=C, L \ \text{or} \ C\cap L$.

\enddefinition

We treat them separately.

\definition{Subcase 1. $g_1(E_1)=C$}
We exclude this case as in Case (a).
\enddefinition

\definition{Subcase 2. $g_1(E_1)=L$ or $C\cap L$}
 
If this case occurs, we prove the statement of Claim 1.
Let $$Z'_{i+1}\overset {g_{i+1}}\to\longrightarrow 
Z'_{i} \dots \dashrightarrow \dots
Z_1 \overset {g_1}\to \longrightarrow Z$$ be the sequence of the
surgeries of MMP such that $Z'_{i}\dashrightarrow X$
consists of flips or divisorial extractions whose centers
are the strict transform of $L$ or $C\cap L$ and
$Z'_{i+1}\to Z'_i$ is a divisorial extraction whose center is $Q$ 
or the strict transform of $C$.

Note that $Z'_i\dashrightarrow  X$ is isomorphic near $Q$.
Then by the similar argument to Case (a),
we can prove that the center of $Z'_{i+1}\to Z'_i$ is not $C$.
Hence we can interchange $Z'_{i+1}\to Z'_i$ and $Z_1 \to X$.
\enddefinition
\qed
\enddemo

We check (1), (2) and (3) of the statement of (4.0)
for $g_1:Z_1 \to X$ as in Claim 1.
(1) is proved by Claim 1 and (2) is clear by the choice of MMP.

Let $C'$ and $L'$ be the strict transforms of $C$ and $L$ on $S_1$
respectively. Let $P_1:=L'\cap D_1$ in case (a) 
(resp. $P_1:=C'\cap D_1$ in case (b1)).

In case (a), $L=\supp B|_S$ and by LC property of 
$\ldc{Z_1}{S_1}{B_1}{E_1}$ along $D_1$, $D_1\not \subset \supp B_1$.
Hence $L'=\supp B_1|_{S_1}$.
So $\ldc{Z_1}{S_1}{B_1}{E_1}$ is LT at points of $D_1$
except $P_1$,
i.e., $P_1=\clc{0}(\ldc{Z_1}{S_1}{B_1}{E_1}|_{E_1})\cap D_1$.

In case (b1), near $Q$, $C=\supp B|_S$.
Hence by the similar argument to case (a),
$P_1=\clc{0}(\ldc{Z_1}{S_1}{B_1}{E_1}|_{E_1})\cap D_1$.

If there is no point in $\clc{0}(\ldc{Z_1}{S_1}{B_1}{E_1}|_{E_1})$
except $P_1$, by setting $n=1$, we finish the poof of (4.0).

Assume that there is a point $Q_1$ in
$\clc{0}(\ldc{Z_1}{S_1}{B_1}{E_1}|_{E_1})$ other than $P_1$.
We first show the uniqueness of $Q_1$.
By (1.7) and the existence of $P_1$,
$\llc(\ldc{Z_1}{S_1}{B_1}{E_1}|_{E_1})$ is connected.
Hence $P_1, Q_1 \in L_1$, where 
$L_1:=\llc(\ldc{Z_1}{S_1}{B_1}{E_1}|_{E_1})-D_1$.
\proclaim{Claim 2}
$L_1$ irreducible.
\endproclaim
\demo{Proof}
If $L_1$ is not irreducible, by $\rho(Z_1/X)=1$ and $g_1(E_1)=Q$,
$D_1$ is ample.
So any component of $L_1$ intersects with $D_1$.
By just above, any component of $L_1$ passes through $P_1$.
But this is a contradiction since there is at least three reduced boundary
of 
$\ldc{Z_1}{S_1}{B_1}{E_1}|_{E_1}$ at $P_1$.
\qed
\enddemo

Consider the degree of $\ldc{Z_1}{S_1}{B_1}{E_1}|_{E_1}|_{L_1}$.
$P_1$ and $Q_1$ contribute degree $1$ respectively to its boundary.
By this, we know the uniqueness of $Q_1$.

Next Claim 3 complete the proof of (3) for $g_1$.

\proclaim{Claim 3}
$L_1=\supp B_1|_{E_1}$
\endproclaim

\demo{Proof}
Assume that there is another component $L'_1$ of $\supp B_1|_{E_1}$ 
except $L_1$.
Since $D_1$ is ample, $D_1\cap L'_1\not = \phi$.
As above, $L'_1$ cannot pass through $P_1$ by LC
property of $\ldc{Z_1}{S_1}{B_1}{E_1}|_{E_1}$ at $P_1$.
This means that there is a component of $\supp B_1|_{S_1}$ 
except $L'$ in case (a) (resp. $C'$ in case (b1)).
But this contradicts that $L=\supp B|_S$ in case (a) 
(resp. $C=\supp B|_S$ near $Q$ in case (b1)).
\qed
\enddemo

Next we prove the existence of the sequence by induction.
Assume that we can construct the sequence satisfying the statement of
(4.0) until i-th stage:
$$Z_i\to Z_{i-1}\to \dots \to X.$$
We decompose the surgeries before $Z_i$ as follow:
$$Z_{j+1}\overset {g_{j+1}} \to \dashrightarrow Z_j \dashrightarrow 
\dots Z_{i+1} \to Z_i,$$
where $Z_j\dashrightarrow Z_i$ consists of flips or divisorial contractions 
which are isomorphic near $Q_i$ and $Z_{j+1}\dashrightarrow Z_j$ is
not isomorphic near $Q_i$.

\proclaim{Claim 4}
$Z_{j+1} \dashrightarrow Z_j$ is a divisorial contraction whose center 
is $Q_i$.
\endproclaim
\demo{Proof}
Assume the contrary. Then $Z_{j+1}\dashrightarrow Z_j$ is a divisorial 
contraction of (the strict transform of) $L_i$
or flip.

First assume that $Z_{j+1}\dashrightarrow Z_j$ is a divisorial 
contraction of (the strict transform of) $L_i$.
Here we denote the strict transform of $E_i$ on $Z_j$ by $E_i$
and the strict transform of $E_i$ on $Z_{j+1}$ by $E'_i$.
The argument is almost the same as one in Claim 1.
Since $\ldb{Z}{E'_i}{E_{j+1}}$ is LT by the choice of MMP, 
$E'_i\cap E_{j+1}$ is irreducible.
Hence $g'_{j+1}:=g_{j+1}|_{E'_i}: E'_i \to E_i$ is isomorphism.
By subadjunction, 
$$\lda{E'_i}{\diff{E'_i}(S_{j+1}+B_{j+1}+ E_{j+1} + E')}
=\pa{{g'}_1}{E_i}{\diff{E_i}(S_j+B_j+ E)},
\tag 4.1.2$$
where 
$E$ is the reduced sum of exceptional divisors of $Z_j \dashrightarrow X$ 
except $E_i$ 
and $E'$ is the strict transform of $E$.
Let $L'_i$ be the divisor corresponding to $L_i$ by the isomorphism
$E'_i\simeq E_i$.
Then $L'_i=E'_i\cap E_{j+1}$.
Hence $L'_i\not \subset \supp B_{j+1}$.
On the other hand, 
$L_i=\supp B_j|_{E_i}$ near $Q_i$ by the assumption of induction (3).
Hence $\supp B_{j+1} \cap E'_i$ must be empty near $Q_i$.
So $Q_i$ is not contained in 
$\clc{0}(\lda{E'_i}
{\diff{E'_i}(S_{j+1}+B_{j+1}+ E_{j+1} + E')})$,
a contradiction to the assumption of the induction (3).

Next assume that $g_{j+1}$ is a flip.
By assumption, there is a connected component of the flipped curve
through $Q_i$ and we call it $l$.
Since $Z_j\dashrightarrow Z_i$ is isomorphic near $Q_i$,
(the strict transform of) $l$ on $Z_i$ passes through $Q_i$.
On $Z_i$, $l$ is a exceptional curve of $Z_i\to X$ so
on $Z_j$, $l\subset E_i$ and $l$ is not contained in another
exceptional divisor of $h_j:Z_j\to X$ or $S_j$.
On the other hand we have ${h_j}^* B.l=0$ and $B_j.l<0$.
Hence there exists a $h_j$-exceptional divisor $E_k$ such that $E_k.l>0$.
We also have $h_j^*S.l=0$ so there is a divisor $F$ which is $S_j$ or 
$h_j$-exceptional such that $F.l<0$  
Hence $F$ must be $E_i$ and in particular $E_i.l<0$.
Then we see that $l=L_i$ and $l=\supp B_j|_{E_i}$ near $l$ by the following 
argument:

$l=L_i$ follows from the fact that 
near $Q_i$, $L_i=\supp B_j|_{E_i}$ and $l\subset \supp B_j$.
For the latter half, we have only to consider the degree of
$\lda{E_i}{\diff{E_i}(S_j+B_j+ E)}|_l$.

Here we consider the situation before the flip.
We use the notation as in the FIGURE (VI).
By the above numerical consideration, $E_i^-.l^->0$ and $E_k^-.l^-<0$.
Then any component of $l^-$ is contained in $E_i^-$.
For if a component $l_1^-$ of $l^-$ is contained in $E_i^-$,
$l_1^- = E_i^-\cap E_k^-$ since $\ldb{Z_{j+1}}{E_i^-}{E_k^-}$ is LT
and hence $E_i^-\cap E_k^-$ is irreducible.
Since $E_k^-$ is normal and $l_1^-$ is contractible in $E_k^-$ , 
we have $E_i^-.l^-=(l_1^-)^2_{E_k^-}<0$, a contradiction.
By this, the natural map $E_i\dashrightarrow E_i^-$ is a morphism.
Let this morphism be $c:E_i\to E_i^-$.
Here we consider the following subadjunction:
$$\lda{E_i}{\diff{E_i}(S_{j}+B_{j}+ E)}
=\pa{c}{E_i^-}
{\diff{E_i^-}{(S_{j+1}+B_{j+1}+ E^-)}},
\tag 4.1.3$$
where $E^-$ is the strict transform of $E$.
Note that $l=L_i$ is a log crepant curve for this.
Since $l=\supp B_j|_{E_i}$ near $l$, 
$\supp B_{j+1}|_{E_i^-}$ is at most a point $c(l)$ near $c(l)$.
But by the $\Bbb Q$-factoriality of $Z_{j+1}$, $\supp
B_{j+1}|_{E_i^-}$ is in fact empty near $c(l)$.
Hence near $c(l)$, 
$$\lda{E_i^-}{\diff{E_i^-}{(S_{j+1}+B_{j+1}+ E^-)}}$$
is LT.
On the other hand there is a log crepant curve for this log divisor,
so there must be $F^-$ which is $S_{j+1}$ or $E_l^-$ for some $l \not = k$
such that $c(l)$ is contained in it and $c(l)$ is a smooth point of $E_i^-$.
Hence $l$ intersects $F$ and $E_k$ outside $Q_i$.
But this gives a contradiction by considering the degree of the boundary of
$\lda{E_i}{\diff{E_i}(S_{j}+B_{j}+ E)}|_l$.
\qed
\enddemo
By replacing $Z_{i+1}\dashrightarrow Z_i$ by $Z_{i+1}\dashrightarrow Z_i$,
we can assume that $g_i:Z_{i+1}\to Z_i$ is a divisorial extraction 
whose center is $Q_i$. For this $g_i$, (2) is clear by the choice of MMP
and the check of (3) is similar to the check for $g_1$.
\qed
\enddemo

\definition{Definition (4.1) (Good sequence, Invariant $\lambda$)}
We call a good sequence 
a sequence of divisorial extractions as in (4.0).
We define  
$$\multline
\lambda := \\
 \min \{n| n \ \text{is the length of a good sequence
for a} \
\Bbb Q  \ \text{-factorial LT model for} 
\ldb{X}{S}{B}\}.
\endmultline$$ 
\enddefinition

We remark that it is a well defined finite number by (4.0).

\proclaim{ (4.2) Conclusion of case (a, $\lambda$) and (b1, $\lambda$)}
Fix a good sequence whose length is $\lambda$.
We apply BTM starting by $g_1$. So we will use the notation as in BTM.
We used the same notation in (2.5) and (4.0) for different objects
but in this conclusion (4.2), the notation in (4.0) is not used except
$g_1 : Z_1 \to X$ and the notation by which we denote objects on $Z_1$.
So no confusion is caused because the same notation in (2.5) and (4.0)
stands for the same object on $Z_1$.

Our conclusion is explained in the following flow chart:  
\endproclaim
\newpage

\hfil Case (a) \hfil 
\hskip1in
$$
\fbox{2.5in}{
Choose a good sequence with length $n$ and run BTM 
starting by $g_1$}
$$

\vskip1pt
$$\downarrow$$

\hskip1in
$$
\fbox{2.8in}
{$R_1$ is a beginning flipping ray. So the flip of $R_1$ exists}
$$

\vskip1pt

$$ 
\oversetbrace\to{
\hbox{\hsize=4.4in
\fbox{1.1in}{$R_2$ is a divisorial ray. \break
DONE}\hskip5pt
\fbox{1.1in}{$R_2$ is a flipping \break 
ray of type 2. \break 
DONE}
\hskip5pt
\fbox{1.1in}{$R_2$ is a flipping 
\break
ray and $S_2 \cap C_2 = \phi$.}
\hskip5pt
\fbox{1.1in}{$R_2$ is a final \break
flipping ray. 
\break
DONE}
}}
$$

\vskip0.5in

$$
\hskip0.7in
\oversetbrace\to{
\hbox{\hsize=1.5in
\fbox{0.7in}{$P \not \in C_1$.}
\hskip1.7in
\fbox{0.7in}{$P \in C_1$.}
}
}$$
\vskip1pt

$$
\oversetbrace\to{
\hbox{\hsize=2.4in
\fbox{1.2in}{
$C_2 \cap \clc{0} \not = \phi$.
Then the connected component $C'_2$ of $C_2$ containing
$\clc{0}$ is of \break
type (a, $\lambda ' < \lambda$).\break 
(Another \break
component \break
is of type 3 or 4. \break
DONE)
}
\hskip.1in
\fbox{1.0in}{$C_2 \cap \clc{0}= \phi$.
Then $C_2$ is of type 3 or 4. \break
DONE}
}}
\oversetbrace\to{
\hbox{\hsize=2.4in
\fbox{1.2in}{A component of \break
$C_1^+$ is contained in LLC.
Then $C_1^+$ is \break
irreducible. }
\hskip.1in
\fbox{1.1in}{No component \break
of $C_1^+$ is contained in LLC.
Then $C_2$ is of type 3. \break
DONE}}}
$$
\vskip0.5in

$$\hskip.6in
\oversetbrace\to{
\hbox{\hsize=2.4in
\fbox{1.0in}{$C_2 \cap \clc{0}= \phi$. Then $C_2$ is of type 3 or 4.\break
DONE}
\hskip.1in
\fbox{1.1in}{$C_2 \cap \clc{0}\not = \phi$}}}
$$
\vskip0.5in
$$\Downarrow \hskip1.0in$$
\hskip 3.2in

$$
\fbox{1.0in}{Reduction to \break
(a, $\lambda ' < \lambda$)}
\hskip1.0in
\oversetbrace\to{
\hbox{\hsize=2.4in
\fbox{1.1in}{$L'_1$ exists. Then the connected \break
component 
$C'_2$ of \break
$C_2$ containing $e$ is of type (b). 
}
\hskip.1in
\fbox{1.1in}{$L'_1$ does not exist. Then the connected component 
$C'_2$ of \break
$C_2$ containing $e$ is of type (c)}}}
$$

\vskip1in

\hskip1.5in 
$$\hskip1.3in\Downarrow \hskip 1.5in\Downarrow$$

\hskip3.3in
 
$$
\hskip1.7in
\hbox{\hsize=2.4in
\fbox{1.1in}{Reduction to \break
case (b)}
\hskip0.4in
\fbox{1.1in}{Reduction to \break
case (c)}}
$$
\vskip0.5in
$$\clc{0}:=\clc{0}(\ldc{Z_2}{S_2}{B_2}{E_2}|_{E_2})$$
$$\llc:=\llc (\ldc{Z_2}{S_2}{B_2}{E_2}|_{E_2})$$
$$P:= L' \cap S_1 \cap E_1$$
\newpage

\hfil Case (b1) \hfil 

$$
\hskip1in
\fbox{2.5in}{
Choose a good sequence with length $n$ and run BTM 
starting by $g_1$}
$$
\vskip1pt
$$\downarrow$$
$$
\hskip1in
\fbox{2.8in}
{$R_1$ is a beginning flipping ray. So the flip of $R_1$ exists}
$$
\vskip1pt
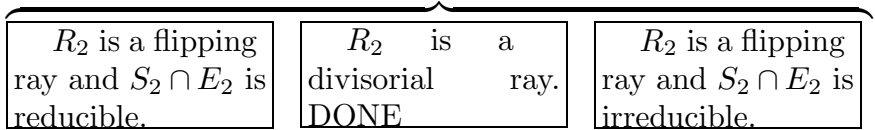
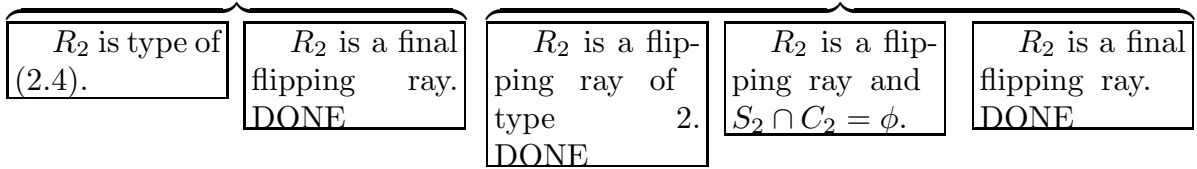
$$
\oversetbrace\to{
\hbox{\hsize=4.4in
\fbox{1.1in}{$R_2$ is a flipping \break
ray and $S_2 \cap E_2$ is reducible.}
\hskip5pt
\fbox{1.1in}
{$R_2$ is a \break
divisorial ray.\break
 DONE}
\hskip5pt
\fbox{1.1in}
{$R_2$ is a flipping \break
ray and $S_2 \cap E_2$ is irreducible.}
}
}
$$
\vskip 0.5in
$$
\oversetbrace\to{
\hbox{\hsize=1.8in
\fbox{0.9in}{$R_2$ is type of (2.4).}\hskip5pt
\fbox{0.9in}{$R_2$ is a final flipping ray.\break
 DONE}}}\hskip5pt
\oversetbrace\to{
\hbox{\hsize=2.7in
\fbox{0.9in}{$R_2$ is a flipping ray of \break
type 2.\break
 DONE}\hskip5pt
\fbox{0.9in}{$R_2$ is a flipping ray and \break
 $S_2 \cap C_2 = \phi$.}
\hskip5pt
\fbox{0.9in}{$R_2$ is a final flipping ray. \break
DONE}}}
$$
\vskip0.5in

$$\downarrow \hskip3.6in$$

$$
\fbox{1.1in}{$R_3$ is a flipping \break
ray of type 2. \break
DONE}
\hskip1.0in
\oversetbrace\to{
\hbox{\hsize=4in
\fbox{1.2in}{A component of \break
$C_1^+$ is contained in LLC.}
\hskip1.0in
\fbox{1.2in}{No component of $C_1^+$ is contained in LLC. 
Then $C_2$ is of type 3. \break
DONE}
}}
$$
\vskip0.5in
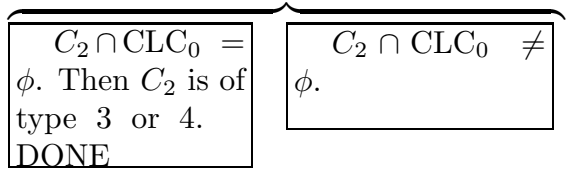
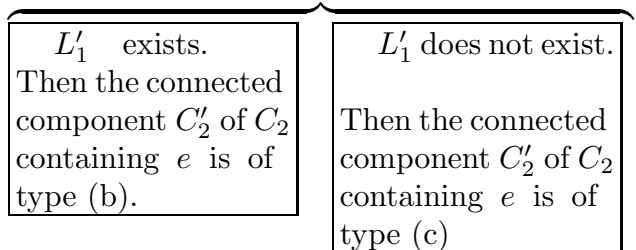
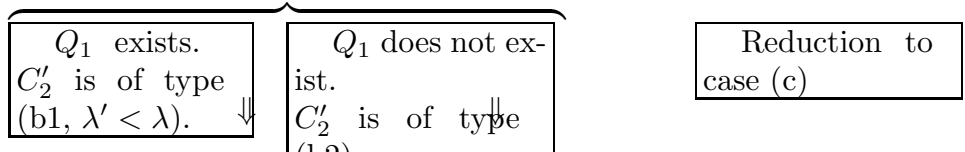
$$
\oversetbrace\to{
\hbox{\hsize=2.1in
\fbox{1.0in}{$C_2 \cap \clc{0}= \phi$.
Then $C_2$ is of type 3 or 4. 
\break
DONE}
\hskip.1in
\fbox{1.1in}{
$C_2 \cap \clc{0} \not = \phi$. \break 
}
}}
$$
\vskip0.7in
$$
\hskip2.0in\oversetbrace\to{
\hbox{\hsize=2.4in
\fbox{1.2in}{$L'_1$ exists. \break
Then the connected \break
component 
$C'_2$ of 
$C_2$ containing $e$ is of \break
type (b). 
}
\hskip.1in
\fbox{1.2in}{$L'_1$ does not exist. \break
Then the connected \break
component 
$C'_2$ of $C_2$ containing $e$ is of \break
type (c)}
}}
$$
\vskip0.7in
$$\hskip3.5in \Downarrow$$
$$
\hskip1.0in
\oversetbrace\to{
\hbox{\hsize=2.1in
\fbox{1.0in}{$Q_1$ exists. \break
$C'_2$ is of type \break
(b1, $\lambda ' < \lambda$).}
\hskip.1in
\fbox{1.1in}{
$Q_1$ does not exist. \break
$C'_2$ is of type \break
(b2).
}
}}
\hskip0.5in
\hbox{
\fbox{1.1in}{Reduction to \break
case (c)}}
$$

$$\Downarrow \hskip1.0in \Downarrow$$

$$
{
\hbox{\hsize=2.4in
\fbox{1.1in}{Reduction to \break
case (b1, $\lambda ' < \lambda$)}
\fbox{1.1in}{Reduction to \break
case (b2)}}}
$$
\vskip0.5in
$$\clc{0}:=\clc{0}(\ldc{Z_2}{S_2}{B_2}{E_2}|_{E_2})$$
$$e \ \text{is the unique element in} \ \clc{0} \text{outside} \ S_2$$
$$\llc:=\llc (\ldc{Z_2}{S_2}{B_2}{E_2}|_{E_2})$$
\newpage

\proclaim{Case (a)}
\endproclaim

Since $g(E)$ is a point, $\supp R_1$ is the strict
transform of $C$. Clearly we have $E.R_1>0$. Hence by $f^*S.R_1<0$,
$S_1.R_1<0$. So $R_1$ is a beginning flipping ray and hence the flip
exists.

\proclaim{Claim 1}
\roster
\item $S_2\cap E_2$ is irreducible;
\item $\lda{Z_2}{E_2}$ is PLT.
\endroster
\endproclaim
\demo{Proof}
\roster
\item
First assume that $P\not \in C_1$.
If $S_2\cap E_2$ is reducible,
a component $\o{{C_1}^+}$
of $C^+_1$ is contained in $S_2\cap E_2$.
Then $\o{C^+_1}\in \llc (\ldc{Z_2}{S_2}{B_2}{E_2})$,
which in turn shows that $C_1$ contains an element of  
$\text{CLC} \ (\ldc{Z_1}{S_1}{B_1}{E_1})$ since $\ldc{Z_1}{S_1}{B_1}{E_1}$
is numerically trivial over $Y$.
But this contradicts the assumption of Subcase 1.

Next assume that $P \in C_1$.
If $S_2\cap E_2$ is reducible,
a component $\o{C^+_1}$
of $C^+_1$ is contained in $S_2\cap E_2$.
On the other hand by $P\in C_1$, $B_1.C_1>0$ and hence $B_2.C^+_1<0$.
So $\o{C^+_1} \subset S_2\cap E_2 \cap \supp B_2$.
But this contradicts LC property of $\ldc{Z_2}{S_2}{B_2}{E_2}$ along 
$\o{C^+_1}$.
\item
By the proof of (1), $E_2$ is normal and
$\llc (\ldc{Z_2}{S_2}{B_2}{E_2}|_{E_2}) \subset \supp {B_2}|_{E_2} \cup
\supp {S_2}|_{E_2}$.
Hence $\lda{Z_2}{E_2}|_{E_2}$ is PLT, which in turn shows that 
$\lda{Z_2}{E_2}$ is PLT by the inversion of adjunction.
\endroster
\qed
\enddemo	

Consider the next extremal ray $R_2$.
$R_2$ may be a divisorial ray but in this case we can obtain
the flip of $C$ by contracting $E_2$.
Assume that $R_2$ is a flipping ray.
Since $S_2\cap E_2$ is irreducible and 
all exceptional curves over $Y$ are contained in $E_2$,
$R_2$ is a middle flipping ray or a final flipping ray.
we may consider only the case that $E_2.R_2 < 0$ and 
$\supp R_2 \cap S_2 = \phi$
(see the last remark in (2.5)).
Note that by the numerical property (2) in (2.5), $B_2.R_2 < 0$.
Hence $R_2$ is a special flipping ray by Claim 1 (2).
Here we divide into two cases.

\definition{Subcase 1. $P\not \in C_1$}
See FIGURE (VIII).
\enddefinition

\proclaim{Claim 2}
$C_1^+$ is irreducible.
\endproclaim
\demo{Proof}
Easy.
\qed
\enddemo

Note that $\llc (\ldc{Z_1}{S_1}{B_1}{E_1}|_{E_1}) = S_1 |_{E_1}
\ \text{or} \ S_1 | _{E_1} \cup L_1$.
By the proof of Claim 1,   
$\llc (\ldc{Z_2}{S_2}{B_2}{E_2}|_{E_2}) = S_2 |_{E_2}
\ \text{or} \ S_2 | _{E_2} \cup L'_1$, where $L'_1$ is the strict transform
of $L_1$. In particular $C_2 \not \subset 
\llc (\ldc{Z_2}{S_2}{B_2}{E_2}|_{E_2})$.

If $C_2 \cap \clc{0}(\ldc{Z_2}{S_2}{B_2}{E_2}|_{E_2}) = \phi$,
$\ldb{Z_2}{E_2}{B_2}|_{E_2}$ is PLT or exceptional near $C_2$.
If the former case occurs, $\ldb{Z_2}{E_2}{B_2}$ is PLT near $C_2$
by the inversion of adjunction 
and hence $C_2$ is of type 3 and the flip exists by (2.2).
If the latter case occurs, $C_2$ is of type 4 and the flip exists by (2.3). 

Hence we assume below that
$C_2\cap \clc{0}(\ldc{Z_2}{S_2}{E_2}{B_2}) \not =\phi$.
By the proof of Claim 1,
we know that there is no $\clc{0}(\ldc{Z_2}{S_2}{E_2}{B_2})$ on $C_1^+$.
Hence $\clc{0}(\ldc{Z_2}{S_2}{E_2}{B_2}) = \{Q_1, P \}$.
By the assumption that $C_2 \cap S_2 = \phi$, $Q_1 \in C_2$. 
Let $C'_2$ be the connected component of $C_2$ containing $Q_1$.
We will show that the flip of $C'_2$ is of type $(a, \lambda ' < \lambda)$.
See FIGURE (IX).
(Another component is of type 3 or 4. So we have only to consider $C'_2$.)
First we prove that $C'_2$ is irreducible.
Denote the strict transform of $L_1$ on $E_2$ by $L'_1$.
Note that $L'_1.C'_2>0$ and $L'_1.C_1^+=0$.
Since $\rho(Z_1/Y)=2$ and $L'_1=\supp B_2|_{E_2}$, 
only components of $C_1^+$ are numerically trivial for
$L'_1$ among curves in $E_2$. By similar reason, 
we see that $(S_2|_{E_2})^2 >0$. 
Note also that every component of $C'_2$ intersects $C_1^+$
because $C_1^+$ is irreducible and before the flip of $C_1$,
every component of $C'_2$ intersect $S_1$.
So after contracting a component of $C'_2$, 
$L'_1$ becomes ample.
By these properties, we can apply (1.9) after
contracting one component of $C'_2$
by setting $B_1 =S_2\cap E_2$, $B_2=L'_1, C_1=C'_2 \ \text{and} \  C_2=\phi$
(every left side is the notation in (1.9)).
Then $C'_2\cap L'_1 \not \in \clc{0}(\ldb{Z_2}{E_2}{B_2})$.
But by the assumption that $Q_1 \in C'_2$, 
$C'_2$ must be empty after the contraction of
a component of it, i.e., $C'_2$ is irreducible. 

We can easily check that $(Z_2, E_2, B_2, C_2)$ satisfies the assumption of 
case (a) (in this case $L'_1$ corresponds to $L$).
Note that the flip $Z_1 \dashrightarrow Z_2$ is isomorphism near $Q_1$.
Hence the above good sequence except $g_1$ is
also a good sequence with length $\lambda - 1$. 
So we finish the reduction to the case $(a, \lambda ' < \lambda)$.

\definition {Subcase 2. $P\in C_1$}
See FIGURE (X).
\enddefinition

Note that $S_2.C_2=0$ and $S_2.C_1^+ >0$.
By $\rho(Z_2/Y_1)=2$, only components of $C_2$ is numerically trivial for $S_2$
among curves in $E_2$.
 Let $a: E_2 \to E_1$ be the contraction which is the restriction of 
$Z_2 \dashrightarrow Z_1$.
Then
$C_1^+$ is the exceptional curve of $a$.

Assume that any component of $C_1^+$ is not contained in
$\llc (\ldc{Z_2}{S_2}{E_2}{B_2}|_{E_2})$. 
Since $C_1^+ \cap S_2 \ni P$, 
$C_1^+ \cap \llc (\ldc{Z_2}{S_2}{E_2}{B_2}|_{E_2}) = \{P \}$ 
by (1.6). If $L_1$ exists, 
its strict transform $L'_1$ must pass through $P$. Hence
$\ldc{Z_2}{S_2}{E_2}{B_2}|_{E_2}$ has two reduced boundaries at $P$.
But $C_1^+ \subset \supp B_2|_{E_2}$, a contradiction to 
LC property of $\ldc{Z_2}{S_2}{E_2}{B_2}|_{E_2}$. 
So $L_1$ does not exist and $\ldc{Z_2}{S_2}{E_2}{B_2}|_{E_2}$ is
PLT near $C_2$. 
Hence $\ldb{Z_2}{E_2}{B_2}$ is PLT near $C_2$ 
by the inversion of adjunction and
$C_2$ is of type 3. So the flip of $C_2$ exists by (2.2).

Assume that a component of $C_1^+$ is contained in
$\llc (\ldc{Z_2}{S_2}{E_2}{B_2}|_{E_2})$. 
By $S_2.C_1^+ >0$, all the components of $C_1^+$ pass through $P$.
Furthermore by $C_1^+\subset \supp B_2$, 
$C_1^+$ must be irreducible and furthermore $\supp B_2 |_{E_2} - C_1^+$
does not intersect $S_2$.
Then $\llc (\ldc{Z_2}{S_2}{E_2}{B_2}|_{E_2}) = 
S_2 | _{E_2} \cup C_1^+ \ \text{or} S_2 | _{E_2} \cup C_1^+ \cup L'_1$,
where $L'_1$ is the strict transform of $L_1$ (if it exists).
So if $C_2 \cap \clc{0}(\ldc{Z_2}{S_2}{B_2}{E_2}|_{E_2}) = \phi$,
$C_2$ is of type 3 or 4 as in Subcase 1.
Assume below that 
$C_2 \cap \clc{0}(\ldc{Z_2}{S_2}{B_2}{E_2}|_{E_2}) \not = \phi$.
Note that outside $S_2 |_{E_2}$, elements of
$\clc{0}(\ldc{Z_2}{S_2}{B_2}{E_2}|_{E_2})$ are only on $C_1^+$ and
in fact only one point exists on $C_1^+$ outside $S_2 |_{E_2}$.
We denote this point by $e$.
Note that if $L_1$ exists, $L_1 \subset C_2$ and $e \in L_1$.
Hence in any case $e\in C_2$. 
Let $C'_2$ be the connected component of $C_2$ containing $e$.
We will show that $\supp B_2 |_{E_2} - C_1^+ = C'_2$.

By  $\supp B_2 |_{E_2} - C_1^+ \cap S_2 = \phi$,  
$\supp B_2 |_{E_2} - C_1^+\subset C_2$.
After contracting only  $\supp B_2 |_{E_2} - C_1^+$, 
$C_1^+$ becomes ample by (1.8) and 
ampleness of $\supp B_1 |_{E_1}$ on $E_1$.
Hence we can apply (1.9) by setting
$B_1=S_2 \cap E_2$, $B_2=C_1^+$, $C_1=C'_2$ and $C_2=\phi$
after the contraction of $\supp B_2 |_{E_2} - {C^+}_1$. 
We obtain that $C'_2$ cannot pass through $\clc{0}(\ldb{Z_2}{E_2}{B_2})$
after the contraction of $\supp B_2 |_{E_2} - C_1^+$
but it means $C'_2$ becomes empty 
after the contraction of $\supp B_2 |_{E_2} - C_1^+$,
i.e., $C'_2=\supp B_2 |_{E_2} - C_1^+$. 

If $L'_1$ exists, then $L'_1 = \supp B_2 |_{E_2} - C_1^+$ and
$C'_2$ is of type (b).
(In this case $C_1^+$ corresponds to $L$ in (3.2).
Furthermore if $Q_1$ exists, $C'_2$ is of type (b, $\lambda ' < \lambda$)
or if $Q_1$ does not exist, $C' _2$ is of type (b2).)

If $L'_1$ does not exist, $C'_2$ is of type (c).
(In this case $C_1^+$ corresponds to $L$ in (3.2).)
See FIGURE (XI).
\definition{Case (b1)}
See FIGURE (XII).
\enddefinition

We know as in case (a) that $R_1$ is a beginning flipping ray
such that $\supp R_1$ is the strict transform of $C$ on $Z_1$.
We examine the next extremal ray $R_2$.
It may be a divisorial ray but in this case we obtain the flip of $C$
by contracting $E_2$. We may assume that $R_2$ is a flipping ray.
We treat by dividing into two cases.

\definition{Subcase 1. $S_2\cap E_2$ is reducible}
See FIGURE (XIII).
\enddefinition

Let $D$ be the strict transform of $S_1\cap E_1$ on $S_2$.
First we prove that $S_2 \cap E_2= D\cup l$,
where $l$ is an irreducible component of $C_1^+$.

Let $S_2 \cap E_2= D\cup l_1 \cup l_2 \dots \cup l_n$,
where $l_i$ is an irreducible component of $C_1^+$ and they form a chain
with this order.
By LC property of $\ldc{Z_2}{S_2}{B_2}{E_2}$,
$l_i\not \subset \supp B_2$.
On the other hand $l_i$'s form a chain and only $l_n$ 
intersect $\supp B_2$. Hence $n$ must be $1$.
We newly denote this $l_1$ by $l$.
Note that $B_2.D=0$. On the other hand $B_2.C_1^+>0$.
Hence $B_2.C_2 \leq 0$.
By the numerical properties (1) and (2) in (2.5),
$E_2.C_2<0$ and $S_2.C_2>0$ or $E_2.C_2>0$ and $S_2.C_2<0$.
In the former case $R_2$ is a final flipping ray so we are done.
In the latter case $C_2 \subset S_2$ 
but contractible curves in $S_2$ are only $D$ and $l$ 
and $B_2.D=0$ and $B_2.l>0$.
So in the latter case $B_2.C_2=0$.
In this case we immediately see that $D=C_2$ 
by $B_2.C_2=0$, $B_2.C_1^+ >0$ and $B_2.D=0$.  
Furthermore we can easily check 
that $D$ satisfies the assumption of (2.4).
We will apply (1.10) 
by setting $a_1= h_2|_{E_2}$, $a_2= g_3|_{E_3}\circ \mu$ and
$K_S + B = K_{X_2} + g(S_2) + g(B_2) + g(E_2)|_{g(E_2)}$, 
where $\mu$ is the normalization of $E_3$. 
For this we need to check the assumptions of (1.10).
The former half of the assumptions can be checked by (2.4)
(1) and (2).
Assumption (1) is clear. Assumption (2) follows from assumption (6) of (2.4).
Assumption (3) follows from (2) of Theorem (2.4).
Let $l'$ be the strict transform of $l$ on $E_3$.
Hence if $(l')^2 \geq 0$, then $l^2 \geq 0$. But $l$ is a component of 
$C_1^+$ and hence $l^2 < 0$, a contradiction. 
So $(l')^2 < 0$. Since $l' = S_3 \cap E_3$, this shows that $S_3.l' < 0$.
Note that $S_3.C_2^+ >0$. So for the next extremal ray $R_3$,
we have $S_3.R_3 < 0$. Hence it is a flipping ray and 
its support is contained in $S_3$.
Note that $l'=S_3 \cap E_3$ is irreducible 
and it is the unique projective curve in $S_3$.
Hence $R_3$ is of type 3 and we are done (see (2.5)).

\definition {Subcase 2 $S_2 \cap E_2$ is irreducible}
\enddefinition

By the irreducibility of $S_2\cap E_2$, we may assume as in case (a)
that $C_2 \cap S_2 =\phi$. By the numerical property (2) in (2.5), we have
$B_2.C_2 >0$. 

\proclaim{Claim 3}
$\lda{Z_2}{E_2}$ is PLT.
\endproclaim
\demo{Proof}
If $B_2.C_1^+ \leq 0$, then $C_1^+ \subset \supp B_2$ and hence
we can prove the assertion as in Claim 1.
If $B_2.C_1^+ >0$, $B_2|_{E_2}$ is ample since $B_2.C_2 >0$.
We can prove similar statements to Claim 1 and Claim 2 in (5.2)
and deduce from these that $\ldb{Z_2}{S_2}{E_2}$ is LT.
Hence we obtain the assertion.
\qed
\enddemo
  
$$\text{Note that some component of} \ C_1^+ \ \text {is contained in} \
\supp B_2 \ \text{by} L' \subset \supp B_2, \tag *$$
where $L'$ is the strict transform of $L$ on $E_2$.
Furthermore
by the $g_2$-ampleness of $S_2$, any component of $C_1^+$ passes through $P$.

Assume first that no component of $C_1^+$ is contained 
in $\llc(\ldc{Z_2}{S_2}{B_2}{E_2}|_{E_2})$. 
By (1.6),
$C_1^+ \cap \llc(\ldc{Z_2}{S_2}{B_2}{E_2}) = P$.
Hence if $L_1$ exists, its strict transform $L'_1$ cannot passes through $P$
by $*$ and (1.1). But $L'_1$ intersects $C_1^+$, 
a contradiction to (1.6). Hence $L_1$ does not exist.
Consequently $\ldc{Z_2}{S_2}{B_2}{E_2}|_{E_2}$ is PLT near $C_2$ 
by (1.7) and $C_2$ is of type 3.

Assume next that some component of $C_1^+$ is contained 
in $\llc(\ldc{Z_2}{S_2}{B_2}{E_2}|_{E_2})$.
We know that such a component is unique. Denote it by $C_1^{'+}$.
Note that
$\llc (\ldc{Z_2}{S_2}{E_2}{B_2}|_{E_2}) = 
S_2 | _{E_2} \cup C_1^{'+} \ \text{or} S_2 | _{E_2} \cup C_1^{'+} \cup L'_1$,
where $L'_1$ is the strict transform of $L_1$.
So if $C_2 \cap \clc{0}(\ldc{Z_2}{S_2}{B_2}{E_2}|_{E_2}) = \phi$,
$C_2$ is of type 3 or 4 as in case (a) Subcase 2.
Assume below that 
$C_2 \cap \clc{0}(\ldc{Z_2}{S_2}{B_2}{E_2}|_{E_2})\not = \phi$.
Note that outside $S_2 |_{E_2}$, elements of
$\clc{0}(\ldc{Z_2}{S_2}{B_2}{E_2}|_{E_2})$ (except $Q_1$ if it exists) 
is only on $C_1^{'+}$ and
in fact only one point is on $C_1^{'+}$ outside $S_2 |_{E_2}$.
We denote this point by $e$ (clearly $e \in \supp B_2|_{E_2}$). 
Note that if $L_1$ exists, $L_1 \subset C_2$ and $e \in L_1$.
Hence in any case $e\in C_2$.
Let $C'_2$ be the connected component of $C_2$ containing $e$.
We will show that $\supp B_2 |_{E_2} - C_1^{'+} = C'_2$.
By  $\supp B_2 |_{E_2} - C_1^{'+} \cap S_2 = \phi$,  
$\supp B_2 |_{E_2} - C_1^{'+}\subset C_2$.
After contracting $\supp B_2 |_{E_2} - C_1^{'+}$ and $C_1^+ - C_1^{'+}$, 
$C_1^+$ becomes ample by (1.8)
and ampleness of $\supp B_1 |_{E_1}$ on $E_1$.
Hence we can apply (1.9) by setting
$B_1=S_2 \cap E_2$, $B_2=C_1^{'+}$, $C_1=C'_2$ and $C_2=\phi$
after the contraction of $\supp B_2 |_{E_2} - C_1^{'+}$. 
We obtain that $C'_2$ cannot pass through $\clc{0}(\ldb{Z_2}{E_2}{B_2})$
after the contraction of $\supp B_2 |_{E_2} - C_1^{'+}$
but it means $C'_2$ becomes empty 
after the contraction of $\supp B_2 |_{E_2} - C_1^{'+}$,
i.e., $C'_2=\supp B_2 |_{E_2} - C_1^{'+}$. 

If $L'_1$ exists, then $L'_1 = \supp B_2 |_{E_2} - C_1^{'+}$ and
$C'_2$ is of type (b).
(In this case $C_1^{'+}$ corresponds to $L$ in (3.2).)
Furthermore if $Q_1$ exists, $C'_2$ is of type (b1, $\lambda ' < \lambda$).
Hence the induction works well.

If $L'_1$ does not exist, $C'_2$ is of type $(c)$.
(In this case $C_1^+$ corresponds to $L$ in (3.2).)
See FIGURE (XIV).
\head 5. Preliminaries for case (b2) and (c)
\endhead

\definition{Definition (5.0) (Invariant $\delta$)}
Consider case (b2) or (c).
We define $$\multline
\delta := \sharp \{E| E \ \text{is a log crepant divisor for} \\
\ldb{X}{S}{B} \ \text{and} \ d(E, S) \leq 1 \}.
\endmultline$$
We can prove $\delta < 	\infty$. See [FA, 4.12.1 Lemma].
\enddefinition

The following lemmas is used for checking the goodness or semi-goodness
of an extraction.

\proclaim{Lemma (5.1) (LT Lemma)}
Let $(Z, \ldc{Z}{S}{B}{E})$ be a $\Bbb Q$-factorial
LC $3$-fold,
where $S$ and $E$ is irreducible surfaces.  
Assume that $S \cap E$ is irreducible 
(resp. a union of two irreducible curves) and set
$D_1 =S \cap E$ (resp. $D_1 \cup D_2 =S \cap E$).  
\roster
\item $D_1\simeq \Bbb P^1$. $S$ and $E$ are generically normal crossing
along $D_1$
(resp. $D_1\simeq \Bbb P^1$. $S$ and $E$ are generically normal crossing
along $D_1$ and $D_2$);
\item $\clc{0}(\ldc{Z}{S}{B}{E}|_{E}) \cap D_1$ consists of one point $P$
and $P \in \supp B|_E$ 
(resp. $\clc{0}(\ldc{Z}{S}{B}{E}|_{E}) \cap D_1$ consists of one point $P$
and $P=D_1\cap D_2$);
\item $S$ is normal.
\item
$M:=\text{LLC} \ (K_Z+\-{g}S+\-{g}B+E|_{E^{\mu}})$ is connected and 
we can write
$$M=\mu^{-1}D_1, B_1,..., B_m$$ (they form a chain with this order) (resp. 
$$M=\mu^{-1}D_1,\mu^{-1}D_2, B_1,....,B_m$$
(they form a chain with this order)),
where $\mu$ is the normalization of $E$ and $B_i$ is an irreducible curve
contained in $\supp B|_E$.
\endroster

Then $\ldb{Z}{S}{E}$ is LT (resp. $K_Z+S+E$ is LT outside $P$).
\endproclaim 

\demo{Proof}
By assumption (4), $E$ is normal.
Also by assumption (4), $\text{CLC}_0(K_Z+S+E|_E)$ is contained in $D_1$
(resp. $D_1$ and $D_2$)
and does not contain 
$D_1\cap \text{Supp} \ B|_E$ 
(resp. $D_2\cap \text{Supp} \ B$).
Hence 
$K_Z+S+E|_E$ is LT
(resp. $K_Z+S+E|_E$ is LT outside $P$), 
which in turn show that $K_Z+S+E$ is LT
(resp. $K_Z+S+E$ is LT outside $P$)
by the inversion of adjunction ([S, 5.13 proposition]).
\qed
\enddemo

\proclaim{Lemma (5.2) (Lemma for goodness)}
Consider the situation of $(b2)$ or $(c)$ in (3.0).
Let $g:Z\to X$ be a (primitive) log crepant divisorial extraction
for $\ldb{X}{S}{B}$.
Let $E$ be the exceptional divisor and $\cup_{i=1}^n D_i:=\-{g}S\cap E$,
where $D_i$'s form a chain and $D_1$ intersects with $\-{g}C$.
Then $n=1\ \text{or}\ 2$.

Furthermore we assume that the following hold:

If $n=1$ (resp. $n=2$),
\roster
\item $D_1\simeq \Bbb P^1$. $\-{g}S$ and $E$ are generically normal crossing
along $D_1$
(resp. $D_1\simeq \Bbb P^1$. $\-{g}S$ and $E$ are generically normal crossing
along $D_1$ and $D_2$);
\item $\ldc{Z}{\-{g}S}{\-{g}B}{E}|_{E}|_{D_1}=
\ldc{D_1}{\frac12 P_1}{\frac12 P_2}{P}$ 
(resp. we assume the same thing and also assume that $P=D_1\cap D_2$);
\item $\-{g}S$ is normal.
\endroster
Then $g$ is a good contraction (resp. a semi-good contraction).
\endproclaim

\demo{Proof}
Note that $\-{g}B$ is $g$-ample.
So if $n\geq 3$, $D_2,.., D_{n-1}$ are contained in fibers but
do not intersect $\-{g}B$, a contradiction. 
Hence $n\leq 2$.

Assume that $n=1$ (resp. $n=2$) and above $(1)$, $(2)$ and (3) hold.
We have only to check the $(2)$ of the definition of good extraction
(resp. semi-good extraction) in (3.3). For this, we use Lemma (5.1).
We have to check assumption (4) of Lemma (5.1).

Let $\mu:E^{\mu}\to E$ be the normalization of $E$.

\proclaim{Claim 1}
$\text{LLC} \ (K_Z+\-{g}S+\-{g}B+E|_{E^{\mu}})$ is connected.
\endproclaim
\demo{Proof}
In the case $g(E)\not=L$, $E$ is a projective surface such that
$\ldc{Z}{\-{g}S}{\-{g}B}{E}|_{E}\equiv 0$.
So by (1.7) and the assumption $(2)$,
we immediately obtain
that $\text{LLC} \ (K_Z+\-{g}S+\-{g}B+E|_{E^{\mu}})$ is connected.

Next we treat the case $g(E)=L$.
By the assumption $(1)$, $n=2$ in this case.

Assume that $\text{LLC} \ (K_Z+\-{g}S+\-{g}B+E|_{E^{\mu}})$ is not
connected.
Then there is a point over $Q$ contained in
$\text{LLC} \ (K_Z+\-{g}S+\-{g}B+E|_{E^{\mu}})$ and not in $D_1$.
In particular, the fiber is reducible.
So we can contract bimeromorphically all the components of the fiber over $Q$
which are not contained in
$\text{LLC} \ (K_Z+\-{g}S+\-{g}B+E|_{E^{\mu}})$.
But this contradicts (1.6).
Hence we finish the proof of Claim 1.
\qed
\enddemo

Let $M:= \ \text{LLC} \ (K_X+\-{g}S+\-{g}B+E|_{E^{\mu}})$.
By Claim 1 and $(2)$, 
this is a tree of curves and $\mu^{-1}D_1$ is an end of
the tree.

\proclaim{Claim 2}
$$
M=
\cases
\mu^{-1}D_1, B_1,.., B_m \ 
\text{(with this order) in case} \ n=1 \\
\mu^{-1}D_1,\mu^{-1}D_2, B_1,....,B_m \ \text{(with this order)
in case} \ n=2 \ \text{and} \ g(E)\not =L \\
\mu^{-1}D_1,\mu^{-1}D_2 \ \text{(with this order)
in case} n=2 \ \text{and} \  g(E)=L.
\endcases
$$
In the above, $B_i$ is an irreducible curve
contained in $\supp \-{g}B|_{E^{\mu}}$.
\endproclaim

\demo{Proof}
Note that in case (b2), $g(E) = Q \ \text{or} \ C$ and
in case (c), $g(E) = Q \ \text{or} \ L$.
First assume that $g(E)\not = L$.
If $n=1$ (resp. $n=2$), $\-{g}B|_{E^{\mu}}$ is ample
(resp. ample in the cases $g(E)=Q$ or
nef and only trivial for $D_1$ in the case $g(E)=C$).
(For the case $g(E)=C$, see Theorem (4.0) Claim 2 in detail.) 
Hence $\supp {\-{g}B|_{E^{\mu}}}$ is connected.
(In the ample case, this is well known. In the case $g(E)=C$,
$D_1$ coincides with the support of the next flipping ray
(see Theorem (4.0) Claim 2).
Hence after contracting $\-{\mu}D_1$ in $E^{\mu}$, 
$\-{g}B|_{E^{\mu}}$ becomes ample. Hence we obtain the connectedness
as usual.)

Let $\mu^{-1}D_1, B_1,..., B_n$
(resp. $\mu^{-1}D_1, C_1,..., C_m, \mu^{-1}D_2, B_1,....,B_n$)
be the irreducible components of $M$ (they form a tree with this order). 
We will show that 
$B_1,.., B_n$ is exactly in $\text{Supp} \ \-{g}B|_{E^{\mu}}$
\newline 
(resp. $B_1,.., B_n$ is exactly in $\text{Supp} \ \-{g}B|_{E^{\mu}}$
and $\{ C_i \}$ is an empty set).
Let $\text{Supp} \ \-{g}B|_{E^{\mu}}=B'+B''$ be the decomposition such that
$B'$ is the union of components of $\text{Supp} \ \-{g}B|_{E^{\mu}}$ contained
in $M$.
Let $B'_k$ be any connected component of $B'$.
By the connectedness of $\supp {\-{g}B|_{E^{\mu}}}$,
$B'_k\cap B''\not=\phi$ so an end of $B'_k$ 
is $B_n$ since $\deg (K_X+\-{g}S+\-{g}B+E|_{E^{\mu}}|_{B_j})=0$
for $1 \leq j \leq n$ 
and for $j\not =n$ $B_j$ intersects already two component of 
$\llc (K_Z+\-{g}S+\-{g}B+E|_{E^{\mu}})$.
In particular we know that $B'$ is actually connected  
and $B'\subset B_1\cup ...\cup B_n$.
But by the fact $\-{g}B|_{E^{\mu}}$ is ample
(resp. ample in the cases $g(E)=Q$ or
nef and only trivial for $D_1$ in the case $g(E)=C$), 
$\mu^{-1}D_1$
(resp. $C_i$ and $\mu^{-1}D_2$) 
must intersect $B'$ so we are done.

Next we treat the case $g(E)=L$.
Let $\mu^{-1}D_1, C_1,..., C_m, \mu^{-1}D_2$
be the irreducible components of $M$, 
where $C_i$ is a component of a fiber of $g$.
No component of $\supp \-{g}B|_{E^{\mu}}$ pass through $P$.
In fact otherwise on $\-{g}S$ (this is normal by the assumption $(3)$), 
$K_Z+\-{g}S+\-{g}B+E|_{\-{g}S}$ has two reduced boundary $D_1$ and $D_2$
and another boundary, a contradiction.
Hence $\{ C_i \}$ must be empty 
since $\-{g}B|_{E^{\mu}}$ is $g \circ \mu$-ample.
Now we finish the proof of Claim 2
\qed
\enddemo

By these 2 claims, we can apply (5.1) and we are done.
\qed
\enddemo

\definition{Remark (5.3)}
We can prove the normality of $E$ even if we replace (1) 
by the following weaker assumption in case $n=2$:

(1)'  
$D_1\simeq \Bbb P^1$. $\-{g}S$ and $E$ are generically normal crossing
along $D_1$.

In fact the assumption that $\-{g}S$ and $E$ are generically normal crossing
along $D_2$ is used only in the last paragraph of the above proof. 
\enddefinition

\head 6. Treatment of case (b2, $\delta$) 
\endhead
We use the notation as in (3.2).
\proclaim{Theorem (6.0)}
There is a good extraction or semi-good extraction $g:Z\to X$ such 
that $g(E) =C$, where $E$ is the exceptional divisor of $g$. 
Furthermore if $g$ is good one of the following holds:
 
\roster
\item $\supp B$ and $S$ are not generically simply tangent along $C$ and
$d(E, S)\leq 1$
or 
\item $\supp B$ and $S$ are generically simply tangent along $C$
and $\pb{g}{X}{S}{B}|_E$ is PLT over a general point of $C$.
\endroster
If $g$ is semi-good, let $D_1 \cup D_2 := \-{g}S \cap E$,
where $D_1$ is the strict transform of $C$ on $\-{g}S$.
Then $D_1$ coincides the support of a flipping ray
satisfying the assumption of (2.4). 
\endproclaim
\demo{Proof}
\proclaim{Claim 1}
Except the case that 
$\supp B$ and $S$ are generically simply tangent along $C$,
there is a log crepant exceptional divisor for $K_X+S+B$
with multiplicity in $S$ $\leq 1$.
\endproclaim
\demo{Proof}
It is sufficient to prove the statement after cutting $X$ by
a general hyperplane section at a general point of $C$.
Let $x$ be a general point of $C$ and $H$ a general hyperplane section
of $X$ through $x$.
Let $\mu:H_t\to H$ be the MLT for $(H, \ldb{X}{S}{B}|_{H}, x)$ and 
$F_i \ (1\leq i \leq n)$ exceptional curves, where $F_1$ intersects 
$\-{\mu}(S|_H)$.
First assume that
there are $2$ reduced components of the boundary of $\ldb{X}{S}{B}|_{H}$.
If $x$ is smooth, then the simple blow up at $x$ gives what we want.
If $x$ is singular, then MLT coincides with MRS, so by (1.3),
we are done.
Next assume that
there is $1$ reduced component of the boundary of $\ldb{X}{S}{B}|_{H}$.
If $x$ is singular, then MRS is dominated by MLT, so we are done by (1.3).
If $x$ is smooth and $n\geq 2$, $d(F_1, S|_H)=1$. Hence we are done.
If $x$ is smooth and $n=1$, then by the description of MLT, 
we encounter the following
$3$ possibilities: See FIGURE (I)

Only in the middle case, we have no log crepant curve
as we want.
Hence we obtain Claim 1.
\qed
\enddemo

We go back to $3$-dimensional situation.
Let $h_t:X^t\to X$ be a $\Bbb Q$-factorial LT model
for $\ldb{X}{S}{B}$ near $S$ and we denote the exceptional divisor
corresponding to $E_1$  by $E$.
Run $N:=h_t^*(\ldb{X}{S}{B})-\epsilon E$-MMP and
let $g':Z'\to X$ be the end result.
Then $N':=K_{Z'}+S'+B'+E'-\epsilon E-\delta {g'}^*S$
($0<\delta << 1$) is KLT and $g'$-nef and $g'$-big, 
its sufficient multiple is $g'$-free by
Kawamata's base point free theorem.  Let $h:Z' \to Z$ be the morphism 
defined by such a multiple and $g:Z\to X$ be the natural morphism. 
Then exceptional divisor of $g$ is
exactly $E$. 
Indeed,
for a general curve in $E$, $N$ is nef
while MMP
since $N\equiv -\epsilon E$ over $X$.
So $E$ is $g$-exceptional.
On the other hand some $h_t$-exceptional divisor $F$ 
($\not =E$) are $g$-exceptional,
a general curve in $F$, 
$h(N')$ is numerically nonpositive, a contradiction. 

We will prove that this $g$ satisfies the statement of Theorem.
By (5.2), we can write
$\-{g}S\cap E=D_1 \ \text{or} \ D_1\cup D_2$.
By the choice of $E$, $\-{g}S$ and $E$ are generically normal crossing 
along $D_1$.
By the assumption of (b2), $(2)$ of Lemma (5.2) holds.
Note that $\-{g}S$ is normal since $S$ is normal and $\-{g}S$ is normal
along intersection curves with $g$-exceptional divisors.
Hence by (5.2), $g$ is good if $\-{g}S\cap E$ is
irreducible. In this case we check that (1) or (2) holds.
If $\supp B$ and $S$ are not generically simply tangent along $C$,
(1) holds by Claim 1.
If $\supp B$ and $S$ are generically simply tangent along $C$,
the restriction of $g$ over a general hyperplane section of $X$ at a general
point of $C$ coincides with MLT. Hence (2) of Theorem (6.0) holds.

Consider the case $\-{g}S\cap E$ is reducible in the below. 
Run the BTM starting from $g$. We use the notation as in (2.5).

\proclaim{Claim 2}
$\-{g}B.R_1=0$, $E.R_1>0$ and $\-{g}S.R_1<0$.
$R_1$ is a flipping ray with $\supp R_1=D_1$
which satisfies the assumption of (2.4).
\endproclaim

\demo{Proof}
It is easy to see that $\-{g}B.D_1=0$.
On the other hand, $\-{g}B$ is positive for fibers of $g$.
Hence $\-{g}B.R_1\leq 0$ since we play a $2$-ray game.
Considering $g^*B$ and $g^*S$,
we obtain that $E.R_1>0$ and $\-{g}S.R_1<0$.
By this $\supp R_1 \subset \-{g}S$.
Clearly $\supp R_1 \subset E$.
By these, we see that $\supp R_1 = D_1$.
The check that it satisfies the assumption of (2.4) is easy
so we left it to readers.
\enddemo

By the property (2) of (2.4), (1) of Lemma (5.2) holds.
Hence $g$ is semi-good and we are done.

\qed
\enddemo

\definition{(6.1) Conclusion of case (b2)}
We fix a divisorial extraction $g: Z\to X$ as in (6.0)
and apply BTM starting by this $g$. We will use the notation as in (2.5).
Our conclusion is explained in the following flow chart:
\enddefinition
\newpage
\hfil Case (b2) \hfil 
$$
\fbox{2.5in}{
Choose a primitive log crepant extraction $g$ as in (6.0) 
and run BTM 
starting by $g$}
\hskip1in
$$
\vskip0.3in
$$\hskip0.5in\oversetbrace\to{
\hbox{\hsize=3.4in
\fbox{1.1in}{(A)
\break
$g$ is semi-good. \break
Then $R_1$ is of \break
the type \break
as in (2.4).}
\hskip1.0in
\fbox{0.8in}{(B)\break
$g$ is good.}
}
}
\hskip2.5in$$
\vskip0.3in
$$\downarrow\hskip3.0in$$

\vskip0.1pt
$$
\fbox{1.1in}{$R_2$ is a flipping \break
ray of type 2. \break
DONE}\hskip5pt
\oversetbrace\to{
\hbox{\hsize=4.4in
\fbox{1.2in}{$R_1$ is a flipping 
\break
ray and \break
$\pb{g}{X}{S}{B}|_{g^*H}$ \break
is PLT}
\hskip1pt
\fbox{1.1in}{$R_1$ is a divisorial \break 
ray. \break 
DONE}
\hskip1pt
\fbox{1.1in}{$R_1$ is a flipping 
\break
ray and $\pb{g}{X}{S}{B}|_{g^*H}$ is not PLT}
}}$$
 
\vskip0.5in

$$
\oversetbrace\to{
\hbox{\hsize=2.4in
\fbox{1.2in}{$\text{LLC} = S_1 \cup E_1$.
Then $C_1$ is of type 3. \break
DONE}
\hskip.1in
\fbox{1.2in}{$\text{LLC}=S_1 \cup E_1 \cup l_1$ }}}
\oversetbrace\to{
\hbox{\hsize=2.4in
\fbox{1.0in}{$C_1 \cap \clc{0}= \phi$.
Then $C_2$ is of type 3 or 4. \break
DONE}
\hskip.1in
\fbox{1.1in}{
$C_2 \cap \clc{0} \not = \phi$.
Then the connected component $C'_1$ of $C_1$ containing\break
$\clc{0}$ is of \break
type (b2, $\delta ' < \delta$). 
(Another \break
component \break
is of type 3 or 4. \break
So we are done)
}}}
$$
\vskip0.5in

$$
\oversetbrace\to{
\hbox{\hsize=2.4in
\fbox{1.0in}{$C_1 \cap \clc{0}= \phi$. Then $C_1$ is of type 3 or 4.\break
DONE}
\hskip.1in
\fbox{1.1in}{$C_1 \cap \clc{0}\not = \phi$.
Then the connected \break
component 
$C'_1$ of \break
$C_1$ containing \break
$\clc{0}$ is of type (c).}}}
$$
\vskip1.2in

$$\hskip1.8in\Downarrow \hskip1.5in \Downarrow$$
\vskip0.2in
$$\hskip2.2in
\fbox{1.1in}{Reduction to \break
case (c)}\hskip0.3in
\fbox{1.0in}{Reduction to \break
case (b2, $\delta ' < \delta$)}
$$
\vskip0.5in

$$\clc{0}:=\clc{0}(\ldc{Z_1}{S_1}{B_1}{E_1}|_{E_1})$$
$$\llc:=\llc (\ldc{Z_1}{S_1}{B_1}{E_1}|_{E_1})$$
\newpage

\definition{(A)}
First we assume that $g$ is semi-good. See FIGURE (XVI).
As in the proof of (6.0), 
$R_1$ is a flipping ray satisfying the assumption of (2.4)
and $\supp R_1=D_1$.
After the flip $D_1$, we have $S_2.{C_1}^+ >0$ and $E_2.{C_1}^+ <0$.
Then we obtain the following:
\enddefinition

\proclaim{Claim 1}
$S_2.R_2<0$.
\endproclaim

\demo{Proof}
We derive a contradiction assuming that $S_2.R_2\geq 0$.
Let $\nu:E_2^{\nu}\to E_2$ be the normalization.
First we see that $({D_2}^+)^2>0$ on $E_2^{\nu}$ as follow:

Since ${D_2}^+=S_2\cap E_2$, and $S_2$ and $E_2$ are generically
normal crossing along ${D_2}^+$, $({D_2}^+)^2=S_2.{D_2}^+$.
On $Z_2$, two extremal rays over $Y$ are $R_2$ and the ray generated by
$D_1^+$.
They satisfy $S_2.R_2\geq 0$ and $S_2.D_1^+>0$.
If ${D_2}^+ \subset \supp R_2$, ${D_2}^+$ is contractible in $E_2$.
But by ${D_2}^+=S_2\cap E_2$, we have $S_2.{D_2}^+<0$, a contradiction.
Hence ${D_2}^+\not \subset R_2$, which in turn show that $S_2.{D_2}^+>0$.

Next we deduce from $({D_2}^+)^2>0$ on $E_2^{\nu}$ that
$(D_2)^2>0$ on $E_1$. (This give us a contradiction
since $D_2$ is contained in a fiber of $g$.)
Let $Q$ be the point as in Theorem (2.4) (3).
Note that $Q$ is a smooth point.
\comment
Hence considering the degree of $\ldc{Z_2}{S_2}{B_2}{E_2}|_{D^*}$,
we know that on $D^*$
there is at most one non LT point for 
$\ldc{Z_2}{S_2}{B_2}{E_2}|_{E_2^{\nu}}$.
\endcomment
On the other hand on $D_1$ there is one non PLT point for 
$\ldc{Z}{\-{g}S}{\-{g}B}{E}|_{E}$ by the assumption of $(b2)$.
By these and (1.10),
\comment
$E$ is dominated by MLT along $D^*$ for
$\ldc{Z_2}{S_2}{B_2}{E_2}|_{E_2^{\nu}}$.
When we go from $E^{\nu}$ to its MLT, there is no exceptional curves
over $Q$. 
Hence 
\endcomment
we have $(D_2)^2>0$ on $E_1$.
As we saw above, this give a contradiction and prove Claim 3.
\qed
\enddemo

Hence it is a flipping ray and 
its support is contained in $S_2$.
Note that $S_2 \cap E_2$ is irreducible 
and it is the unique projective curve in $S_2$.
Hence $R_2$ is of type $2$ and we are done (see (2.5)).
\definition{(B)}
Next we assume that $g$ is good.
We may assume that $R_1$ is a flipping ray.
Since $S_1\cap E_1$ is irreducible, we may assume that
$C_2 \cap S_2 = \phi$ (see Conclusion of (a) and (b1)).
Let $x$ be a general point of $C$ and $H$ a general hyperplane section
of $X$ through $x$.
Let $F_1,\dots F_n$ be the exceptional curves of MLT for
$(H, \ldb{X}{S}{B}|_{H}, x)$.
We consider by dividing into two cases.
\enddefinition

\definition {Subcase 1. $\ldc{Z}{S_1}{B_1}{E_1}|_{g^*H}$ is PLT}
See FIGURE (XVII).
\enddefinition
In this case $g_1$ coincides with MLT over $H$. 
Hence
$$\llc(\ldc{Z}{S_1}{B_1}{E_1}|_{E_1})= (S_1\cap E_1) \ \text{or} \ 
l_1\cup \dots  l_m \cup (S_1\cap E_1),$$
where $l_i$ is a irreducible curve contained in a fiber of $g$.
In the former case $R_1$ is a flipping ray of type 3 so we are done by (2.2).
We assume that the latter case occurs.
In this case $m=1$ holds.
For otherwise some $l_i$ would not intersect with $S_1\cap E_1$
(consider the degree of the boundary of
$\ldc{Z}{S_1}{B_1}{E_1}|_{E_1}|_{S_1 \cap E_1}$).
But this contradicts the $g$-ampleness of $S_1$.
By considering the degree of the boundary of 
$\ldc{Z}{S_1}{B_1}{E_1}|_{E_1}|_{l_1}$,
we can prove that there is at most one point in
$\clc{0}(\ldc{Z}{S_1}{B_1}{E_1}|_{E_1})$
outside $S_1\cap E_1$. 

If $C_1$ does not pass through a point in
$\clc{0}(\ldc{Z}{S_1}{B_1}{E_1}|_{E_1})$
outside $S_1\cap E_1$, $R_1$ is a flipping ray of type 3 or 4 so we are done
by (2.2) and (2.3).
We assume that $C_1$ passes through the point $Q_1$ in
$\clc{0}(\ldc{Z}{S_1}{B_1}{E_1}|_{E_1})$
outside $S_1\cap E_1$.
Note that 
$(\supp B_1|_{E_1}-l_1).S_1=0$.
Hence by $S_1.l_1>0$ and $S_1.C_1=0$ and $\rho(Z_1/Y)=2$,
we have $\supp B_1|_{E_1}-l_1 \subset C_1$.
Let $C'_1$ be the connected component of $C_1$ containing
$\supp B_1|_{E_1} -l_1$.
By $Q_1\in C'_1$ and (1.9), we can prove as before
$C'_1=\supp B_1|_{E_1}-l_1$.
Since $\supp B_1|_{E_1}-l_1$ is not contained in
$\llc(\ldc{Z}{S_1}{B_1}{E_1}|_{E_1})$,
the type of $C'_1$ is (c) (in this case $l_1$ corresponds to $L$).

\definition {Subcase 2. Contrary to Subcase 1}
See FIGURE (XVIII).
\enddefinition
In this case, (1) holds in (6.0).

\proclaim{Claim 2}
$\supp B_1|_{E_1}$ is contained in 
$\llc(\ldc{Z}{S_1}{B_1}{E_1}|_{E_1})$ and is a section or 
a union of a section and an irreducible fiber.
\endproclaim
\demo{Proof}
Recall that in the construction of $g$, we take $E=E_1$.
Hence the fiber of $g$ over $x$ intersects $\supp B_1|_{E_1}$ at one 
point and the intersection point is contained in 
$\llc(\ldc{Z}{S_1}{B_1}{E_1}|_{E_1})$
(cf. (1.4)). So $\supp B_1|_{E_1}$ is a union of
a section contained in $\llc(\ldc{Z}{S_1}{B_1}{E_1}|_{E_1})$ and 
(possibly empty) components of fibers. Assume that the fiber-part is not empty.
Then it must pass through $P$ by (5) in Set up (3.0) 
and $g_1$-ampleness of $S_1$ 
and hence the section-part
cannot pass through $P$ by (1.1).
By (1.7) and the property (3) of a good extraction in (3.3),
$\llc(\ldc{Z}{S_1}{B_1}{E_1}|_{E_1})$ is connected, so some component of the
fiber-part is contained in $\llc(\ldc{Z}{S_1}{B_1}{E_1}|_{E_1})$.
By this and (1.1), the fiber-part must be irreducible
and contained in $\llc(\ldc{Z}{S_1}{B_1}{E_1}|_{E_1})$. 
\qed
\enddemo

We can prove as in Subcase 1 that there is at most one point in
$\clc{0}(\ldc{Z}{S_1}{B_1}{E_1}|_{E_1})$
outside $S_1\cap E_1$. 

If $C_1$ does not pass through a point in
$\clc{0}(\ldc{Z}{S_1}{B_1}{E_1}|_{E_1})$
outside $S_1\cap E_1$, $R_1$ is a flipping ray of type 3 or 4 so we are done
by (2.2) and (2.3).

We assume that $C_1$ passes through the point $Q_1$ in
$\clc{0}(\ldc{Z}{S_1}{B_1}{E_1}|_{E_1})$
outside $S_1\cap E_1$. 
In this case we will prove that the former case in Claim 2 cannot occur.
Indeed otherwise $Q_1$ would sit on a section contained in 
$\llc(\ldc{Z}{S_1}{B_1}{E_1}|_{E_1})$. Let $m$ be the fiber 
containing $Q_1$ and $m'$ its connected component containing $Q_1$
and intersecting $S_1$ at a time. After extracting a log crepant curve
from $Q_1$, $m'$ become contractible but 
this contradicts (1.6).

Hence we denote the fiber-part of $\llc(\ldc{Z}{S_1}{B_1}{E_1}|_{E_1})$ by $l$.
Let $C'_1$ be the connected component of $C_1$ through $Q_1$.
By the similar argument to Subcase 1 using (1.9),
we can prove $C'_1=\supp B_1|_{E_1}-l$.
Furthermore $C'_1$ is of type (b2, $\delta ' < \delta$).

\head 7. Treatment of case $(c)$
\endhead

\definition{(7.0) Set up for (c)}
If $\delta=0$, we set 
$$d:=\min \{d(E,S)| E\ \text{is a log crepant exceptional divisor for} 
\ \ldb{X}{S}{B}\}.$$
Let $$I:= \{i| E_i \ \text{is a log crepant exceptional divisor for} 
\ \ldb{X}{S}{B} \ \text{with} \ d_i\leq 1\}$$
if $\delta>0$ or
$$I:= \{i| E_i \ \text{is a} \ 
\text{is a log crepant exceptional divisor for} 
\ \ldb{X}{S}{B} \ \text{with} \ d_i=d\}$$ if $\delta=0$, where
$d_i:=d(E_i, S)$.
Let $h: Z_t\to X$ be a $\Bbb Q$-factorial LT model of 
$K_X+S+B$ near $S$ which extracts any divisor $E_i$ with $i\in I$. 
We apply 
$N:= K_{Z_t}+S_t+B_t+E_t-\epsilon (\sum\limits_{i\in I} d_iE_i)$-MMP 
over $X$,
where $S_t$ (resp. $B_t$) is the strict transform of $S$ (resp. $B$)
and $E_t$ is the reduced exceptional divisor of $h$.
Since any flip while MMP is a flip of type 1, we can carry out
the program.
Let $g': Z'\to X$ be the end result. 
Since 
$N':=K_{Z'}+S'+B'+E'-\epsilon (\sum\limits_{i\in I} d_iE_i)-\delta {g'}^*S$
($\delta << 0$) is KLT and $g'$-nef and $g'$-big, 
its sufficient multiple is $g'$-free by
Kawamata's base point free theorem.  Let $h:Z' \to Z''$ be the morphism 
defined by such a multiple and $g'':Z''\to X$ be the natural morphism. 
Then exceptional divisors of $g''$ are
exactly the divisors $E_i$'s such that $i\in I$. 
Indeed,
for any $E_i$ with $i\in I$ and a general curve in $E_i$, $N$ is positive
while MMP
since $N\equiv \epsilon(S_t+\sum\limits_{j\not\in I} d_j E_j)$ over $X$.
So any $E_i$ with $i\in I$ is $g''$-exceptional.
On the other hand some $h_t$-exceptional divisors $E_j$'s 
with $j\not\in I$ are $g''$-exceptional,
for at least one $E_j$ and a general curve in $E_j$, 
$h'(N')$ is numerically nonpositive, a contradiction. 
Note that $\-{g''}S$ is normal since $S$ is normal and $\-{g''}S$ is normal
along intersection curves with $g''$-exceptional divisors.
Let the chain $\cup_{i\in I}E_i\cap {g''}^{-1}S$ be 
$\cup_{k=1}^n D_k$, where $D_1$ intersects ${g''}^{-1}C$
and $D_n$ intersects ${g''}^{-1}L$ (if $L$ is the center of some $E_i$
with $i\in I$, $D_n$ is equal to ${g''}^{-1}L$).
Then we encounter the following two possibilities:
\roster
\item"(1)" $\-{g''}C \cap D_1$ is contained in 
$\clc{0}({g''}^*(K_X+S+B)|_{{g''}^{-1}S})$;
\item"(2)" $\-{g''}C \cap D_1$ is not contained in 
$\clc{0}({g''}^*(K_X+S+B)|_{{g''}^{-1}S})$.
\endroster

We can express above two possibilities in a more concrete way:
\roster
\item"(1)"
\item"(1.1)"
$D_1\not= {g''}^{-1}L$ (and hence $D_1\simeq \Bbb P^1$) and
${g''}^*(K_X+S+B)|_{{g''}^{-1}S}|_{D_1}=K_{D_1}+P'+P$,
where $P:=D_1\cap D_2$ if $n \geq 2$ and $P:= D_1 \cap \-{g''}L$ if $n=1$

or
\item"(1.2)"
$D_1= {g''}^{-1}L$ and $P':= {g''}^{-1}C\cap {g''}^{-1}L$;
\item"(2)"
$D_1\not= {g''}^{-1}L$ (and hence $D_1\simeq \Bbb P^1$) and
${g''}^*(K_X+S+B)|_{{g''}^{-1}S}|_{D_1}=K_{D_1}+\frac12 P_1+\frac12 P_2+P$,
where $P:=D_1\cap D_2$ if $n\geq 2$ and $P:= D_1 \cap \-{g''}L$ if $n=1$.
\endroster
See FIGURE (XIX).
In any case, we call $E$ the exceptional divisor which cut out $D_1$
and $d:=d(E, S)$ (it is consistent with the above definition of $d$). 
We remark that $\-{g''}S$ and $E$ cross normally along $D_i$.
In fact otherwise 
there is a log crepant divisor $F$ for $\ldb{X}{S}{B}$ whose center is $D_i$
with $d(F, \-{g''}S+E)\leq 1$.
If $\delta>0$, $1\geq d(F, \-{g''}S+E)\geq d(F, \-{g''}S+dE)$, i.e., 
$d(F,S)\leq 1$. But this contradicts the property of $h$.
If $\delta=0$,
$$d\geq d(F, d\-{g''}S+dE)\geq d(F, \-{g''}S+d(E,S)E), \ \text{i.e.,} \ 
d(F,S)\leq d.$$ 
But this again contradicts the property of $h$.
\enddefinition
We remark the following:

\proclaim{Claim}
Above division into cases (1-1), (1-2) and (2) depends only on $f$.
\endproclaim

\demo{Proof}
By running ${g''}^*(K_X +S +B)-\epsilon(\-{g''}S +dE)$-MMP and
taking the LC model as in (7.0), we obtain the primitive
extraction of $E$ from $g''$. 
We can check the above procedure is isomorphism 
near $P'$ in case (1) of (7.0) (resp. near $P_1$ and $P_2$
in case (2) of (7.0).)

Assume that the case (2) occurs in the above.
Note that even after making $g''$ a primitive extraction of $E$, 
the property (2) holds for this primitive contraction.

Assume furthermore that
there is another $\t{g''}$ with the similar properties to $g''$ but
the case (1) occurs.
Then we can obtain the primitive extraction of the strict transform of $E$
from $\t{g''}$ and after making $\t{g''}$ a primitive extraction of $E$,
the property (1) for this primitive contraction.
Hence we obtain two different primitive extraction of $E$.
But this contradicts [FA, Lemma 6.2].

If the case (1) occurs, then by above argument, (2) does not occur
for any $\t{g''}$. 

It is determined only by $f$ whether $L$ is a center of an exceptional divisor 
$E_i$ for $i \in I$.
Hence it is detremined only by $f$ whether (1-1) or (1-2) occurs.  
This completes the proof.
\qed
\enddemo

\newpage
We treat the case (1) and (2) separatedly in the below.
We first treat the case $(1)$. We deduce the following
from the condition of $(1)$.
\proclaim{Claim }
For any log crepant divisor $F$ for $K_X+S+B$ whose center
is $P'$ or a curve through $P'$, $d(F, {g''}^{-1}S+E)>1$.
\endproclaim
\demo{Proof}
Note that among $E_i$ with $i\in I$, $E$ is the unique divisor 
whose center contains $P'$.
Hence if $\delta>0$, 
$d(F, {g''}^{-1}S+E)\geq d(F, {g''}^{-1}S+dE)=d(F, S)>1$ and 
if $\delta=0$, 
$d(F, {g''}^{-1}S+E)\geq d(F, \frac1d ({g''}^{-1}S+dE))=\frac 1d d(F, S)>1$.
\qed
\enddemo
Note that $\ldb{Z''}{\-{g''}S}{\-{g''}B}$ and 
$\ldb{Z''}{E}{\-{g''}B}$ is PLT at $P'$ by
normality of $\-{g''}S$ and $E$ at $P'$,
and the inversion of adjunction (cf.
the subadjunction formula in (1) above).  
By these, we obtain the following singularity.

\definition{(*)}
(Z, P) is an algebraic germ of a normal $3$-fold singularity,
$S_1$ and $S_2$ are irreducible $\Bbb Q$-Cartier divisors  and $B$ is
a $\Bbb Q$-Cartier boundary such that $\llcorner B \lrcorner = \phi$
with the following properties:
\roster
\item"(i)" $K_Z+S_1+S_2+B$ is LC and 
$K_Z+S_i+B$ is PLT for $i = 1, 2$;
\item"(ii)" $2(K_Z+S_1+S_2+B)\sim 0$;
\item"(iii)" $C:=S_1\cap S_2$ is a smooth curve and 
$K_Z+S_1+S_2+B|_{S_i}|_C= K_C+P$, i.e., 
$K_Z+S_1+S_2+B|_{S_i}$ is not PLT at $P$ for $i = 1,2$;
\item"(iv)" 
for any log crepant divisor $F$ with respect to $K_Z+S_1+S_2+B$ whose center
is $P$ or a curve through $P$, $d(F, S_1+S_2)>1$.
\endroster
\enddefinition

We examine the $(*)$ for a while and prove
\proclaim{Theorem (7.1) (Classification of (*))}
We consider the object $(*)$.
Then one of the following holds:
\roster
\item $P$ is smooth on $Z$, $S_1$ and $S_2$;
\item locally analytically $(P\in S_1\cup S_2 \subset Z)\simeq
(o\in (xy=0) \subset \frac13(2,2,1))$.
The weighted blow up $h:Z'\to Z$ with weight $\frac13(2,2,1)$ is log crepant
for $\ldc{Z}{S_1}{S_2}{B}$. Let $F'$ be the exceptional divisor of $h$.
Then $\pc{h}{Z}{S_1}{S_2}{B}|_{F'}=\ldc{F'}{M_1}{M_2}{\frac12 M_3+\frac12M_4}$,
where all $M_i$'s are irreducible curves, $M_i:=\-{h}S_i\cap F'$
for $i=1,2$, $M_3:=\supp \-{h}B\cap F'$ and $M_4$ is an ODP curve singularity.
On $F'$, $\pc{h}{Z}{S_1}{S_2}{B}$ is LT;
See FIGURE (XX).

\item locally analytically $(P\in S_1\cup S_2 \subset Z)\simeq
(o\in (xy=0) \subset \frac12(1,0,1))$.
\endroster
\endproclaim

We begin with some lemmas.
\proclaim{Lemma (7.2)}
Assume that $Z$ is not smooth at $P$.
Let $\t{B}$ be a $\Bbb Q$-Cartier boundary such that
\roster
\item"(i)" $K_Z+S_1+S_2+\t{B}$ is LC and 
$K_Z+S_i+\t{B}$ is PLT for $i = 1, 2$;
\item"(ii)" $2(K_Z+S_1+S_2+\t{B})\sim 0$.
\endroster
Then 
$K_Z+S_1+S_2+\t{B}|_{S_i}|_C= K_C+P$, i.e., 
$K_Z+S_1+S_2+\t{B}|_{S_i}$ is not PLT at $P$ for $i= 1, 2$.
\endproclaim
\demo{Proof}
Assume the contrary to what we want, i.e., 
$K_Z+S_1+S_2+\t{B}|_{S_i}$ is PLT at $P$ for $i=1, 2$.
By the existence of $\t{B}|_{S_i}$ and the description of index $2$
PLT surface singularity (see (1.4)), 
we know that $S_i$ is smooth at $P$ and $\diff{S_i}(0)=0$.
Hence by (1.11), $Z$ is also smooth, a contradiction.
\qed
\enddemo
\proclaim{Lemma (7.3)}
(Z, P) is analytically $\Bbb Q$-factorial.
$K_Z+S_1+S_2+B$ is LT outside $P$ and curves where $S_i$ 
and $\supp B$ are simply tangent.
In particular $\ldb{Z}{S_1}{S_2}$ is LT.
\endproclaim
\demo{Proof}
If $(Z,P)$ is not $\Bbb Q$-factorial, by (1.13), we can
take an analytically $\Bbb Q$-factorial LT model 
$\pi:\t{Z}\to Z$ for $K_Z+S_1+B$ near $S_1$. 
By (i) of $(*)$, 
$K_{Z_1}+S_1+B$ is PLT so this modification is a small morphism.
Since $S_1$, $S_2$ and $B$ is $\Bbb Q$-Cartier, the fiber over $P$ is
contained in the inverse images of $S_1$, $S_2$ and $B$.
But this contradicts (1.2).
We will prove the rest.
Let $P'\not=P$ be a point on $Z$ 
and $Z'$ a general hyperplane section of $Z$ through $P'$.
On $Z'$, we express the restricted divisors with $'$.
Note that the property (iv) of (*) holds for $K_{Z'}+S'_1+S'_2+B'$ and $P'$.
If $P'$ is singular on $Z'$ and 
$K_{Z'}+S'_1+S'_2+B'$ is not LT, 
there is a log crepant curve $F$ 
for $K_{Z'}+S'_1+S'_2+B'$ such that 
$d(F, S'_1+S'_2)\leq1$ by (1.3), a contradiction.
If $P'$ is nonsingular on $Z'$ and $K_{Z'}+S'_1+S'_2+B'$ is not LT, 
then by (iv) of (*) and (1.4), 
$P'$ is contained in only one $S'_i$ and $S'_i$ is simply tangent
to $\text{Supp}\ B'$ at $P'$.  
So we obtain the result.
\qed
\enddemo
\proclaim{Lemma (7.4)}
$K_Z+S_1+S_2$ is $1$-complemented.
Assume that $Z$ is not smooth.
Let $\overline{B}$ be an $1$-complement.
Then possibly after replacing $B$ while keeping $(*)$,
$a_l(K_Z+S_1+S_2+B)\geq a_l(K_Z+S_1+S_2+\overline{B})$ holds for any
exceptional divisor.
\endproclaim
\demo{Proof}
By (7.3), $\ldb{Z}{S_1}{S_2}$ is LT.
Hence there exists a good resolution
$r: Z'\to Z$ for $\ldb{Z}{S_1}{S_2}$ 
such that all log discrepancy for $\ldb{Z}{S_1}{S_2}$ are positive
and $\-{r}S_1\cap \-{r}S_2\not = \phi$.
Set $\ldc{Z'}{\-{r}S_1}{\-{r}S_2}{G}=\pb{r}{Z}{S_1}{S_2}$.
Then all the coefficients in $G$ are $<1$ and 
the left hand side (we call this $D$) is LT (subboundary version).
Let $H$ be a very ample divisor on $Z'$ and $r^* r_* H := H+\sum g_iG_i$,
where $G_i$'s are exceptional divisors.
Then $-\sum g_iG_i$ is $r$-ample.
Hence for $\epsilon << 1$, 
$D':=-(\ldc{Z'}{\-{r}S_1}{\-{r}S_2}{G+\epsilon\sum g_iG_i})|_{\-{r}S_1}$
is LT and $r|_{\-{r}S_1}$-ample.
Furthermore $\-{r}S_2|_{\-{r}S_1}\not = \phi$ and is not contained in
a fiber of $r$. 
So by [FA, 19.11, Corollary], 
there is a $1$-complement for $D'$ near any fiber
intersecting $\-{r}S_2|_{\-{r}S_1}$.
Hence there is a  $1$-complement for $\ldb{Z}{S_1}{S_2}|_{S_1}$.
This is lifted to a $1$-complemented for $\ldc{Z'}{\-{r}S_1}{\-{r}S_2}{G}$
by inversion of adjunction.

Next we assume that $Z$ is not smooth.
By (7.2), if we replace $B$ by a more general element 
such that $a_l(K_Z+S_1+S_2+B)\geq a_l(K_Z+S_1+S_2+\overline{B})$ 
for any exceptional divisor,
the conditions of $(*)$ are preserved.
Hence we are done.

\qed
\enddemo

The next lemma is a key to the proof of (7.1).
\proclaim{Lemma (7.5)}
Assume that
 $Z$ is not smooth.
Let $F$ be a exceptional divisor over $P$ 
(not over a curve through $P$) with $a_l(F, K_Z+S_1+S_2)\leq 1$.
Then $F$ satisfies one of the following:
\roster
\item $d(F, S_1+S_2)>1$;
\item 
there is a log crepant divisorial extraction $Z'\to Z$
for $K_Z+S_1+S_2+B$ whose exceptional divisor $F'$ 
contains a generically ODP curve singularity of $Z'$ and 
$F$ is extracted by the simple blow up of the curve singularity.
$a_l(F, K_Z+S_1+S_2+B)=\frac12$ and 
$d(F, S_1+S_2)=\frac12 d(F', S_1+S_2)>\frac12$.
\endroster
\endproclaim
\demo{Proof}
By (ii) of $(*)$ and the assumption that the center of $F$ is $P$, 
$a_l(F, K_Z+S_1+S_2+B)=0 \ \text{or} \ \frac12$.
If the former case occurs, $F$ satisfies $(1)$ by (iv) of $(*)$.
We assume that the latter case occurs.
Let $h:Z_t\to Z$ be a $\Bbb Q$-factorial LT model
for $K_Z+S_1+S_2+B$ near $S_1 \cup S_2$
and $F_j$'s exceptional divisors.
$F$ lands on some $F_j$.
Let $(S_i)_t$ be the strict transform of $S_i$.

\definition{Case 1. $F$ is contracted to a point $z$ on $Z_t$}
Then $z$ is contained in at most two of $F_j$'s and $(S_i)_t$'s.
In fact, if $z$ is contained in three of $F_j$'s and $(S_i)_t$'s,
$z$ is smooth and is not contained in any other component
of the boundaries of $h_t^*(K_Z+S_1+S_2+B)$ 
by (1.12).
Hence the index of $h_t^*(K_Z+S_1+S_2+B)$ at $z$ is $1$.
But this contradicts $a_l(F, K_Z+S_1+S_2+B)=\frac12$.
So we divide into two cases as follows:
\enddefinition

\definition{Subcase 1.1. $z$ is contained in two of $F_j$'s and $(S_i)_t$'s}
Let such two components $D_1$ and $D_2$.
Near $z$, $\pc{h}{Z}{S_1}{S_2}{B}|_{D_i}=
\ldc{D_i}{\diff{D_i}(0)}{D_{2-i}|_{D_i}}{\-{h}B|_{D_i}}$ holds.
It is LT and of index $2$ at $z$.
By (1.4) and (1.12),
we easily deduce that
$z\in Z$ is a $\frac12(1,1,1)$ singularity if $z\not \in \-{h}B$
or $z$ is smooth point of $Z$ if $z\in \-{h}B$.
In any case $D_1+D_2$ is Cartier at $z$.
On the other hand
$d(D_1+D_2, S_1+S_2)>1+1=2$ by the assumption (iv) of $(*)$.
Hence $d(F,S_1+S_2)>1$. So $F$ satisfies $(1)$.
\enddefinition

\definition{Subcase 1.2. $z$ is contained in only $F_j$}
We prove that this case does not occur. 
Take a $1$-complement $\o{B}$ as in (7.4).
Since we may assume that $B$ is general more than $2\o{B}$ (see (7.4)),
$a_l(F, \ldc{Z}{S_1}{S_2}{\o{B}})\leq \frac 12$,
which in turn show that
$a_l(F, \ldc{Z}{S_1}{S_2}{\o{B}})=0$ because the index of 
$\ldc{Z}{S_1}{S_2}{\o{B}}$ is $1$.
By the similar way, we know that
$a_l(F_j, \ldc{W}{S_1}{S_2}{\o{B}})=0$.
Assume that $z\not \in \supp \-{h}\o{B}$.
Then $\lda{W}{F_j}$ is PLT and of index $1$
so $a_l(F, \ldc{W}{S_1}{S_2}{\o{B}})\geq 1$, a contradiction.
Assume that $z\in \supp \-{h}\o{B}$ and let $D:=\-{h}\o{B}\cap F_j$.
Then $$\ldc{W}{\-{h}S_1}{\-{h}S_2}{\-{h}\o{B}+F_j}|_{F_j}=
\ldb{F_j}{D'}{D} \ \text{near} \ z,$$
where $D'$ is $\diff{F_j}(0)- (\text{contribution of} \ D)$. 
But the index of this log divisor is $1$ and $\lda{W}{F_j}$ is PLT
so $D'=0$.
Hence the index of $\lda{F_j}{D}$ is $1$ and 
so $\lda{F_j}{D}$ is PLT.
Hence  $a_l(F, \ldc{Z}{S_1}{S_2}{\o{B}})\geq 1$, a contradiction.
Now we exclude this case.
\enddefinition

\definition{Case 2. $F$ is contracted to a curve $l$}
By (1.2), $l$ is contained
in at most two of $F_j$'s and $(S_i)_t$.
If $l$ is contained in $F_j$ and $F_k$ (or $(S_i)_t$), 
then $l$ is not contained
in any component of $\text{Supp} \ B_t$. 
Hence the index of $h_t^*(K_Z+S_1+S_2+B)$ along $l$ is $1$, a contradiction
to $a_l(F, K_Z+S_1+S_2+B)=\frac12$.
So $l$ is contained in only one $F_j$. 
For a general point $P'$ on $l$ and a 
general hyperplane section $Z'_t$ of $Z_t$ through $P'$,
$K_{Z'_t}+B'_t+F_j:=h_t^*(K_Z+S_1+S_2+B)|_{Z'_t}$
is LT of index 2 and not of index 1 by the existence
of $F$ whose log discrepancy is $\frac12$.
So by (1.4),
$P'\not\in \supp {B'_t}$ and $P'$ is an ODP on $Z'_t$ or
$P'\in \supp {B'_t}$ and $P'$ is smooth on $Z'_t$. 
Note that $$d(F, S_1+S_2)=d(F, (S_1)_t+(S_2)_t+d(F_j, S_1+S_2)F_j)
=d(F, F_j)d(F_j, S_1+S_2).\tag 7.5.1$$ 
For the latter case, $d(F, F_j)$ is an integer.
Hence $F$ satisfies $(1)$ by $(7.5.1)$.
For the former case, if $d(F, S_1+S_2)\leq 1$ then $d(F, F_j)=\frac12$ 
by (7.5.1) and hence
$F$ is the exceptional divisor of the simple blow up of $l$. 
So $F$ satisfies $(2)$.
\enddefinition
Now we finish the proof of (7.5).
\qed
\enddemo

\demo{Proof of Theorem (7.1)}
Since $\ldb{Z}{S_1}{S_2}$ is LT and $S_1$ and $S_2$ are 
$\Bbb Q$-Cartier, we have $(P\in S_1\cup S_2 \subset Z)\simeq
(o\in (xy=0) \subset \frac1s(t_1,t_2,t_3))$, where $t_i$ is an nonnegative
integer less than $s$ and $(s,t_1,t_2)=1$ by (1.12).

First we treat the case $P$ is an isolated singularity.
\proclaim{Claim 1}
If $P$ is an isolated singularity, then $P$ is a terminal singularity.
\endproclaim
\demo{Proof}
We may assume that $P$ is singular.
Let $F$ be any exceptional divisor over $P$.
If $a_l(F, \ldb{Z}{S_1}{S_2})\leq 1$, 
then by (7.5), $a_l(F, K_Z)=a_l(F, \ldb{Z}{S_1}{S_2})+d(F, S_1+S_2)>1$.
Hence for any $F$, we have $a_l(F, K_Z)>1$ and we are done.
\qed
\enddemo

If $P$ is smooth on $Z$, we obtain (1) by (1.12).
Assume below that $P$ is singular point of $Z$.
Then by (1.12), the index of $P$ is greater than $1$.
Hence by Kawamata's minimal discrepancy theorem [Ka 2],
there is a exceptional divisor $F$ over $P$ with 
$a_l(F,K_Z)=1+\frac1r$, 
where $r$ is the index of $P$.
Since $P$ is terminal, $a_l(F, \ldb{Z}{S_1}{S_2})\in \frac {\Bbb Z}r$
and $a_l(F, \ldb{Z}{S_1}{S_2})<1+\frac1r$ by $a_l(F,K_Z)=1+\frac1r$.
Hence $a_l(F, \ldb{Z}{S_1}{S_2})\leq 1$ and we apply (7.5).
If $F$ satisfies (1) of Lemma (7.5), 
then $d(F, S_1+S_2)\geq 1+\frac1r$ since $P$ is terminal.
On the other hand, $a_l(F,\ldb{Z}{S_1}{S_2})>0$ since 
$a_l(F,\ldc{Z}{S_1}{S_2}{B})\geq 0$ and the center of $F$ on $Z$ is $P$.
But  these contradict  
$$d(F, S_1+S_2)+a_l(F,\ldb{Z}{S_1}{S_2})=1+\frac1r. \tag 7.1.1 $$
So $F$ satisfies (2).
By $a_l(F, K_{Z}+S_1+S_2+B)=\frac12$, 
$a_l(F, K_{Z}+S_1+S_2)>\frac12$. On the other hand $d(F,S_1+S_2)>\frac12$.
By these and (7.1.1),
$r$ must be odd (we set $r=2r'+1$) and 
$a_l(F, \ldb{Z}{S_1}{S_2})=d(F, S_1+S_2)=\frac{r'+1}{2r'+1}$.
Furthermore by the theorem of Morrison-Stevens [MS, Theorem 2.4]
and that $S_1+S_2$
is not a Cartier divisor (since $d(F, S_1+S_2)$ is not an integer),
we can take $s=2r'+1$, $t_1+t_3=2r'+1$ and $t_2=1$.
By the weighted blow up with a weight $\frac1{2r'+1}(t_1,1,2r'+1-t_1)$,
$F$ is extracted by the uniqueness of the exceptional divisor with
minimal discrepancy for terminal quotient singularity (see [Ka3]). 
So $\frac{r'+1}{2r'+1}=d(F, S_1+S_2)=\frac{t_1+1}{2r'+1}$.
Hence $t_1=r'$. 
We prove $r'=1$, 
which in turn show that (2) of (7.1) holds with a weight change.
It suffices to prove that if $r'\geq 2$, 
$F$ must be log crepant for $\ldc{Z}{S_1}{S_2}{B}$ (this is a contradiction).
Let $F''$ be the curve 
which is extracted by the restriction of the above weighted blow up to $y=0$.
Let $F_1,F_2,...., F_{2r'}$ be the exceptional curve of the MRS
of the singularity of the origin on $y=0$, where $F_1$ intersects the strict
transform of $x=y=0$.
Then $F''$ is $F_{r'+1}$ on the MRS by the choice of weight.
Hence if $r'\geq 2$, 
$F''$ must be log crepant for $\ldc{Z}{S_1}{S_2}{B}|_{S_2}$ 
by (1.3),
which in turn show that $F$ is log crepant for $\ldc{Z}{S_1}{S_2}{B}$,
a contradiction. The latter half of (2) is an easy calculation.

Next we treat the case that $P'$ is non-isolated singularity.
We will prove that $P$ satisfies (3) of (7-1).
Since no curve singularity through $P'$ is contained in 
$\llc(\ldc{Z}{S_1}{S_2}{B})$ by (7.3) and 
any curve singularity through $P'$
is a generically ODP curve since the index of
$\ldc{Z}{S_1}{S_2}{B}$ is $2$.
\comment
By this we know that there is only one curve singularity through $P'$
since otherwise $(s,t_1,t_2)\not =1$.
\endcomment
Let $\t{Z}:=\pmb{Spec} (\Cal O_Z \oplus \Cal O_Z(-S_1 -S_2))$.
Then $\pi:\t{Z}\to Z$ is a double covering such that $\t{P}:=\-{\pi}P$
is an isolated singularity on $\t{Z}$.
It is easy to see that $\ldc{\t{Z}}{\t{S_1}}{\t{S_2}}{\t{B}}:= 
\pc{\pi}{Z}{S_1}{S_2}{B}$ satisfies the condition (*). 
If $\t{P}$ satisfies $(1)$, we obtain what we want.
We show that $\t{P}$ does not satisfies $(2)$.
Assume that $\t{P}$ satisfies $(2)$.
Then it is easy to see that 
locally analytically $(P\in S_1\cup S_2 \subset Z)\simeq
(o\in (xy=0) \subset \frac16(5,2,1))$.
Take the weighted blow up $h: W\to Z$ with the weight $\frac16(5,2,1)$ 
and let $F$ be the exceptional divisor.
The restriction of $h$ to $S_1$ is the MRS of $P$ 
which extracts a $(-3)$-curve 
so $h$ is log crepant for $\ldc{Z}{S_1}{S_2}{B}$ by (1-3).
Then at the intersection of $F$ and $\-{h}S_2$,
there is a singularity $Q$ such that
$(Q\in F\cup \-{h}S_2 \subset W) \simeq
(o\in (xy=0) \subset \frac15(2,1,3))$.
Since the index of $\pc{h}{Z}{S_1}{S_2}{B}$ is $2$ at
the singularity, $\supp \-{h}B$ passes through it and neither 
$\pc{h}{Z}{S_1}{S_2}{B}|_F$
or $\pc{h}{Z}{S_1}{S_2}{B}|_{\-{h}S_2}$
is PLT.
Hence the weighted blow up with the weight $\frac15(2,1,3)$
is log crepant for $\pc{h}{Z}{S_1}{S_2}{B}$ since the restriction
of this weighted blow up to $\-{h}S_2$ is log crepant for
$\pc{h}{Z}{S_1}{S_2}{B}|_{\-{h}S_2}$ by (1.3).
Let $F'$ be the exceptional divisor of this weighted blow up.
Then $d(F', S_1+S_2)=d(F', \-{h}S_2+\frac76 F)=\frac23$.
But this contradicts the property of $(*)$.

\qed
\enddemo

\comment
\proclaim{Lemma}
while BTM, the flips of the strict transform of $C$.....
the non-increasing of LLC.
the description of $B$.
\endproclaim
\demo{Proof}
By an abuse of notation, we also call $C$, $D_1$ and $D_2$ 
the strict transform of $C$, $D_1$ and $D_2$.  
In a neighborhood of $C$ on $Z$, $C= \supp B_1|_{S_1}$ 
and $C$ is contractible on $S_1$. Hence there exists a component $C'$ 
of $C$ such that $B_1.C'<0$. 
On the other hand, for a $1$-dimensional fiber $l$ for $g$, $B_1.l>0$,
so we must have $B_1.C_1<0$. Hence $C_1\subset S_1\cap \supp B_1$,
which in turn shows that $C_1\subset C$. We can easily deduce that
$E_1.C_1>0$ and $S_1.C_1<0$ by $B_1.C_1<0$  
so the flip of $C_1$ is a beginning flip.
After the flip of $C_1$, 
Any component of $C^+_1$ is not contained in $S_2$ by Lemma ?.
If some components of $C$ left on $S_2$, then as above there exists
a component $C"$ of $C$ such that $B_2.C"<0$.
On the other hand, for the flipped curve $C^+_1$, $B_2.C^+_1>0$,
so we must have $B_2.C_2<0$. (Note that $g_2$ cannot be a divisorial
contraction.) Hence as above we see that $C_2\subset C$
and any component of $C^+_2$ is not contained in $S_3$.
By the numerical property (1) in (2.5),
we can easily see that $E_2.C_2>0$ and $S_2.C_2<0$ so the flip of $C_2$
is a beginning flip.
Since the number of the components of $C$ on $X$ is at most $2$,
any component of $C$ lefts on the strict transform of $S$.
\qed
\enddemo
\endcomment

\proclaim{Lemma (7.6) (Reducible case)}
Consider the case $(c)$.
Let $g: Z\to X$ be a primitive log crepant divisorial extraction for
$\ldb{X}{S}{B}$ and $E$ the exceptional divisor.
Assume the following conditions:
\roster
\item $\-{g}S\cap E$ is reducible 
(and hence by Lemma (5.2), we can write $\-{g}S\cap E=D_1\cup D_2$);
\item $D_1\simeq \Bbb P^1$. $\-{g}S$ and $E$ are generically normal crossing
along $D_1$ (Note that we assume nothing about $D_2$ here.);
\item $\ldc{Z}{\-{g}S}{\-{g}B}{E}|_{D_1}=
\ldc{D_1}{\frac12 P_1}{\frac12 P_2}{P}$;
\item $\-{g}S$ is normal.
\endroster

Run the BTM starting by $g$. 

Then 
$D_1$ becomes a flipping ray satisfying the assumptions of (2.4)
after the flip of all the strict 
transforms of the components of $C$.
It satisfies the following numerical properties:
$B'.D_1=0, S'.D_1<0 \ \text{and} \ E'.D_1>0$,
where $B'$, $S'$ and $E'$ are the strict transform of $B$, $S$ and $E$
(by abuse of notation, we express the strict transform of $D_1$ also
by $D_1$).
In particular $g$ is a semi-good extraction by the property (2) of
(2.4).
\endproclaim

\demo{Proof}
We use the notation as in (2.5).
In particular we write $S_1$, $B_1$ and $E_1$ instead of $\-{g}S$, $\-{g}B$
and $E$.
By an abuse of notation, we also call $C$, $D_1$ and $D_2$ 
the strict transforms of $C$, $D_1$ and $D_2$.  
In a neighborhood of $C$ on $Z$, $C= \supp B_1|_{S_1}$ 
and $C$ is contractible on $S_1$. Hence there exists a component $C'$ 
of $C$ such that $B_1.C'<0$ and $E_1.C' >0$. 
On the other hand, for an exceptional curve $l$ for $g$, 
$B_1.l>0$ and $E_1.l<0$,
so we must have $B_1.C_1<0$, $E_1.C_1>0$ and $S_1.C_1 <0$. 
Hence $C_1\subset S_1\cap \supp B_1$,
which in turn shows that $C_1\subset C$. By the above numerical property,
$C_1$ is a beginning flipping ray and the flip exists.

\proclaim{Claim}
After the flip of $C_1$, 
Any component of $C^+_1$ is not contained in $S_2$.
\endproclaim
\demo{Proof}
If some component $(C^+_1)'$ is contained in $S_2$,
it contained in $S_2\cap E_2$. 
In particular $(C^+_1)'\subset \llc(\ldc{Z_2}{S_2}{B_2}{E_2})$.
Hence $C_1$ contains an element of  $\text{CLC} \ (\ldc{Z_1}{S_1}{B_1}{E_1})$
since $(\ldc{Z_2}{S_2}{B_2}{E_2}).C^+_1=0$.
But this contradicts the assumption of case (c) and (3) of this lemma.
\qed
\enddemo

Note that $E_2$ must be normal since $E_1$ is normal by Remark (5.3)
and the proof of Claim.
If some components of $C$ left on $S_2$, then as above there exists
a component $C''$ of $C$ such that $B_2.C''<0$ and $E_2.C''>0$.
On the other hand, for the flipped curve $C^+_1$, 
$B_2.C^+_1>0$ and $E_2.C^+_1<0$,
so we must have $B_2.C_2<0$, $E_2.C_2>0$ and $S_2.C_2<0$.  
Hence as above we see that $C_2\subset C$
and any component of $C^+_2$ is not contained in $S_3$.
By the above numerical property,
we can easily see that $E_2.C_2>0$ and $S_2.C_2<0$ so the flip of $C_2$
is a beginning flip. (After the flip, $E_3$ is normal.)
Since the number of the components of $C$ on $X$ is at most $2$,
no component of $C$ lefts on the strict transform of $S$.
After the flips of the strict transform of $C$, we call $Z'$ the image of $Z$
and $S'$, $B'$ and $E'$ the strict transform of $S$, $B$ and $E$.
Let $R$ be the next extremal ray.
We show that $R$ is a flipping ray stated as in the statement of this lemma.
We claim that first $B'.D_1=0$. 
After the flips of the strict transform of $C$, 
$\supp B'$ intersects $S'$ at most in finite points near $D_1$
since the flipped curves are not contained in $S'$
but by $\Bbb Q$-factoriality, $\supp B'\cap S'$ must be empty near $D_1$.
Hence we obtain $B'.D_1=0$.
On the other hand, $B'$ is positive for the flipped ray one before $R$.
Hence we must have $B'.R\leq 0$.

We prove that $B'.R<0$ does not occur. Assume that $B'.R<0$.
Then we see that neither $D_1$ or $D_2$ is not contained in $\supp R$
since neither $D_1$ or $D_2$ is contained in $\supp B'$.
Hence by the numerical property (2) in (2.5),
we must have $S'.R>0$.
Hence by the numerical property (1) in (2.5),
$E'.R<0$, i.e., $R$ is a final flipping ray (it is not divisorial
since $B'.R<0$). Note that $D_1 \cap \supp R = \phi$ since $\supp R \subset
\supp B'$ and $D_1 \cap \supp B = \phi$.
We see inductively that after the flip of $R$,  
an extremal ray is negative for $B'$ and $E'$ and positive for $S'$
as long as $D_1$ is not a flipping curve while BTM
(hence it is not a divisorial ray and its support does not intersects $D_1$).
Hence after finite number of such a flip, an extremal ray $R''$ such
that $B''.R''=0$ appears.
We show  that $R''$ is not divisorial.
If $R''$ is divisorial, after the contraction we obtain the flip of $C$. 
Let $C^+$ be the flipped curve. Then $B^+.C^+<0$ since $B.C>0$. 
If $E''$ is contracted to a point, $C^+$ passes through the contracted
point.
But by $B''.R''=0$, $C^+$ cannot be contained in $\supp B^+$, a contradiction.
If $R''$ is contracted to a curve, the image curve is not contained in 
$\supp B^+$ since $B''.R''=0$.
Hence the center of $E$ on $X$ cannot be $L$ and the image of $E''$ is 
contained in $C^+$. But this again contradicts $B^+.C^+<0$.

Hence $R''$ is not divisorial and $D_1$ is a flipping curve (I don't know
if $\supp R'' =D_1$.) Here we consider the flipping contraction of $R''$
in the analytic category. First we flip only $D_1$. 
Since $E'$ is normal and above procedures 
after the flip of the strict transform of $C$
are isomorphic near $D_1$, $E''$ is normal near $D_1$.
Hence we can check that $D_1$ satisfies
the assumption of (2.4).
Hence by (2.4), after the flip, $D_1$ go outside $E''$.
Next we flip other components of $\supp R''$. While such flips,
$D_1$ never become a flipping curve again.
In fact the numerical property (1) in (2.5) 
still holds in the analytic category.
While such flips, $E''.D_1$ continues to be positive and hence $S''.D_1$ 
continues to be negative. So $(\lda{Z''}{E''}).D_1$ continues to be positive.
But such flipping curves (other component of $\supp R''$) are negative
for $\lda{Z''}{E''}$. Hence we obtain the assertion.

After the flip of $R''$ (and possibly the flips of some final flipping rays), 
finally $D_1$ become a curve contained in $C^+$. But
$D_1$ is not contained in $\supp B^+$, a contradiction to $B^+ . C^+ <0$.

Hence we have $B'.R=0$. we immediately see that $\supp R \subset D_1$.
By the numerical property (2) in (2.5),
we must have $S'.R>0$ or $E'.R>0$.
We can show the former case cannot occurs just as above.
Hence the latter case occurs and 
We can easily check $R$ satisfies the assumptions of (2.4)
by the assumptions of this Lemma.
\qed
\enddemo

\proclaim{Lemma (7.7) (Reducible case 1)}
Consider the situation as in (7.6).
Assume furthermore that $g(E)=Q$.

Then we have the description as explained in the following flow chart:
\endproclaim

\newpage
$$\hskip1in
\fbox{2.5in}{$D_1$ become a flipping ray of the type in (2.4).
Consider the next extremal ray $R'$}
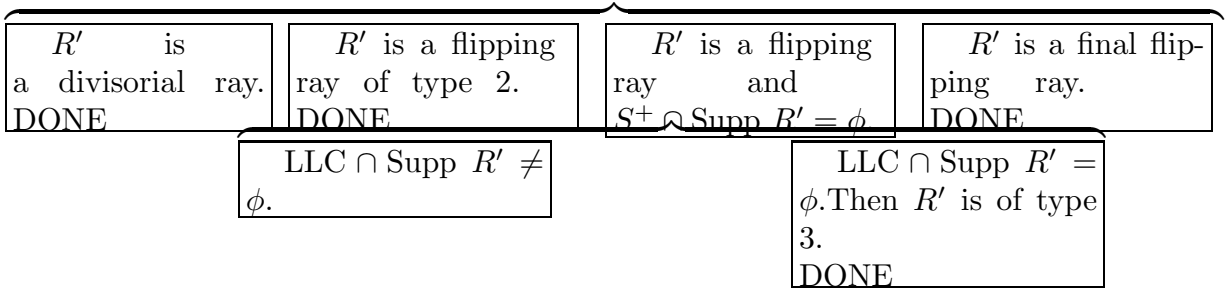
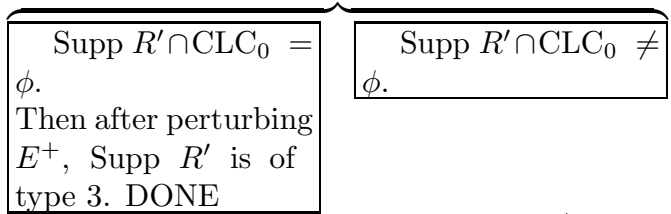
$$

$$\oversetbrace\to{
\hbox{\hsize=4.7in
\fbox{1.1in}{$R'$ is \break
a divisorial ray.\break
 DONE}\hskip5pt
\fbox{1.2in}{$R'$ is a flipping \break
ray of type 2. \break
DONE}
\hskip5pt
\fbox{1.2in}{$R'$ is a flipping \break
ray and \break
$S^+ \cap \supp R' = \phi$.}
\hskip5pt
\fbox{1.2in}{$R'$ is a final flipping ray. 
\break
DONE}
}}
$$
\vskip1.0pt

$$\hskip0.5in\oversetbrace\to{
\hbox{\hsize=4in
\fbox{1.3in}{
$\text{LLC} \cap \supp R' \not = \phi$.}
\hskip1.0in
\fbox{1.3in}{
$\text{LLC} \cap \supp R' = \phi$.Then $R'$ is of type 3. \break
DONE}}}
$$
\vskip0.5in
$$
\oversetbrace\to{
\hbox{\hsize=2.6in
\fbox{1.3in}{$\supp R' \cap \clc{0}= \phi$.
\break Then after perturbing $E^+$, $\supp R'$ is of \break
type 3. DONE}
\hskip.1in
\fbox{1.3in}{
$\supp R' \cap \clc{0} \not = \phi$.}}}
\hskip1.0in
$$
\vskip0.5in
$$\hskip1.0in\downarrow$$
$$\fbox{4.7in}{\roster
\item $Q':=D_2\cap \-{g}L$ is smooth on $\-{g}S$, $E$ and $Z$.
$(D_2)^2_E=2$.
$\supp B_1|_E$ and $D_2$ are simply tangent at $Q'$;
\item $P$ is smooth on $\-{g}S$, $E$ and $Z$;
\item $d=\frac{3}{-(D_2)_{\-{g}S}^2 -1}$;
\item 
$d(F, S)=d+1$ and $a_l(F, S)=0$ 
for the exceptional divisor $F$ of the simple blow up along $D_1$ or $D_2$
and for any other log crepant divisor $F\not= E$ 
for $\ldb{X}{S}{B}$, $d(F, S)>d+1$.
\endroster}
$$

\vskip1.5in
$$\llc :=\llc(\ldc{Z^+}{S^+}{B^+}{E^+}|_{E^{+\nu}}).$$
$$\clc{0}:=\clc{0}(\ldc{Z^+}{S^+}{B^+}{E^+}|_{E^{+\nu}}).$$

See also FIGURE (XXI).

\demo{Proof}
When $D_1$ becomes a flipping ray, we use the same notation as in (7.6)
($S'$, $E'$.., etc).
After the flip of $D_1$,
we call the flipped curve $D^+_1$ and $S^+$, $B^+$, $E^+$ and $D^+_2$ 
the strict transforms of $S$, $B$, $E$ and $D_2$.
Let $R'$ be the next extremal ray.
Since $S^+\cap E^+$ is irreducible,
we can assume that $\supp R'\cap S^+ =\phi$ and $E^+.R'<0$ 
(see the flow chart).
Furthermore 
we can assume that
$\supp R'\cap \llc(\ldc{Z^+}{S^+}{B^+}{E^+}|_{E^{+\nu}})\not= \phi$, 
where $\nu:E^{+\nu}\to E^+$ is the normalization of $E^+$. 
Let $D^* = \-{\mu}D_1^+$.
\proclaim{Claim 1}
$\llc(\ldc{Z^+}{S^+}{B^+}{E^+}|_{E^{+\nu}})=D^*\cup D^+_2$.
\endproclaim
\demo{Proof} 
If $\llc(\ldc{Z^+}{S^+}{B^+}{E^+}|_{E^{+\nu}})$ has 
another component except $D^*$ and $D^+_2$,
$\llc(\ldc{Z}{S}{B}{E}|_E)$ has another component $D'$ 
except $D_1$ and $D_2$, and 
$D'$ must intersect $D_2$ at a point ($\not = P$)
by (3.4). Furthermore 
$\llc(\ldc{Z^+}{S^+}{B^+}{E^+}|_{E^{+\nu}})=D^*\cup D^+_2\cup D^{'+}$,
where $D^{'+}$ is the strict transform of $D'$, and
$D^*\cap D^{'+}=\phi$ and $D^+_2\cap D^{'+}\not=\phi$.
Note that $(D^+_2)^2_{E^{+\nu}}=S^+.D_2>0$ since $S^+.R'=0$ and
$S^+.D^+_1>0$.
Hence 
$aD^+_2-(\lda{E^{+\nu}}{\{\diff{E^{+\nu}}({S^+}+{B^+})\}})
\equiv aD^+_2 +\llcorner \diff{E^{+\nu}}({S^+}+{B^+})\lrcorner$
is nef and big for $a>>0$. Furthermore  
$\lda{E^{+\nu}}{\{\diff{E^{+\nu}}({S^+}+{B^+})\}}$ is KLT 
by Theorem (1.7) (2),
so by Kawamata's base point free theorem, $D^+_2$ is semi-ample.
\comment
Let $l$ be a connected component of $\supp R'$ such that
$l \cap \llc(\ldc{Z^+}{S^+}{B^+}{E^+}|_{E^{+\nu}})\not= \phi$ and
$m$ the fiber of the morphism defined by a sufficient multiple of $D^+_2$ 
containing $l$. 
\endcomment

On the image $\o{E}$ of $E^{+\nu}$ 
by the morphism defined by a sufficient multiple of $D^+_2$,
take a 
$\lda{\o{E}}{(1-\epsilon)(\o{\diff{E^{+\nu}}({S^+}+{B^+})}
-\o{D^{'+}}})$-extremal ray
$A$ ($0<\epsilon <<1$ is a sufficiently small rational number) such that  
$\o{D^{'+}}.A>0$. We can take such a ray since $\o{D^{'+}}^2 >0$
(cf. Claim 3 in the proof of (1.9)). 
Note that $\o{D_2^+}$ is ample.
If $A$ is birational, by $\o{D^{'+}}.A>0$ and the ampleness of $\o{D^+_2}$,
$\supp A$ intersects $\o{D^+_2}$ and $\o{D^{'+}}.A>0$.
Hence by (1.6), $\supp A \cap \o{D^+_2} \cap \o{D^{'+}}$
is one point. We take this procedure till the extremal ray is not birational.
Note that the situation does no change while this procedure.
So we may assume that 
the image of $\o{E}$ by the extremal contraction associated to $A$ is a curve. 
Since $\o{D^{'+}}.A>0$, $\o{D^*}$ is not contained in a fiber.
(Note that 
$\o{D^*}$ does not intersect $\o{D^{'+}}$ by applying (1.6)
to the morphism defined by sufficiently multiple of $D^+_2$.)
Hence a general fiber of the contraction of $A$ intersects 
$\llc(\ldc{Z^+}{S^+}{B^+}{E^+}|_{E^{+\nu}})$ with at least three points.
But this is impossible since for a general fiber $r$, $K_{E^{+\mu}}.r = -2$.

The contraction of $A$ cannot contract $\o{E}$ to a point since
$\o{D^*}$ and $\o{D^{'+}}$ does not intersect.
 
This is a final contradiction and we finish the proof of Claim 1.
\qed
\enddemo
\definition{Subcase 1}
$\clc{0}(\ldc{Z^+}{S^+}{B^+}{E^+}|_{E^{+\nu}})\cap \supp R'=\phi$
\enddefinition
Let $l$ be any connected component of $\supp R'$.
In this case we prove that after perturbing $E^+$ near $l$,
the flip of $l$ become a flip of type 3.

By assumption of Subcase 1,
$\ldc{Z^+}{S^+}{B'}{E^+}|_{E^{+\nu}}$ is exceptional near $l$.
By (2.4) (4), the nonnormal locus of $E^+$ is $D^*$.
Hence $l$ is not contained the nonnormal locus of $E^+$ and
taking a more general element $E'$ than $E^+$ in $|E^+\supset l|$ near $l$,
$E'$ is normal.
Since $\llc(\ldc{Z^+}{S^+}{B^+}{E^+}|_{E^{+\nu}})$ is not contained in 
$\supp B^+|_{E^{+\nu}}$ near $l$, 
$\ldc{Z^+}{S^+}{B'}{E'}|_{E'}$ has at most isolated LLC near $l$.
If $\ldc{Z^+}{S^+}{B'}{E'}|_{E'}$ has isolated LLC near $l$,
$\llc(\ldc{Z^+}{S^+}{B^+}{E^+}|_{E^{+\nu}})$ also has isolated LLC
near $l$. But this contradicts the assumption of Subcase 1.

Then $\ldc{Z^+}{S^+}{B'}{E'}|_{E'}$ is KLT and by the inversion of
adjunction, 
$\ldc{Z^+}{S^+}{B'}{E'}$ is PLT.

\definition{Subcase 2}
$\clc{0}(\ldc{Z^+}{S^+}{B^+}{E^+}|_{E^{+\nu}})\cap \supp R'\not=\phi$
\enddefinition
In this case we will prove (1), (2), (3) and (4) of this Lemma.
We first prove (1).
Since $D^*\cap D^+_2$ is contained in
$\clc{0}(\ldc{Z^+}{S^+}{B^+}{E^+}|_{E^{+\nu}})$,
$D^*$ has exactly one point in
$\clc{0}(\ldc{Z^+}{S^+}{B^+}{E^+}|_{E^{+\nu}})$ other than $D^*\cap D^+_2$.
We call this $Q''$.
Let $l$ be the connected component of $\supp R'$ containing $Q''$.

\proclaim{Claim 2}
There is no log crepant primitive extraction $\mu:E'\to E^{+\nu}$ for
$\ldc{Z^+}{S^+}{B^+}{E^+}|_{E^{+\nu}}$ over $Q'$ such that
$(\-{\mu}D^+_2)^2\geq 0$ and $\supp {\-{\mu}(B^+|_{E^{+\nu}})}\cap 
\-{\mu}D^+_2=\phi$.
\endproclaim

\demo{Proof}
Assume that there exists such a $\mu$. 
As in the proof of Claim 1, we can prove that $\-{\mu}D^+_2$ is semi-ample.
Let $m$ be
the fiber of the morphism defined by a sufficient multiple of $\-{\mu}D^+_2$ 
containing $l$.
There are components of $\supp {\-{\mu}(B^+|_{E^{+\nu}}}$ intersecting
$\-{\mu}l$ since $B^+.l>0$.
Let $b$ be the union of such components of 
$\supp {\-{\mu}(B^+|_{E^{+\nu}})}$.
By the assumption for $\mu$,
$\supp \-{\mu}(B^+|_{E^{+\nu}})\cap \-{\mu}D^+_2=\phi$. 
Hence $b$ intersects the exceptional curve for $\mu$ and 
$\-{\mu}l\cup b \subset m$.
But $\-{\mu}l$ contains $Q''$ and $b$ intersects
the exceptional curve for $\mu$, so we obtain a contradiction by
applying (1.6) for the birational contraction
of $m$ after extracting some log crepant curve from $Q''$.
\qed
\enddemo

Assume that $Q'$ is a singular point of $E^{+\nu}$.
Since $Q'$ is the unique singular point of $E^{+\nu}$ on $D^+_2$ 
by Theorem (2.4) (3), $(\-{m}D^+_2)^2 >0$, where $m$ is the MRS along $D^+_2$.
We can make $m$ primitive and satisfying the assumption of Claim 2, 
a contradiction.   
So we see that $Q'$ is a smooth point
on $E^{+\nu}$.
Furthermore we see that
$(D^+_2)^2=1$ and
$\supp B^+|_{E^{+\nu}}$ and $D^+_2$ are simply tangent
at $Q'$.
By the last property, 
$E^+$ contains no curve singularity near $Q'$.
By the property of semi-good extraction,
$\ldb{Z^+}{S^+}{E^+}$ is LT near $Q'$.
Hence by (1.11), $Z^+$ is also smooth at $Q'$.
Since $\ldb{Z^+}{S^+}{E^+}$ is LT at $Q'$
and $E^+$ contains no curve singularity near $Q'$,
$S^+$ is also smooth by (1.12). 
We finish the proof of (1).

Next we prove (2).
It suffices to prove (2) on $Z'$ since 
the flips before $R$ are isomorphism near $P$.
Let $t:\+n{\t{E}}\to \+n{E}$ be MLT of 
$\ldc{Z^+}{S^+}{B^+}{E^+}|_{E^{+\nu}}$ at $Q''$ and
$F_1,..., F_n$ exceptional divisors of $t$, where $F_1$ intersects
$\-{t}D^*$. By (2.4), $E'$ is obtained from $\+n{\t{E}}$ by contracting
$\-{t}D^*, F_1,..., F_{n-1}$.
Then we claim the following:

\proclaim{Claim 3}
$n=1$ and $\-{t}D^*$ is a $(-1)$-curve.
\endproclaim
(2) follows from this claim and (1.11) immediately . 
\demo{Proof of Claim 3}
Let $h:\+n{E}\to \+n{\o{E}}$ be the morphism
defined by a sufficient multiple of $D^+_2$. 
We see that $\+n{\o{E}}\simeq \Bbb P^2$ since $h(D^+_2)\simeq \Bbb P^1$
is an ample Cartier divisor on $\+n{\o{E}}$ and $h(D^+_2)^2=1$ by (1)
and (2).
Also we know that $h(D^*)$ is a line on $\+n{\o{E}}$.
We claim that $\cup^n_{i=1} F_i\cup \-{t}D^*$ intersects strict transforms
of $h$-exceptional curves only at points on $F_n$.
In fact, if the strict transform $F$ of
a $h$-exceptional curve intersects
$\cup^n_{i=1} F_i\cup \-{t}D^*$ at points on 
$F_i$ with $i\not= n$ (resp. $\-{t}D^*$), then
$F$ is not contained in components of 
the boundary of $\pc{t}{Z^+}{S^+}{B^+}{E^+}|_{E^{+\nu}}$
by considering the degree of $\pc{t}{Z^+}{S^+}{B^+}{E^+|_{E^{+\nu}}}|_{F_i}$
and so $F$ is a $(-1)$-curve on the MRS of $Q''$.
Hence also $F$ does not intersect $\supp \diff{\+n{E}}({S^+}+{B^+})$ 
except $F_i$ with $i\not= n$ (resp. $\-{t}D^*$).
By these, the strict transform of $F$ on $E'$ intersects $D_1$
only at $P$ and does not intersect $\supp {B'|_{E'}}$.
But this is a contradiction since $B^-|_{E^-}$ is nef on $E'$
and numerically trivial only for $D_1$.
Hence we obtain the description as follows:

FIGURE (XXII).

Let $\ldc{\+n{\o{E}}}{h(D^*)}{h(D^+_2)}{\o{B}}=
h_*(\ldc{Z^+}{S^+}{B^+}{E^+}|_{E^{+\nu}})$.
and let $s:\+n{\t{E}}\to G$ the contraction of the strict transforms
of $h$-exceptional curves. By Claim 1, they are not log crepant
curves for 
$\ldc{\+n{\o{E}}}{h(D^*)}{h(D^+_2)}{\o{B}}$. 
Hence we obtain the MLT for
$\ldc{\+n{\o{E}}}{h(D^*)}{h(D^+_2)}{\o{B}}$ possibly after contracting
the images of some $F_i$'s on $G$. By the above,
such curves intersect $\cup^n_{i=1} F_i\cup \-{t}D^*$ at points
only on $F_n$. Hence after contraction, the image of 
$\supp \-{(\nu\circ t)}B^+|_{E^+}$ intersects the image of
$\cup^n_{i=1} F_i\cup \-{t}D^*$ at points only on the image of $F_n$.
So if we contract the image of $F_n$, 
$s_*(\pc{t}{Z^+}{S^+}{B^+}{E^+|_{E^{+\nu}}})$ become non LT. 
Hence $G$ is MLT for $\ldc{\+n{\o{E}}}{h(D^*)}{h(D^+_2)}{\o{B}}$.
On the other hand
$\supp \o{B}$ is tangent to $h(D^+_2)$ by (1).
Since $\ldc{\+n{\o{E}}}{h(D^*)}{h(D^+_2)}{\o{B}}$ is numerically trivial,
$\supp \o{B}$ is a smooth conic. Furthermore by the existence of
$Q''$, $\supp \o{B}$ is also tangent to $h(D^*)$.
Hence MLT for the image of $Q''$ is described as in (1.4).
From these, we obtain Claim 3.

\comment
If $n\geq 2$, then $h\circ t$ contains exactly one blow up 
on the strict transform of $D^*$.  
So $(\-{t}D^*)^2=0$, which contradicts the contractibility of $\-{t}D^*$.
Hence $n=1$.
Furthermore we prove that $h\circ t$ contains exactly two blow up 
on the strict transform of $D^*$. This completes the proof of Claim 3.
Note that $h(P')$ is a smooth point on $\+n{\o{E}}$.
On the MRS of points on $F_1$,
its exceptional curves and
$h$-exceptional curves intersecting $F_1$ are $(-2)$-curves
since 
$h$-exceptional curves intersecting $F_1$ are $(-1)$-curves or $(-2)$-curves
and the index of $\ldc{Z^+}{S^+}{B^+}{E^+}|_{E^{+\nu}}$ is $2$. 
Decompose the contraction $h'$ from the MRS as above
to $\+n{\o{E}}$ into the sequence of contractions of $(-1)$-curves.
Then at most one $h'$-exceptional intersects
the strict transform of $D^*$ when  
the strict transform of $D^*$ become a $(-1)$-curve. 
Furthermore
such a $h'$-exceptional is $(-2)$-curve over a singular point on $F_1$
and hence it does not intersect $h'$-exceptional except $F_1$.
This gives what we want.
\endcomment
\qed
\enddemo

For (3), we consider 
$0=(S_1+dE).D_2=\{(D_1+D_2).D_2\}_E+d\{(D_1+D_2).D_2\}_{S_1}$.
By (1) and (2), $(D_1.D_2)_E= 1$, $(D_1.D_2)_{S_1}= 1$, 
$(D_2)^2_E=2$ so (3) follows.
(4) is easily deduced from (1) and (2).
\qed
\enddemo

\proclaim{Lemma (7.8) (Reducible case 2)}
Consider the situation as in (7.6).
Assume furthermore that $g(E)=L$.
\comment
\item $D_1$ becomes a flipping ray as in (7.6)
(i.e., it satisfies the numerical properties 
$B'.D_1=0, S'.D_1<0 \ \text{and} \ E'.D_1>0$).
\endroster
\endcomment

Then 
\roster
\item $P:=D_1\cap D_2$ is smooth on $S_1$ and $Z$ and is an ODP or a smooth
point on $E$.
\item 
$d(F, S)=d+1$ and $a_l(F, S)=0$ 
for the exceptional divisor $F$ of the simple blow up along $D_1$ or $D_2$
and for any other log crepant divisor $F\not= E$ 
for $\ldb{X}{S}{B}$, $d(F, S)>d+1$.
\endroster
\endproclaim
\demo{Proof}
When $D_1$ becomes a flipping ray, we use the same notation as in 
(7.6) ($S'$, $E'$.., etc).

\proclaim{Claim}
$P$ is an ODP or a smooth point on $E$.
\endproclaim

\demo{Proof}
Let $l$ be the fiber of $E\to L$ containing $D_1$ and $l'$
any other component of $l$ intersecting $D_1$. 
Note that $B_1.l'>0$.
So if $l'$ passes through $P$, after contraction of $l'$,
$\ldc{Z}{S_1}{B_1}{E}|_E$ has two reduced boundaries 
$D_1$ and $D_2$ and some other boundary contained in  $\supp B_1$ at $P$,  
a contradiction.
Hence $l'$ does not pass through $P$.
Hence contractions of components of $l$ except $D_1$ does not change
the singularity of $P$ on $E$.
Furthermore by the fact that 
any component of $l$ except $D_1$ is not contained in
$\llc (\ldc{Z}{S_1}{B_1}{E}|_E)$ and 
$\ldc{Z}{S_1}{B_1}{E}|_E$ is numerically trivial for $g$,
the property $(3)$ of Lemma (7.6) does not change.
so we may assume that $D_1$ is an irreducible fiber. 
By this, we easily show that $P$ is an ODP or a smooth point on $E$.
\qed
\enddemo

Assume that $P$ is smooth on $E$.
Since $\ldc{Z}{S_1}{B_1}{E}|_E$ has two reduced boundary 
$D_1$ and $D_2$ at $P$, 
$E$ contains no singular curve near $P$.
So $P$ is smooth on $Z$
by (1.11) and LT property of $\lda{Z}{E}$.
Furthermore by (2.4), on $S_1$, 
$P$ is a smooth point or a singular point of index $\geq 3$. 
So $P$ must be smooth on $S_1$. 

If $P$ is an ODP on $E$, it suffices to prove that $P$ is smooth on $S_1$
(by (1.11)).
If $P$ is singular on $S_1$,
$P$ is a singular point of index $\geq 3$. 
Hence the index of $\lda{Z}{S_1}$ at $P$ is greater than or equal to $3$.
Take the index $1$ cover $\pi$ of $\lda{Z}{S_1}$. Then the inverse image
of $E$ has more than or equal to $2$ components since $P$ is an ODP on $E$.
But this contradicts LC property of $\pb{\pi}{Z}{S_1}{E}$
after taking a small $\Bbb Q$-factorization for $\pa{\pi}{Z}{S_1}$. 
Now we finish the proof of (1).
(2) easily follows from (1) and 
the property of semi-good extraction (2) in (3.3).
\qed
\enddemo

\proclaim{Lemma (7.9) (Irreducible case)}
Consider the case $(c)$.
Let $g: Z\to X$ be a primitive log crepant divisorial extraction 
for $\ldb{X}{S}{B}$ satisfying the following conditions:   
\roster
\item $D:=\-{g}S\cap E$ is irreducible and $\Bbb P^1$.
\item $\-{g}S$ and $E$ are generically normal crossing
along $D$;
\item $\ldc{Z}{\-{g}S}{\-{g}B}{E}|_{D}=
\ldc{D}{\frac12 P_1}{\frac12 P_2}{P}$;
\item $\-{g}S$ is normal.
\endroster
Then $g$ is a good extraction (immediately by Lemma (5.2)).
Run the BTM starting by $g$ 
(we use the notation as in (2.5)).
Then surgeries while BTM are described as in the following flow chart:
\endproclaim
\newpage
$$
\hskip1in
\fbox{2.5in}{
Flip the strict transform of $\-{g}C$. \break
Consider the next extremal ray R'.}
$$
\vskip0.3in
$$
\oversetbrace\to{
\hbox{\hsize=4.6in
\fbox{1.1in}{$R'$ is a \break
divisorial ray. \break
DONE}\hskip5pt
\fbox{1.2in}{$R'$ is a flipping \break ray of type 2.\break  DONE}
\hskip5pt
\fbox{1.1in}{$R'$ is a flipping \break
ray and \break
$S' \cap \supp R' = \phi$.}
\hskip5pt
\fbox{1.2in}{$R'$ is a final \break
flipping ray. \break
DONE}
}}
$$
 
\vskip0.3in

$$
\hskip1in
\oversetbrace\to{
\hbox{\hsize=3.1in
\fbox{1.4in}{$\supp R' \cap \clc{0}= \phi$.
Then $\supp R'$ is of type 3 or 4. \break
DONE}
\hskip.1in
\fbox{1.7in}{
$\supp R' \cap \clc{0} \not = \phi$.
Then the connected component $C'$ of $\supp R'$ containing
$\clc{0}$ is of type (c)
}}}
$$
\vskip1.0in
$$\llc :=\llc(\ldc{Z'}{S'}{B'}{E'}|_{E'})$$
$$\clc{0}:=\clc{0}(\ldc{Z'}{S'}{B'}{E'}|_{E'})$$

\demo{Proof}
By the similar reason to the proof of (7.6),
first we must flip the strict transform of $C$. 
After the flips, 
we call $Z'$ the image of $Z$ and $S'$, $B'$ and $E'$ 
the strict transform of $S$, $B$ and $E$.
We can prove as Claim in the proof of (7.6)
that any component of the flipped curve 
is not contained in $S'$.
Let $R'$ be the next extremal ray.
Since $S'\cap E'$ is irreducible, we may assume that
$R'$ is a flipping ray such that $\supp R'\cap S'=\phi$, $E'.R'<0$ and 
$B'.R'>0$
(see the above flow chart).
Furthermore we may assume that there exists the unique component 
$M'$ of $\llc(\ldc{Z'}{S'}{B'}{E'}|_{E'})$ 
and there exists a point $P'\in \clc{0}(\ldc{Z'}{S'}{B'}{E'}|_{E'})$ on $M$.
(see the flow chart).
Let $C'$ be the connected component of $\supp R'$ through $P'$.
(Needless to say, 
the flip of another component of $\supp R'$ is of type 3 or 4.) 
We show that the flip of $C'$ is of type $(c)$.  
By the assumption of (5) in Set up (3.0)
and Lemma, $\supp B'|_{E'}-M'$ does not
intersect $S'$.
Hence $\supp B'|_{E'}-M'$ is contained in $\supp R'$ since $S'.R'=0$ and
$S'$ is positive for the flipped curve one before $R'$.
Furthermore 
by the existence of $P'$ and the connectedness of $\supp B'|_{E'}$,
$P'\in \supp B'|_{E'}-M'$, i.e., $\supp B'|_{E'}-M'\subset C'$. 
We show that $\supp B'|_{E'}-M'=C'$.
Assume that there is a component $C''$ of $C'$ which is not contained in  
$\supp B'|_{E'}-M'$ and contains $P'$.
After contracting $\supp B'|_{E'}-M'$, $M'$ becomes ample
since $B'.R'>0$ and $B'$ is positive for the flipped curve one before $R'$.
So by (1.9), $P'\not\in C''$, a contradiction.
Hence $\supp B'|_{E'}-M'=C'$. 
Replacing $C$ by $C'$, $S$ by $E'$, $B$ by $B'$ and $L$ by $M'$,
we can easily check the condition of case $(c)$.
(Note that $\ldb{Z'}{E'}{S'}$ is LT since LT
property of $\ldb{Z}{S}{E}$ is preserved 
by the flip of the strict transform of $C$.)
\qed
\enddemo

\proclaim{Lemma (7.10)}
Consider the case (c).
Let $g: Z\to X$ be a primitive log crepant divisorial extraction for
$\ldb{X}{S}{B}$ and $E$ the exceptional divisor.
Let $D_1$ be the component of $\-{g}S\cap E$ which intersects $\-{g}C$.
Let $C'$ be a component of $C$.
Assume that $\-{g}S \cap E$ is generically normal crossing along $D_1$
and $\-{g}C'\cap D_1$ is a smooth point on $\-{g}S$.
Then $\-{g}C'$ is a $(-1)$-curve or a $(-2)$-curve on the MRS of $\-{g}S$.
More precisely,
$\-{g}C'$ is a $(-2)$-curve on the MRS of $\-{g}S$ 
if and only if $\-{g}S$ is nonsingular on $\-{g}C'$ and 
there is no component of $\diff{\-{g}S}(0)$
and $\-{g}B|_{\-{g}S}-\-{g}C'$ intersecting $\-{g}C'$.

Assume 
in addition that $\ldc{Z}{\-{g}S}{\-{g}B}{E}|_{\-{g}S}$ is PLT at 
$\-g{C'}\cap D_1$ and $\-{g}C'$ is a $(-1)$-curve on the MRS of $\-{g}S$.
Then 
one of the following holds:
\roster
\item there is only one singularity $x\not = \-{g}C'\cap D_1$
on $\-{g}C'$ resolved by a $(-3)$-curve.
There is no component of $\diff{\-{g}S}(0)$
and $\-{g}B|_{\-{g}S}-\-{g}C'$ intersecting $\-{g}C'$; 
\item there is only one singularity $x\not = \-{g}C'\cap D_1$
which is canonical of type $A_m$
(we allow $m=0$).
There is no component of $\-{g}B|_{\-{g}S}-\-{g}C'$ intersecting $\-{g}C'$
and there is only one component of $\diff{\-{g}S}(0)$ 
(an ODP curve singularity of $Z$) or 
$\supp \-{g}B|_{\-{g}S}$ intersecting $\-{g}C'$
which pass through $x$.
\endroster
\endproclaim
\demo{proof}
By $(\ldc{Z}{\-{g}S}{\-{g}B}{E}).\-{g}C'=0$, we have
$$(\ldc{\-{g}S}{\diff{\-{g}S}(0)}{\-{g}B|_{\-{g}S}}{D_1}).\-{g}C'=0.$$
Note that $\-{g}C'$ is contained 
in $\-{g}B|_{\-{g}S}$ with coefficient $\frac12$.
Also note that $P'$ is smooth on $\-{g}S$.
By these, we obtain the equality 
$K_{\-{g}S}.\-{g}C'+\frac12 (\-{g}C')^2+1+\alpha=0$, 
where $\alpha$ is a nonnegative rational number.
On the other hand, we have the equality
$K_{\-{g}S}.\-{g}C'+(\-{g}C')^2=-2+\beta$, 
where $\beta$ is a nonnegative rational number.
By these, we obtain $K_{\-{g}S}.\-{g}C'=-2\alpha-\beta$
and $(\-{g}C')^2=-2+2\alpha+2\beta$.
Hence we obtained the first part of the lemma.

Assume in addition that $\ldc{Z}{\-{g}S}{\-{g}B}{E}|_{\-{g}S}$ is PLT at 
$\-g{C'}\cap D_1$ and $\-{g}C'$ is a $(-1)$-curve on the MRS of $\-{g}S$.
Note that after contraction of $C$, 
$\ldb{X}{S}{B}|_S$ becomes LC of index $2$ at the image of $C$. 
$\-{g}S$ becomes its MLT after contracting $\-{g}C$.
By the description of the MLT of a LC point of index $2$,
the image of $\-{g}C$ by its contraction is an ODP or 
is smooth and there is one of 
an ODP singular curve of $Z$ or $\supp \-{g}B|_{\-{g}S}$
containing it.
If the former case occurs, we have $(1)$.
If the latter case occurs, we have $(2)$.
\qed
\enddemo

\proclaim{Lemma (7.11) (Complement)}
Let $(X, \lda{X}{S})$ is a PLT $d$-fold, where $d\geq 3$ and 
$S$ is a $\Bbb Q$-Cartier
irreducible divisor.
Let $f:X\to Y$ be a projective contraction such that 
$-(\lda{X}{S})$ is $f$-nef and $f$-big.
Assume that the index of $\lda{X}{S}$ is $n$ except finite points.
Then if there is a KLT index $n$ complement $\o{B}$ for $\lda{X}{S}|_S$,
it lifts to a PLT index $n$ complement for $\lda{X}{S}$, i.e., 
there is a PLT index $n$ complement $B$ for $\lda{X}{S}$ such that
$B|_S=\o{B}$.
\endproclaim
\demo{Proof}
Let $X_0$ be the maximum locus where the index of $\lda{X}{S}$ is $n$. 
and $S_0:=S|_{X_0}$. 
Consider the local cohomology exact sequence:
$$H^1(X, \Cal O_X(-n(\lda{X}{S})-S))\to 
H^1(X_0, \Cal O_{X_0}(-n(\lda{X}{S})-S))\to 
H^2_{\{X-X_0\}}(X,\Cal O_X(-n(\lda{X}{S})-S)).$$ 
The first term vanishes by Kawamata-Viehweg vanishing theorem.
The last term vanishes by Cohen-Macaulay property of
$\Cal O_X(-n(\lda{X}{S})-S)$ and $d\geq 3$.
Hence we obtain $H^1(X_0, \Cal O_{X_0}(-n(\lda{X}{S})-S))=0$.
Next consider the exact sequence:
$$0\to \Cal O_{X_0}(-n(\lda{X}{S})-S)\to \Cal O_{X_0}(-n(\lda{X}{S}))\to
\Cal O_{S_0}(-n(\lda{X}{S}))\to 0.$$
Note that the first part of the sequence is exact since $-n(K_X+S)$
is Cartier on $X_0$.
By the above vanishing, we obtain the surjection
$$H^0(X, \Cal O_X(-n(\lda{X}{S}))\to
H^0(S, \Cal O_S(-n\lda{X}{S}).$$
Note that 
$$H^0(X_0, \Cal O_{X_0}(-n(\lda{X}{S}))\simeq H^0(X, \Cal
O_X(-n(\lda{X}{S}))$$
and 
$$H^0(S_0, \Cal O_{S_0}(-n(\lda{X}{S}))\simeq H^0(S, \Cal O_S(-n(\lda{X}{S}))$$
by reflexivility.
Hence there is a member $B'\in |-n(\lda{X}{S})|$ such that $B'|_S=n\o{B}$.
Let $B:=\frac1n B'$.
Then $\ldb{X}{S}{B}|_S=\ldb{S}{\diff{S}(0)}{\o{B}}$.
By the assumption, the right side is KLT.
Hence by the inversion of adjunction, $\ldb{X}{S}{B}$ is PLT.
This completes the proof.
\qed
\enddemo
\proclaim{Lemma (7.12) (Case $\delta=0$)}
Consider the case (c). Assume that $\delta=0$.
Then $\supp B$ and $S$ are generically simply tangent along $L$.
In particular $1 < d \leq 2$.
\endproclaim

\demo{Proof}
If $\supp B$ and $S$ are not simply tangent along $L$,
$\delta$ must be positive by the proof of Claim 1 in the proof of (6.0).
For the latter half, note that there is a log crepant divisor $F$ over $L$
for $\ldb{X}{S}{B}$ with $d(F,S)=2$.
\qed
\enddemo 
\definition{Conclusion of case (c, $\delta$)}

By running ${g''}^*(K_X +S +B)-\epsilon(\-{g''}S +dE)$-MMP and
taking the LC model as in (7.0), we obtain the primitive
extraction of $E$. We can check the above procedure is isomorphism 
near $P'$ in case (1) of (7.0)(resp. near $P_1$ and $P_2$
in case (2) of (7.0)). Hence in case (2), we can check that
$g$ is good or semi-good.

First we exclude some cases.
\enddefinition
\proclaim{Proposition (7.13)}
If the case $(1)$ of (7.0) occurs,
$P'$ is not smooth on $Z$.
\endproclaim
\demo{Proof}
Assume that $P'$ is smooth on $Z$.
Note first that all the components of $\-{g}C$ pass through $P'$.
Let $C'$ be a component of $C$. 
If there is another component of $C$,
then $\-{g}C'$ is $(-1)$-curve 
on the MRS of $\-{g}S$ by (7.10).
The same thing holds for another component of $\-{g}C$.
But this is a contradiction since $\-{g}C$ is contractible and $P'$ is smooth
on $\-{g}S$.
Hence $\-{g}C$ is irreducible.
In particular, near $P'$, 
$\ldc{Z}{\-{g}S}{\-{g}B}{E}|_{\-{g}S}=\ldb{\-{g}S}{\frac12\-{g}C}{D_1}$,
where 
$\-{g}C$ and $D_1$ are normal crossing at $P'$.
But this contradicts the assumption of $(1)$ in (7.0).
So we exclude this case.
\qed
\enddemo

\proclaim{Proposition (7.14)}
If the case $(1)$ of (7.0) occurs and
$(P'\in \-{g}S\cup E \subset Z)\simeq
(o\in (xy=0) \subset \frac12(1,0,1))$, then one of the following holds:
\roster
\item
after the flip of the strict transform of $C$, $D_1$ become a final flipping
ray or the strict transform of $E$ becomes
a exceptional divisor of a divisorial contraction.
In particular the flip of $f$ exists by (2.5);
\item
there is a index $2$ PLT complement near $C$.
In particular the flip of $f$ exists by (2.2).
\endroster

\comment
If the case $(1)$ of (7.0) occurs,
$P'$ is not the singularity of type (3) as in (7.1).
\endcomment
\endproclaim

\demo{Proof}

\proclaim{Claim}
$P'$ is a smooth point on $\-{g}S$ and an ODP on $E$.
\endproclaim

\demo{Proof}
Assume the contrary, i.e., $P'$ is an ODP on $\-{g}S$.
If $\-{g}C$ is irreducible, $\ldc{Z}{\-{g}S}{\-{g}B}{E}|_{\-{g}S}=
\ldb{\-{g}S}{D_1}{\frac12 \-{g}C}$ near $P'$.
But then $\ldc{Z}{\-{g}S}{\-{g}B}{E}|_{\-{g}S}$ must be PLT at $P'$,
a contradiction to the assumption of (1) in (7.0).
Assume that $\-{g}C$ is reducible, i.e., 
it has two components $\-{g}C'$ and $\-{g}C''$.
By (7.10), $\-{g}C'$ (resp. $\-{g}C''$) is a $(-1)$-curve or $(-2)$-curve
on the MRS of $\-{g}S$. 

If $D_1=\-{g}L$, $\-{g}S\simeq S$ and hence $Q$ is an ODP on $S$.
On the other hand, $K_S.C'<0$ (resp. $K_S.C''<0$)
by $(\lda{X}{S}).C'<0$ (resp. $(\lda{X}{S}).C'<0$).
So $C'$ (resp. $C''$) is a $(-1)$-curve 
on the MRS of singularities of $\-{g}S$. 
But on it we have the non-contractible chain consisting of 
the strict transform of $C'$, 
$(-2)$-curve resolving $Q$ and the strict transform of $C''$, 
a contradiction.

If $D_1\not=\-{g}L$, $D_1$ is not a $(-1)$-curve 
on the MRS of singularities of $\-{g}S$. 
For otherwise we have the non-contractible tree consisting of
the strict transform of $\-{g}C'$, 
$(-2)$-curve resolving $P'$ and the strict transform of $\-{g}C''$ 
and the strict transform of $D_1$, a contradiction.
In particular, $K_{\-{g}S}.D_1\geq 0$ and hence $(\lda{Z}{\-{g}S}).D_1\geq 0$.
So we can write $\pa{g}{X}{S}=\ldb{Z}{\-{g}S}{aE}$, where $a$ is a nonnegative
rational number.
Intersecting this with $\-{g}C'$ (resp. $\-{g}C''$), we obtain
$K_{\-{g}S}.\-{g}C'<0$ (resp. $K_{\-{g}S}.\-{g}C''<0$).
Hence
$\-{g}C'$ (resp. $\-{g}C''$) is a $(-1)$-curve 
on the MRS of singularities of $\-{g}S$. 
But this is a contradiction by the same reason as in case $D_1=\-{g}L$.  
\qed
\enddemo

By this Claim, 
there is an ODP curve singularity $D_0$ on $\-{g}S$ through $P'$.
Hence $\-{g}C$ must be irreducible by (1.4).
By the argument of (7.10), $\-{g}C$ is $(-1)$-curve on 
the MRS of $\-{g}S$ and there is no component of $\diff{\-{g}S}(0)$ 
intersecting $\-{g}C$ except $D_0$ and 
there is at most one singular point on $\-{g}C$ which is of type $A_m$.
(See FIGURE (XXIII).)
\comment

Hence $\ldb{Z_2}{S_2}{E_2}$ is of index $2$ at $z$.

Let $\t{Z_2}:= \pmb{Spec} \ (\Cal O_{Z_2}\oplus \Cal O_{Z_2}(-S_2))$ near $z$
and $\pi:\t{Z_2} \to Z_2$ the natural morphism.
Since $D'_0$ is a generically ODP curve, $\-{\pi}S_2$ is smooth and 
Cartier except $\-{\pi}z$. Hence by (1.11),
$\-{\pi}z$ is a smooth point of $\t{Z_2}$.
Hence locally analytically $(z\in S_2\cup E_2 \subset Z_2)\simeq
(o\in (xy=0) \subset \frac12(1,0,1))$.
But $D'_0$ is tangent to $S_2 \cap E_2$, a contradiction. 
\endcomment

Assume that $\-{g}S\cap E$ is reducible.
We prove that (1) holds in this case.

We run the BTM starting from $g$. We use the notation as in (2.5).
We see as in (7.6) that $C_1$ is
the strict transform of $C$.
We prove that $C_1^+$ is not contained in $S_2$.
Note that $h_1(C_1)$ is a smooth point of $h_1(S_1)$
by the above description of $\-{g}C$.
So $\lda{h_1(S_1)}{\frac12 h_1(D_0)}$ is canonical.
On the other hand $\ldb{S_2}{\frac12 D'_0}{B_2|_{S_2}}$ is $g_2$-ample,
where $D'_0$ is the strict transform of $D_0$.
Hence no component of $C_1^+$ is contained in $S_2$.

By this, $S_1 \dashrightarrow S_2$ is a morphism contracting $C_1$.
We denote the image point of $C_1$ by $z$.
We can easily show that $z \not \in \supp B_2$.
By this we have $B_2.D_1=0$, where $D_1$ is the strict transform of $D_1$ 
(by abuse of notation).
We can deduce from this that $D_1$ is contained in $\supp R_2$.
By the numerical properties (1) and (2),
we have $S_2.R_2 < 0$ and $E_2.R_2 >0$, or
$S_2.R_2 >0$ and $E_2.R_2 <0$.
If the latter case occurs, we are done.
Hence it suffices to exclude the former case.
Assume that the former case occurs.
As we proved above, $z$ is smooth on $S_2$.
Let $\mu:\t{S_2}\to S_2$ be the MRS of the point $D_1\cap D_2$.
Let $E_2|_{S_2}=D_1+\alpha D_2$ and 
$\mu^*(D_1+\alpha D_2)=
\-{\mu}D_1+\alpha\-{\mu}D_2+\sum\limits_{i=1}^n f_iF_i$, 
where $F_i$ is exceptional curves on $D_1\cap D_2$. 
Note that $E_2.D_1>0$.
So we obtain 
$(\-{\mu}D_1)^2+\alpha>0$ if $D_1 \cap D_2$ is smooth on $\-{g}S$ 
or
$(\-{\mu}D_1)^2+f_1>0$ if $D_1 \cap D_2$ is singular on $\-{g}S$. 
In any case 
$(\-{\mu}D_1)^2>0$ by $\alpha \leq 1$ and $f_i \leq 1$, 
a contradiction to the contractibility of $\-{\mu}D_1$.
So we are done.

Next we assume that 
$\-{g}S\cap E$ is irreducible. We have two cases, one is that $g(E)=L$
and another is that $g(E)=Q$. In these cases, we show (2) holds.

If $g(E)=L$, $\-{g}S\simeq S$.
Let $K_X +S|_S=K_S+\frac12 g(D_0)+\frac{n-1}n L$,
where $n$ is a natural number.
Since 
$\pa{g}{X}{S}.\-{g}C<0$, we have
$-1+\frac12+\frac{n-1}n<0$. 
Hence $n=1$.  
Near $Q$,
let $\t{X}:= \pmb{Spec} \ (\Cal O_X\oplus \Cal O_X(-(\lda{X}{S})))$. 
The natural morphism 
$\t{X}\to X$ is a double cover of $X$ ramified along $g(D_0)$.
By $n=1$ and (1.11), $\t{X}$ is smooth. 
Hence the index of $\lda{X}{S}$ is $2$ at $Q$ and hence on $C$.
Take a general member $\t{B}$ of $-2(\lda{X}{S})$.
Then near $C$, $\ldb{X}{S}{\frac12\t{B}}$ is PLT and of index $2$.
So the flip exists by (2.2).

Assume that $g(E)=Q$.
Let 
$\pa{g}{X}{S} |_{\-{g}S}=K_{\-{g}S}+\frac12 D_0+\frac{n-1}n\-{g}L+a'D_1$,
where $n$ is a natural number and $a'$ is a positive natural number.
Intersecting this with $\-{g}C$, we obtain 
$(\pa{g}{X}{S}).\-{g}C=a'-\frac12<0$ and so $a'<\frac12$. 
Let $\mu:(\-{g}S)^{\mu}\to \-{g}S$ be the MRS of $\-{g}S$
and $F_1,..., F_m$ exceptional curves over $D_1 \cap g^{-1} L$.
Note that $D_1 \cap g^{-1} L$ is 
the unique possibility of singularity of $\-{g}S$ on $D_1$. 
Write
$$\pc{\mu}{\-{g}S}{\frac{n-1}n \-{g}L}{\frac12 D_0}{a'D_1}=$$
$$\ldc{\-{g}S^{\mu}}{\frac{n-1}n \-{\mu}\-{g}L}{\frac12 \-{\mu}D_0}
{a'\-{\mu}D_1+\sum\limits_{i=1}^{m}f_iF_i}.$$
Intersecting this successively with $\-{\mu}D_1$, $F_1$,...and $F_m$
noting $a'<\frac12$, we know that
$\-{\mu}D_1$, $F_1$,...and $F_m$ are $(-2)$-curves, i.e., $Q$ is a 
canonical singularity of type $A_{m+1}$ on $S$ and $n=1$.
In this case, we can easily construct a index $2$ KLT complement
for $K_X +S |_S$ and hence by (7.11), 
it lifts to a index $2$ PLT complement of $K_X + S$.
This completes the proof of this proposition. 
\qed
\enddemo

Hence we need to treat the following $6$ cases:

\roster
\item"(i)" Treatment of case (1) in (7.0)-$\delta>0$ and $P'$ is a
$\frac13(2,2,1)$-singularity
\item"(ii)" Treatment of (2) in Set up (7.0)-$\delta>0$ and $g$ is good
\item"(iii)" Treatment of (2) in Set up (7.0)-$\delta>0$ and $g$ is semi-good
\item"(iv)" Treatment of (1) in Set up (7.0)-$\delta=0$ and $P'$ is a
$\frac13(2,2,1)$-singularity
\item"(v)" Treatment of (2) in Set up (7.0)-$\delta=0$ and $g$ is good
\item"(vi)" Treatment of (2) in Set up (7.0)-$\delta=0$ and $g$ is semi-good
\endroster

Our conclusion is explained in the following flow chart:
\pagebreak

\hskip1.5in
\hfil $\delta >0$ \hfil  

$$
\hskip1.5in
\fbox{2.5in}{
Choose a good primitive log crepant extraction $g$ \break
such that $d(E, S)\leq 1$ and run BTM starting by $g$}
$$
\vskip0.5in
$$
\oversetbrace\to{
\hbox{\hsize=4.4in
\fbox{1.2in}
{We fall into Case (1) of (7.0) \break
and $P'$ satisfies (3) of (7.1).
DONE}
\hskip0.5in
\fbox{1.2in}{We fall into Case (1)
 of (7.0) \break
and $P'$ satisfies (2) of (7.1)}\hskip0.5in
\fbox{1.2in}{We fall into (2) of (7.0).}
}}$$
\vskip0.5in
$$\downarrow\hskip1.0in$$
$$\fbox{2.0in}
{Let $\t{g}$ be the primitive \break
extraction of $\t{E}$. \break
Then we fall into \break
Case (2) of (7.0). \break
Run the BTM starting by $\t{g}$.} 
\hskip0.5in
\oversetbrace\to{
\hbox{\hsize=3.4in
\fbox{1.1in}{$\-{g}S \cap E$ is \break
reducible.\break
Let $\t{g}$ be the primitive extraction of $\t{E}$. \break
Then we fall into Case (2) of (7.0). \break
Run the BTM \break
starting by $\t{g}$.}
\hskip1pt
\fbox{1.1in}{$\-{g}S \cap E$ is irreducible. \break
Apply (7.9). \break
We may assume \break
that $\t{R}$ is of type (c). \break
Then this is of \break
type (c, $\delta ' < \delta$)}}
}$$
\vskip1.2in
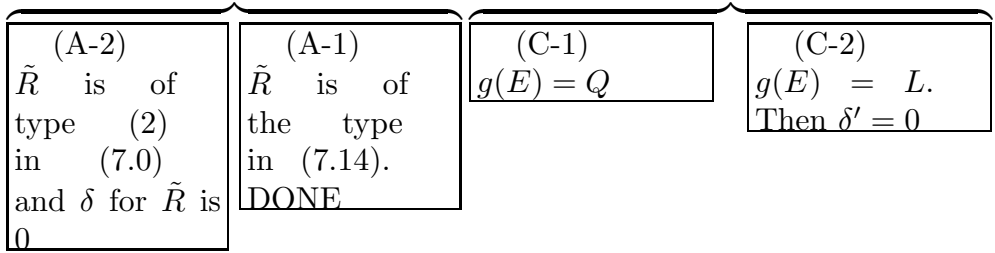
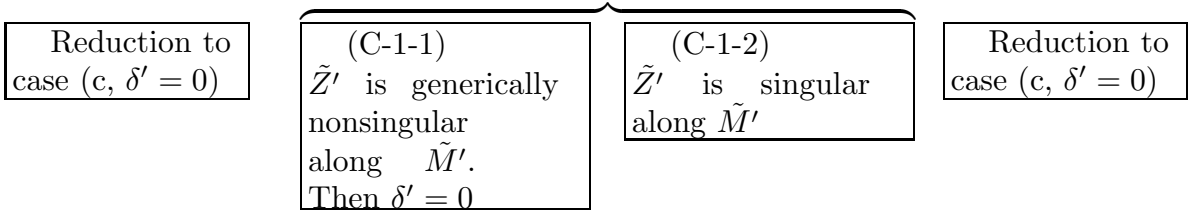
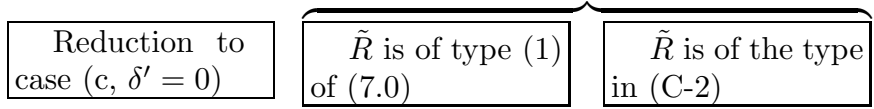
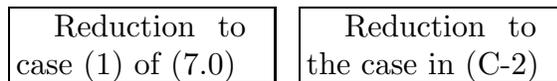
$$\hskip3.5in\Downarrow$$
\vskip5pt
$$
\oversetbrace\to{
\hbox{\hsize=1.8in
\fbox{0.9in}{(B) \break
$\-{\t{g}}S \cap \t{E}$ \break
is reducible.\break
Apply (7.6). \break
Then (1)$\sim$(4) 
\break
of (7.6) do not hold. \break
DONE}
\hskip1pt
\fbox{0.9in}{(A) \break
$\-{\t{g}}S \cap \t{E}$ \break
is irreducible.\break
Apply (7.9). \break
We may \break
assume that $\t{R}$ is of \break
type (c).
}}}
\oversetbrace\to{
\hbox{\hsize=1.8in
\fbox{0.9in}{(D)\break
 $\-{\t{g}}S \cap \t{E}$ is reducible. \break
Apply (7.6). \break
Then (1)$\sim$(4) \break
of (7.6) do not hold. \break
DONE}
\hskip1pt
\fbox{0.9in}{(C) \break
$\-{\t{g}}S \cap \t{E}$ is irreducible. \break
Apply (7.9). \break
We may \break
assume that $\t{R}$ is of \break
type (c, $\delta '$). \break
}
}}
\fbox{0.9in}{Reduction \break
to case \break
(c, $\delta ' < \delta$)}
$$
\vskip1in
$$
\hskip0.5in
\oversetbrace\to{
\hbox{\hsize=1.8in
\fbox{0.9in}{(A-2) \break
$\t{R}$ is of \break
type (2) \break
in (7.0) \break
and $\delta$ for $\t{R}$ is $0$}
\hskip1pt
\fbox{0.9in}{(A-1) \break
$\t{R}$ is of \break
the type \break
in (7.14). \break
DONE}
}}
\hskip1pt
\oversetbrace\to{
\hbox{\hsize=2.0in
\fbox{1.0in}{(C-1) \break
$g(E)=Q$}
\hskip0.1in
\fbox{1.0in}{(C-2) \break
$g(E)=L$. \break
Then $\delta ' =0$
}}}
$$
\pagebreak

$$
\hskip0.5in
\oversetbrace\to{
\hbox{\hsize=1.8in
\fbox{0.9in}{(A-2) \break
$\t{R}$ is of \break
type (2) \break
in (7.0) \break
and $\delta$ for $\t{R}$ is $0$}
\hskip0.1pt
\fbox{0.9in}{(A-1) \break
$\t{R}$ is of \break
the type \break
in (7.14). \break
DONE
}}}
\hskip1pt
\oversetbrace\to{
\hbox{\hsize=2.0in
\fbox{1.0in}{(C-1) \break
$g(E)=Q$}
\hskip0.1in
\fbox{1.0in}{(C-2) \break
$g(E)=L$. \break
Then $\delta ' =0$
}}}
$$

\vskip1.0in
$$\Downarrow\hskip3.2in\Downarrow$$

\vskip0.1in
$$
\fbox{1.0in}{Reduction to \break
case (c, $\delta ' =0$)}\hskip0.2in
\oversetbrace\to{
\hbox{\hsize=2.0in
\fbox{1.2in}{(C-1-1)
\break
$\t{Z'}$ is generically \break
nonsingular \break
along $\t{M'}$. \break
Then $\delta ' = 0$}
\hskip.1in
\fbox{1.2in}{(C-1-2)
\break
$\t{Z'}$ is singular \break
along $\t{M'}$}}}
\hskip0.1in
\fbox{1.0in}{Reduction to \break
case (c, $\delta ' =0$)}
$$
\vskip0.5in

\vskip0.5in
$$\Downarrow\hskip1.5in$$

$$
\fbox{1.1in}{Reduction to \break
case (c, $\delta ' =0$)}
\hskip 0.1in
\oversetbrace\to{
\hbox{\hsize=2.0in
\fbox{1.1in}{$\t{R}$ is of type (1) of (7.0)}
\hskip.1in
\fbox{1.1in}{$\t{R}$ is of the type in (C-2)}}}
$$
\vskip0.5in
$$\hskip1.5in\Downarrow\hskip1.0in\Downarrow$$
$$\hskip1.1in
\fbox{1.1in}{Reduction to \break
case (1) of (7.0)}
\hskip.1in
\fbox{1.1in}{Reduction to \break
the case in (C-2)}
$$

\newpage
\hfil Case (c) and $\delta =0$ \hfil 
$$
\hskip1in
\fbox{2.5in}{
Choose a good primitive log crepant extraction $g$ \break
such that $d(E, S)=d$ and run BTM starting by $g$}$$
\vskip0.5in
$$
\oversetbrace\to{
\hbox{\hsize=4.4in
\fbox{1.2in}
{We fall into Case (1) of (7.0) \break
and $P'$ satisfies (3) of (7.1).
DONE} 
\hskip5pt
\fbox{1.3in}{We fall into \break
Case (1)
 of (7.0). \break
Then $P'$ satisfies (2) of (7.1)}\hskip5pt
\fbox{1.3in}{We fall into \break
Case (2)\break
of (7.0).}
}}$$
\vskip1in
$$\downarrow\hskip1.0in$$
$$\fbox{2.0in}
{Let $\t{g}$ be the primitive \break
extraction of $\t{E}$. \break
Then we fall into \break
Case (2) of (7.0). \break
Run the BTM starting by $\t{g}$.} 
\hskip0.5in
\oversetbrace\to{
\hbox{\hsize=3.4in
\fbox{1.0in}{$\-{g}S \cap E$ is reducible.}
\hskip1pt
\fbox{1.0in}{$\-{g}S \cap E$ is irreducible.\break
Then there is an exceptional\break
complement \break
of index $2$
for \break
$\lda{X}{S}$ and \break
hence $C$
become of \break
type 4. \break
DONE}}
}$$
\vskip1.5in
$$
\oversetbrace\to{
\hbox{\hsize=1.8in
\fbox{1.0in}{(B) $\-{\t{g}}S \cap \t{E}$ is reducible.
\break
Apply (7.6). \break
Then (1)$\sim$(4) of (7.6) \break
do not hold. \break
DONE}
\hskip1pt
\fbox{1.0in}{(A) $\-{\t{g}}S \cap \t{E}$ is irreducible.
\break
Apply (7.9). \break
We may assume that $\t{R}$ is of type (c).
}}}
\oversetbrace\to{
\hbox{\hsize=1.8in
\fbox{1in}{$g(E)=Q$.\break
Then there is a PLT \break
complement of \break
index $2$ for $\lda{X}{S}$ and  \break
hence $C$ \break
become of type 3  \break
DONE}
\hskip1pt
\fbox{1in}{$g(E)=L$.\break
Then there is a PLT \break
complement of \break
index $2$
for $\lda{X}{S}$ and  hence $C$ \break
become of type 3  \break
DONE
}
}}
$$
\pagebreak
$$
\hskip0.5in
\oversetbrace\to{
\hbox{\hsize=1.8in
\fbox{0.9in}{(A-1) $\t{R}$ is of \break
the type \break
in (7.14). \break
DONE}
\hskip1pt
\fbox{0.9in}{(A-2) 
$\t{R}$ is of \break
type (2) \break
in (7.0) \break
and $\delta$ for $\t{R}$ is $0$}}}
$$
\vskip1.0in
$$\hskip1.0in\Downarrow$$
\vskip0.5in
$$\hskip2.5in
\fbox{1.0in}{Reduction to \break
case (c, $\delta ' =0$)\break
and the case of \break
(2) in (7.0)}
\hskip1.0in$$

\newpage

\comment

(i) $\Rightarrow$ (v), (vi) or Proposition (7.14). 
\newline
More precisely, after replacing $g$ by the primitive extraction of the
exceptional divisor of the weighted blow up as in Theorem (7.1) (2) and
running BTM, we show that the middle flip satisfies (v), (iv) or satisfies
the assumption of Proposition (7.14).
\newline
(ii) $\Rightarrow$ (ii) with smaller $\delta$ 
\newline
More precisely, by running BTM for $g$, we show that the middle flip
satisfies (ii) with smaller $\delta$.
\newline
(iii) $\Rightarrow$ (iv), (v) or (vi)

\newline
(iv) $\Rightarrow$ (v) or (vi) 
\newline
More precisely, we show that there exists a good or semi-good extraction
and replacing $g$ by this and running BTM, we show that the middle flip
satisfies (iv), (v) or (vi)?
\newline
(v) $\Rightarrow$ index 2 exceptional complement 
\newline
More precisely, we show that there exists index 2 exceptional complement near $C$.
\newline
(vi) $\Rightarrow$ index 2 PLT complement
\newline
More precisely, we show that there exists index 2 PLT complement near $C$.
\newpage
\endcomment

\definition{Treatment of (1) in Set up (7.0)-$P'$ is 
a $\frac13(2,2,1)$-singularity}
We treat the case $\delta >0$ and $\delta =0$ simultaneously. 
Recall there is a exceptional divisor $F'$ over $Z$ 
as in Theorem (7.1) $(2)$.
We change the notation. We write $\t{E}$ instead of $F'$.
Note that $d(\t{E}, S) >1$.
Running 
$\pb{(g\circ h)}{X}{S}{B}-\epsilon(\-{(g\circ h)}S+\t{d}\t{E})$-MMP over $X$,
we obtain a primitive extraction of $\t{E}$, which we call $\t{g}:\t{Z}\to X$.
Let $\t{D}:= \-{\t{g}}S\cap \t{E}$ if $\-{\t{g}}S\cap \t{E}$ is irreducible
or $\t{D}_1\cup \t{D}_2:=\-{\t{g}}S\cap \t{E}$ 
if $\-{\t{g}}S\cap \t{E}$ is reducible, 
where only $\t{D}_1$ intersects $\-{\t{g}}C$.
(By (5.2), there are such two possibilities.)
Note that outside $E$, the MMP is isomorphism. 
Hence if $\-{\t{g}}S\cap \t{E}$ is irreducible (resp. reducible), 
$\t{D}$ (resp. $\t{D_1}$) is 
the strict transform of $\t{E}\cap\-{(g\circ h)}S$.
Furthermore by the description of (7.1), (7.0) (2) holds for $\t{g}$.
By Theorem (7.1) (2), the assumptions of (7.9)
(resp. (7.6)) holds.
Hence $\t{g}$ is good (resp. semi-good).
\enddefinition
\definition{(A)}
If $\t{g}$ is good, 
by applying (7.9),
we have only consider the case $(c)$ appears while BTM.
We use the notation as in the proof of (7.9)
(replacing $g$, $E$, $D_1$, $C'$, $M$, $M'$,...etc 
by $\t{g}$, $\t{E}$, $\t{D_1}$, $\t{C'}$, $\t{M}$, $\t{M'}$,..etc).
We claim that $\t{C'}$ satisfies the assumption of Proposition (7.14)
(and there exists the flip of $\t{C}'$, hence of $C$) or
$\t{C'}$ is of type (2) in (7.0), $\delta=0$.
Do the same procedure as in (7.0) for $\t{C}'$.
Recall 
$\ldc{\t{Z}}{\-{\t{g}}S}{\-{\t{g}}B}{\t{E}}|_{\t{E}}=
\ldc{\t{E}}{\t{M_1}}{\t{M}}{\frac12\t{M_3}+\frac12\t{M_4}}$,
where $\t{M_1}:=\-{\t{g}}S|_{\t{E}}$, 
$\t{M_3}:=\supp \-{\t{g}}B|_{\t{E}} -\t{M}$ and 
$\t{M_4}$ is a generically ODP curve of $\t{Z}$. 
Since $\t{g}$ is good, 
$\t{M_1}\cap \t{M_3}\not= \t{M_1}\cap \t{M_4}$.
See FIGURE (XXIV).
\enddefinition
\definition{(A-1)}
If we fall into case (1) of (7.0), 
the flip of $\t{C}'$ satisfies the assumption of Proposition (7.14)
by the existence of the generically ODP curve along $\t{M_4}$.
\enddefinition
\definition{(A-2)} 
If we fall into case (2) of (7.0),
assume that there exists a log crepant divisor $F$ for
$\ldc{\t{Z}'}{\t{S}'}{\t{B}'}{\t{E}'}$ whose center contains $\t{P}'$ 
(i.e., is equal to $\t{P}'$ or $\t{M'}$) and $d(\t{E}', F)\leq 1$ 
(to get a contradiction).
Let $l:W\to \t{Z'}$ be a (primitive) log crepant extraction of $F$.
Here we can take $F$ the divisor corresponding to $E$ in (7.0).
Since the flips before $\t{C'}$ do not touch with $\t{M}$,
the center of $F$ on $\t{Z}$ is $\t{P}$ or $\t{M}$ and we can consider
$l$ is defined over $\t{Z}$.
Run 
the $\pb{(\t{g}\circ l)}{X}{S}{B}-\epsilon (\-{(\t{g}\circ l)}S+d(F, S)F)$-MMP
over $X$.
We show possibly after the flip of $\-{l}\t{M}$ that 
$\-{l}\t{E}$ is contracted to a curve.

Let $A$ be the first extremal ray.
By the assumption that we fall into case (2) of (7.0), 
$\-{l}\t{M_3}\cap F|_{\-l{\t{E}}}\not= \-{l}\t{M_4}\cap F|_{\-l{\t{E}}}$.
Hence $\-{(\t{g} \circ l)}B.\-{l}\t{M_4}=0$. 
Since for a $l$-exceptional curve, $\-{(\t{g}\circ l)}B$ is positive,  
we know that $\-{l}\t{B}.A\leq 0$.

If $\-{(\t{g} \circ l)}B.A=0$, 
$\-{l}\t{M_4}\subset \supp A$
since the extremal ray for $F$ is positive for $\-{(\t{g}\circ l)}B$ and we
play $2$ ray game.
By (1.6), $\-l\t{M_4}$ cannot be birationally contracted.
Hence $A$ is a divisorial ray and $\-{l}\t{E}$ is contracted to a curve.

If $\-{(\t{g}\circ l)}B.A<0$, 
$A$ is a flipping ray and $\supp A\subset \-{(\t{g}\circ l)}B$.
But $\-{l}\t{M_3}$ cannot be birationally contracted by (1.6)
since $\-{l}\t{M_3}$ intersects
$F|_{\-l{\t{E}}}$ and $\-{(\t{g}\circ l)}S|_{\-l{\t{E}}}$ at two points.
So $\supp A=\-{l}\t{M_2}$ and $l(F)=\t{P}$.
Since $A$ is
$\pb{(\t{g}\circ l)}{X}{S}{B}-\epsilon(\-{(\t{g}\circ l)}S+d(F, S)F)$-extremal ray,
$\-l{\t{E}}.A<0$ and $F.A>0$. So this is a flip of type 1 and the flip exists
by (2.0).
Let $A'$ be the next extremal ray.
We express all the strict transforms with superscript $+$.
Since $\-{(\t{g}\circ l)}B^+.\-{l}\t{M_4^+}=0$ and for the flipped curve, 
$\-{(\t{g}\circ l)}B^+$ is positive, 
we know that $\-{(\t{g}\circ l)}B^+.A'\leq 0$.
But we see that
$\-{(\t{g}\circ l)}B^+.A'=0$ and $\-{l}\t{E}^+$ is contracted to a curve.
Hence we obtain what we want.

Furthermore by the above argument, we also see that
when $\t{E}$ is contracted to a curve, $\-{l}\t{E}$ and $F$ is generically
normal crossing along 
the irreducible component of $\-{l}\t{E} \cap F$ intersecting $\-{l}\t{C'}$.
The same thing also holds for $\-{l}\t{E}$ and $\-{(\t{g}\circ l)}S$.
Since 
$\-{l}\t{M_4}^+.\-{(\t{g}\circ l)}B^+=0$ and 
it intersects $F$ and $\-{(\t{g}\circ l)}S$ at distinct two points, 
$\-{l}\t{M_4}^+$ is an irreducible fiber by (1.6).
Furthermore it intersects $F$ and $\-{(\t{g}\circ l)}S$ at smooth points
normally.
Hence
$\-{(\t{g}\circ l)}S.\-{l}\t{M_4}=F.\-{l}\t{M_4}=1$.
Since $\-{l}\t{M_4}$ does not intersect singular curve of $\t{Z}$ and neither
does a general fiber,
$n:=-\-{l}\t{E}.\-{l}\t{M_4}$ is a positive integer.
Since $(\-{(\t{g}\circ l)}S+d(\t{E}, S)\-{l}\t{E}+d(F, S)F).\-{l}\t{M_4}=0$ and
$d(F, S)=d(F, \t{E})d(\t{E}, S)$, we obtain the equality
$1-nd(\t{E}, S)+d(F, \t{E})d(\t{E}, S)=0$.
From now on, we treat separately the case $\delta>0$ and the case $\delta=0$.

First we treat the case $\delta>0$.

If $n\geq 2$ in the above equality,
then $$d(\t{E}, S)\leq \frac12(d(\t{E}, S)d(F, \t{E})+1)\leq 
\frac12(d(\t{E}, S)+1)$$ by the assumption $d(F, \t{E})\leq 1$.
So $d(\t{E},S)\leq 1$, 
a contradiction to the choice of $\t{E}$.

If $n=1$, $d(F, S)=d(\t{E}, S)-1\leq \frac 43 -1<1$.
Hence after contracting $\t{E}$, 
we fall into case (2) of (7.0) with $\delta>0$ since
after contracting $\-{l}\t{E}$, $F$ and $\-{(\t{g}\circ l)}S$ 
cross normally along
the image of $\-{l}\t{E}$
by $\-{l}\t{E}.\-{l}\t{M_4}= -1$,  a contradiction
by Claim in (7.0).

Next we treat the case $\delta=0$.

If $n\geq 2$ in the above equality,
we get a contradiction by the same way as the case $\delta>0$.

If $n=1$,
then $d(F, S)=d(\t{E}, S)-1\leq 2-1=1$, 
a contradiction to the assumption $\delta=0$.

Now we finish the treatment of the case $\t{g}$ is good.
\enddefinition
\definition{(B)}
If $\t{g}$ is semi-good, we can apply (7.6)
and (7.7) for $\t{g}$.
Note that $d(\t{E}, S)>1$ so if (4) of Lemma (7.7)
holds, then for all log crepant divisors their multiplicities in $S$
is greater than 1 (resp. d) in case $\delta >0$ (resp. $\delta =0$), 
a contradiction.
Hence in this case, we can construct the flip of $f$ as 
explained in the flow chart.
\qed
\enddefinition

\definition{Treatment of (2) in Set up (7.0), $\delta>0$ and $g$ is good}

In this case, we can apply (7.9).
We have only to consider the case $(c)$ appears while BTM.
We use the notation as in the proof of (7.9).
Then we can check the flip of $C'$ is of type $(c, \delta'<\delta)$.
So we reduce the existence of the flip of $C$ to the flip of type 
$(c, \delta'<\delta)$.
\qed
\enddefinition

\definition{Treatment of (2) in Set up (7.0), $\delta>0$ and $g$ is semi-good}
(See FIGURE (XXV).)
Let $h:W\to Z$ be the simple blow up of $D_1$ and
$\t{E}$ the exceptional divisor of $h$.
Running 
$\pb{(g\circ h)}{X}{S}{B}-\epsilon(\-{(g\circ h)}S+\t{d}\t{E})$-MMP over $X$,
we obtain an primitive extraction of $\t{E}$, which we call $\t{g}:\t{Z}\to X$.
Let $\t{D}:= \-{\t{g}}S\cap \t{E}$ 
if $\-{\t{g}}S\cap \t{E}$ is irreducible
or $\t{D}_1\cup \t{D}_2:=\-{\t{g}}S\cap \t{E}$ 
if $\-{\t{g}}S\cap \t{E}$ is reducible, 
where only $\t{D}_1$ intersects $\-{\t{g}}C$.
(By (5.2), there are such two possibilities.)
Since $h|_{\-{(g\circ h)}S}$ is isomorphism over $D_1$ and
the MMP is isomorphism outside $E$,
the conditions $(2)$ and $(3)$ of a semi-good extraction in (3.3)
are preserved (replacing $E$ by $\t{E}$). 
Hence for $\t{g}$, the assumptions of (7.9)
or (7.6) are satisfied.
\enddefinition
\definition{(C)}
If $\-{\t{g}}S\cap \t{E}$ is irreducible, 
we can apply (7.9) for $\t{g}$ and we may assume 
that a flip of type $(c)$ appears while BTM.
We use the notation as in the lemma 
(replacing $g$, $E$, $D_1$, $C'$, $M'$,..etc 
by $\t{g}$, $\t{E}$, $\t{D_1}$, $\t{C'}$, $\t{M'}$,..etc).
We prove that $\delta$ for the flip of $\t{C}'$ is $0$. 

First consider the case $g(E)=Q$.
\enddefinition
\definition{(C-1-1)}
Assume that $\t{Z'}$ is nonsingular along $\t{M'}$. Then the following claim
holds:
\enddefinition

\proclaim{Claim}
Assume that $\t{Z}'$ is nonsingular generically along $\t{M}'$.
Then there is 
a log crepant divisor $F$ for $\ldb{X}{S}{B}$ over 
$\t{Q}':= \t{M}'\cap (\t{S}'|_{\t{E}'})$ whose multiplicity in $\t{S}'+\t{E}'$
is less than or equal to $1$.
\endproclaim
\demo{Proof}
Note that $S_1\to \t{S'}$ is the contraction of $D_2$.
By (1) and (2) of Lemma (7.7),
 $S_1$ is smooth along $D_2$ and 
by (3) of the same Lemma, $d(E, S)=\frac3{-(D_2)^2_{S_1}-1}\leq 1$.
Let $n:=-(D_2)^2_{S_1}$.
Then by these, $\t{Q}$ is $\frac1n(1,1)$ singularity on $\t{S}$ and $n\geq 4$.
It is easy to see that the index of $\lda{\t{Z}}{\t{S}}$ at $\t{Q}$ is $n$ and
the index $1$ cover for this is smooth by (1.11).
Hence locally analytically $(\t{Q}\in \t{S}\cup\t{E}\subset \t{Z})\simeq
(o\in (xy=0) \subset \frac1n(k,1,1))$, 
where $1\leq k<n$ is a natural number.
The weighted blow up with the weight $\frac1n(k,1,1)$ extracts a log crepant
divisor for $\ldc{\t{Z}}{\t{S}}{\t{E}}{\t{B}}$ since the restriction
of the weighted blow up to $\t{S}$ coincides with the MRS
and it must be log crepant for
$\ldc{\t{Z}}{\t{S}}{\t{E}}{\t{B}}|_{\t{S}}$ by $n\geq 4$.
The multiplicity of the exceptional divisor in $\t{S}+\t{E}$ is
equal to $\frac{k+1}n$ and this is less than or equal to $1$.
So we are done.
\qed
\enddemo
By this claim we complete the reduction 
if $\t{Z}'$ is nonsingular generically along $\t{M}'$.
In fact, let $d(F, \t{S}')=a$ and $d(F, \t{E}')=b$.
Since $a+b\leq 1$, $a<1$ and $b<1$, 
$$d(F, S)=d(F,\t{S}'+d(\t{E}', S)\t{E}') 
=a+d(\t{E}', S)b=a+(d(E, S)+1)b<1+d(E,S).$$
Hence by (4) of Lemma (7.7) for $g$,
it must be $E$.
Furthermore we know by this 
that on $\t{Z}$, the center of the exceptional divisor of the simple blow up
of $D_2$ is also $\t{Q}'$.
Hence for any log crepant divisor for
$\ldc{\t{Z}}{\t{S}'}{\t{B}'}{\t{E}'}$ whose center contains $\t{P}'$, 
its multiplicity in $\t{S}'+(d(E, S)+1)\t{E}'$ is greater than $d(E, S)+1$, 
which in turn show that 
its multiplicity in $\t{E}'$ is greater than $1$, which we want.

\definition{(C-1-2)}
If $\t{Z}'$ is singular along $\t{M}'$,
then there exists a log crepant divisor for 
$\ldc{\t{Z'}}{\t{S'}}{\t{E'}}{\t{B'}}$ over $\t{M}'$ such that its multiplicity
in $\t{E'}$ is less than $1$ (take a general hyperplane section at a general
point of $\t{M}$ and take the MRS).
By (4) of Lemma (7.7), it must be $E$.
Furthermore the exceptional divisor of simple blow up of $D_2$
lands on $\t{Q'}$.
Since $E$ is the unique log crepant divisor for
$\ldb{\t{Z}}{\t{E}}{\t{B}}$ such that its multiplicity in $\t{E}$ is
less than $1$,
the flip of $\t{C}'$ is the flip of the type 
(1) in (7.0) or the type treated next in (C-2).
(Note that $\t{M'}$ corresponds to $L$ in (3.2).)
So we obtain the reduction as explained in the flow chart.
\enddefinition

\definition{(C-2)}
Next consider the case $g(E)=L$.
By (2) of Lemma (7.8),
only $E$ and the exceptional divisor of the simple blow up of $D_1$ or $D_2$ 
are the divisors such that its multiplicity in $S$ is less than
or equal to $d(E, S)+1$. Note that they cannot land on $\t{P}'$ or
$\t{M}'$ since they are not projective.
Hence for any log crepant divisor for
$\ldc{\t{Z}}{\t{S}'}{\t{B}'}{\t{E}'}$ whose center contains $\t{P}'$, 
its multiplicity in $\t{S}'+(d(E, S)+1)\t{E}'$ is greater than $d(E, S)+1$, 
which in turn show that 
its multiplicity in $\t{E}'$ is greater than $1$, which we want.
\enddefinition

\definition{(D)}
If $\-{\t{g}}S\cap \t{E}$ is reducible, 
we can apply (7.6) for $\t{g}$ and assume that
(1), (2), (3) and (4) of Lemma (7.7) hold. 
By (4), 
$d(\t{F}, S)=d(\t{E}, S)+1$ and $a_l(\t{F}, \ldb{X}{S}{B})=0$ 
for the exceptional divisor $\t{F}$ of the simple blow up 
along $\t{D}_1$ or $\t{D}_2$
and for any other log crepant divisor $\t{F}\not= \t{E}$ 
for $\ldb{X}{S}{B}$, $d(\t{F}, S)>d(\t{E}, S)+1$.
In particular,
for any log crepant divisor $\t{F}\not= \t{E}$ 
for $\ldb{X}{S}{B}$, $d(\t{F}, S)>1$.
But this contradicts the existence of $E$.
Hence 
in this case we can construct the flip of $f$ 
as explained in the flow chart.
\qed
\enddefinition

\definition{Treatment of (2) in Set up (7.0), $\delta=0$ and $g$ is good}

By (7.9), we can assume that there is a curve $M$ of 
$\llc(\ldc{Z}{\-{g}S}{\-{g}B}{E}|_E)$ through $Q'$.
\enddefinition
\proclaim{Claim}
$Q'$ is smooth on $\-g{S}$, $E$ and $Z$.
\endproclaim

\demo{Proof}
By (1.12), there is no curve singularity
outside $\-{g}S\cup E$.
If there is a curve singularity through $Q'$ contained in $\-{g}S\cup E$,
it must be $M$ 
(note that $Z$ is smooth generically along $L'$ by (7.12)).
Hence there is a log crepant divisor $F$ for $\ldc{Z}{\-{g}S}{\-{g}B}{E}$
such that $d(F, \-{g}S+E)\leq 1$ over that curve singularity.
Hence $$d(F, S)=d(F, \-{g}S)+d(E, S)d(F, E)< d(d(F, \-{g}S)+d(F,E))
\leq d,$$ a contradiction to the choice of $d$.
By the same way, we can prove that there is no log crepant divisor
for $\ldc{Z}{\-{g}S}{\-{g}B}{E}$ such that $d(F, \-{g}S+E)\leq 1$ over
$Q'$ or a curve through $Q'$.
So $\ldc{Z}{\-{g}S}{\-{g}B}{E}$ satisfies (*) at $Q'$.
Then (1) or (2) of Theorem (7.1) holds 
since both $\-{g}S$ and $E$ contain no curve singularity through $Q'$.
If (2) holds, let $h:Z'\to Z$ be the extraction of $F'$ taken in
(7.1).
The existence of $\-{g}L$ and $M$ contradicts the description of (2).
\comment, on $?$, there are two point contained
in $\clc{0}(\pc{h}{Z}{\-{g}S}{\-{g}B}{E}|_F)$, a contradiction.
Hence (1) holds and this gives the claim.
\endcomment
\qed
\enddemo

By the above claim, there is a index $2$ exceptional complement for
$\lda{X}{S}|_S$.
In fact, we have only to replace $\-{g}L$ by two smooth curve intersecting
$D$ at smooth points transversely with coefficient $\frac12$ and
consider the images of them.
By the argument of (7.11), 
it lifts to a index $2$ exceptional complement
for $\lda{X}{S}$. Hence the flip of $C$ exists (type 4).
\qed

\definition{Treatment of (2) in Set up (7.0), $\delta=0$ and $g$ is semi-good}

First we treat the case $g(E)=Q$.
By (7.7), 
we may assume that (1), (2), (3) and (4) hold.
On the other hand, $2\geq d(E, S)>1$ by (7.12),
Hence $-(D_2)^2_{\-{g}S}= 3$ and $d(E, S)=\frac32$. 
We prove that there is a PLT index $2$ complement near $C$, which in turn
show that the flip of $C$ exists (type 3).
\enddefinition

\proclaim{Claim 1}
\roster
\item
Let $\mu:E^{\mu}\to E$ be the MRS of $E$.
Then $(\-{\mu}D_1)^2\leq 0$.
\item 
One of $P_1$ and $P_2$, say $P_1$ is smooth on $E$ and $\-{g}S$ and 
$\-{g}C$ passes through it. 
$P_2$ is an ODP on $E$ and is smooth on $\-{g}S$. 
$(D_1)^2_E=\frac12$ and $(D_1)^2_{\-{g}S}=-2$.
\item 
$C$ is irreducible.
On $\-{g}C$, $\-{g}S$ has one singular point resolved by one $(-3)$-curve.
\item $E\simeq \Bbb F_{2,0}$. The index of $\lda{X}{S}$ is $2$ at $Q$.
\endroster
See FIGURE (XXVI).
\endproclaim

\demo{Proof}
\roster
\item
Assume that $(\-{\mu}D_1)^2>0$.
Since $\-{\mu}D_1$ and $\-{\mu}D_2$ are nef and big,
the inequality 
$(\-{\mu}D_1)^2.(\-{\mu}D_2)^2\leq (\-{\mu}D_1.\-{\mu}D_2)^2$ holds.
On the other hand,
$(\-{\mu}D_2)^2=2$ and
$\-{\mu}D_1.\-{\mu}D_2=1$
since $E$ is smooth on $D_2$, $(D_2)^2_E=2$ and $D_1.D_2=1$. 
But this contradicts the inequality.
\item
By the ampleness of $\-{g}S|_E$, we have 
$$0<\-{g}S.D_1=((D_1+D_2).D_1)_E=(D_1)^2_{E}+1, i.e.,$$ 
$$(D_1)^2_{E}>-1. \tag *$$
On the other hand,
by 
$0=g^*S.D_1=\-{g}S.D_1+\frac32E.D_1$,
we obtain
$$(D_1)^2_{\-{g}S}=-\frac23(D_1)^2_E-\frac53. \tag **$$ 

We may assume that an irreducible component of
$\-{g}C$ passes through $P_1$. 
Then $P_1$ is smooth 
since $\ldc{Z}{\-{g}S}{\-{g}B}{E}|_{\-{g}S}$ is LT at $P_1$
of index $2$ and $\-{g}C\subset \supp \-{g}B|_{\-{g}S}$.
$P_1$ is smooth also on $E$ by the same reason.
$P_2$ is a smooth point or an ODP on $\-{g}S$ and $E$.
By (*) and (1),
$(D_1)^2_E=0$ (resp. $(D_1)^2_E=-\frac12 \ \text{or} \ \frac12$) 
if $P_2$ is smooth on $E$ 
(resp. if $P_2$ is an ODP on $E$).
Hence by (**), 
$(D_1)^2_{\-{g}S}=-\frac53$ 
(resp. $(D_1)^2_{\-{g}S}=-\frac43\ \text{or} \ -2$)
if $P_2$ is smooth on $E$ 
(resp. if $P_2$ is an ODP on $E$). 
If $P_2$ is an ODP on $\-{g}S$, $(D_1)^2_{\-{g}S}\in \frac{\Bbb Z}2 \ 
\text{and} \ \not\in \Bbb Z$, a contradiction.
So $P_2$ is smooth on $\-{g}S$ and hence $(D_1)^2_{\-{g}S}$ is an integer.
Consequently $P_2$ must be an ODP on $E$
and $(D_1)^2_{\-{g}S}=-2$.
\item
By (2), there is a ODP curve singularity $D_0$ on $\-{g}S$ through $P_2$.
So no irreducible component of $\-{g}C$ passes through $P_2$.  
This prove that $C$ is irreducible.
By (2), 
$(\lda{Z}{\-{g}S}).D_1=(\lda{\-{g}S}{\frac12 D_0}).D_1
=\frac12$. 
On the other hand, by (2), $E.D_1=-1$.
Hence $$\lda{Z}{\-{g}S}+\frac12E=\pa{g}{X}{S}. \tag ***$$
By (7.10), $\-{g}C$ is a $(-1)$-curve or $(-2)$-curve on the
MRS of $\-{g}S$.
But if $\-{g}C$ is a $(-2)$-curve, by $(***)$,
$(\lda{X}{S}).C=\frac12>0$, a contradiction to the speciality of the flipping
contraction $f$.
Hence $\-{g}C$ is a $(-1)$-curve on the MRS of $\-{g}S$ 
and $(1)$ or $(2)$ of Lemma (7.10) holds.
If $(2)$ holds and $m=0$,
by $(***)$, 
we have $(\lda{X}{S}).C=-1+\frac12+\frac12=0$, 
a contradiction to the speciality of the flipping contraction $f$.
If $(2)$ holds and $m>0$,
we obtain on the MRS of $\-{g}S$
a noncontractible chain consisting of
a $(-2)$-curve resolving of $x$, the strict transforms of $\-{g}C$
and $D_1$, a contradiction.
Hence $(1)$ holds and we are done.
\item 
Let $\nu:E\to \o{E}$ be the morphism defined by sufficient multiple of $D_2$.
Since $E$ is smooth on $D_2$ and $(D_2)^2_E=2$, 
$\o{E}\simeq \Bbb F_{2,0} \ \text{or} \ \Bbb P^1 \times \Bbb P^1$.
Since $(D_1)^2_E=\frac12$, $\nu(D_1)^2\geq\frac12$. 
But on $\o{E}$, $\nu(D_1)$ and $\nu(D_2)$ intersect at one point.
Hence $\nu(D_1)^2=\frac12$, which in turn show that 
$\t{E} \simeq \Bbb F_{2,0}$ and $\nu$ is isomorphism
near $D_1$. But since $D_1+D_2$ is ample on $E$, the exceptional curve 
of $\nu$ must intersect $D_1$. So actually $\nu$ is isomorphism, i.e., 
$E\simeq \Bbb F_{2,0}$.
There is no curve singularity of $Z$ on $E$.
In fact, if there exists a curve singularity of $Z$ on $E$,
it intersects $D_1$ or $D_2$ by the ampleness of $D_1+D_2$ on $E$.
But this is impossible by (2) and Lemma (7.7) (1).
Hence $E$ is a Cartier divisor since $E\simeq \Bbb F_{2,0}$
(cf. proof of (1.11)).
By this and $(***)$, 
the index of $\lda{X}{S}$ is $2$ at $Q$.
\endroster
\qed
\enddemo
By Claim 1 (2) and (4), $(\lda{X}{S}).C=-\frac16$.
There exists a Weil divisor $\o{B}$ on $S$
through $x$ (the singularity on $C$) such that
$\lda{S}{\frac12\o{B'}}$ is PLT and has index $2$ and $\o{B'}.C=-\frac13$.
By (7.11), 
it lifts to an index $2$ PLT complement for $\lda{X}{S}$.

Next we treat the case $g(E)=L$.
By (7.8),
(1) and (2) of the lemma hold.
\proclaim{Claim 2}
\roster
\item 
One of $P_1$ and $P_2$, say $P_1$ is smooth on $E$ and $\-{g}S$ and 
$\-{g}C$ passes through it. 
$P_2$ is smooth on $E$ and is an ODP on $\-{g}S$. 
$(D_1)^2_E=0$ and $(D_1)^2_{\-{g}S}=-\frac32$.
$P$ is smooth on $E$.
\item 
$C$ is irreducible.
On $\-{g}C$, 
$\-{g}S$ is smooth and 
there is a ODP curve singularity $D_{-1}\subset \-{g}S$ of $Z$
or 
a component of $\supp \-{g}B$ with coefficient $\frac 12$ in $\-{g}B|_{\-{g}S}$
normally crossing with $\-{g}C$.
\item The index of $\lda{X}{S}$ is $2$ near $C$.
\endroster
See FIGURE (XXVII).
\endproclaim
\demo{Proof}
\roster
\item
Take $P_1$ as in Claim 1. The statement about $P_1$ is proved as in Claim 1.
Let $l$ be a general fiber of $E\to L$.
Then $\-{g}B.l=\frac12$ by (7.12).
On the other hand $\-{g}B.D_1\geq \frac12$ 
since $D_1$ intersect $\-{g}B|_E$ at a smooth point $P_1$.
Since the fiber containing $D_1$ deforms to $l$,
$D_1$ must be numerically equivalent to $l$ in $Z$.
In particular, $D_1$ is an irreducible and reduced fiber.
Hence $P$ and $P_2$ are smooth on $E$.
Furthermore $-\frac12=E.l=E.D_1=(D_1)^2_{\-{g}S}+1$, i.e., 
$(D_1)^2_{\-{g}S}=-\frac32$. 
By this, $P_2$ is an ODP on $\-{g}S$.
\item
The proof is similar to Claim 1 (3).
By (1), $P_2$ is an ODP on $\-{g}S$.
So no irreducible component of $\-{g}C$ passes through $P_2$.  
This prove that $C$ is irreducible.
By $(\lda{Z}{\-{g}S}).D_1=(K_{\-{g}S}).D_1=0$, 
we have $$(\lda{Z}{\-{g}S})=\pa{g}{X}{S}.\tag *$$
By (7.10), $\-{g}C$ is a $(-1)$-curve or $(-2)$-curve on the
MRS of $\-{g}S$.
If $\-{g}C$ is a $(-2)$-curve, by $(*)$, $(\lda{X}{S}).C=0$, a contradiction. 
Hence $\-{g}C$ is a $(-1)$-curve and $(1)$ or $(2)$ of the lemma holds.
If $(1)$ holds,
then we obtain a non-contractible chain on the MRS of $\-{g}S$
consisting of the $(-3)$-curve resolving the singularity on $\-{g}C$,
the strict transforms of $\-{g}C$ and $D_1$ and $(-2)$-curve resolving
$P_2$, a contradiction. 
Hence $(2)$ holds.
But if $m>0$, we have again a noncontractible chain, a contradiction.
Hence we obtain what we want.
\item
By $(2)$, the index of $\lda{Z}{\-{g}S}$ near $\-{g}C$ is $2$. 
Hence by $(*)$, the index of $\lda{X}{S}$ is also $2$ near $C$.
\endroster
\qed
\enddemo
Let $\o{B}$ is a general smooth Cartier divisor near $C$ such that
$\o{B}.C=1$.
By Claim 2 (3), $(\lda{X}{S}).C=-\frac12$.
Hence $(\ldb{X}{S}{\frac12\o{B}}).C=0$.
By Claim 2 (3), the index of $\ldb{X}{S}{\frac12\t{B}}$ is $2$.
Taking $\t{B}$ generally, this is PLT so we are done.
\qed

\enddocument

\endRefs

\enddocument